\documentclass[12pt]{article}
\usepackage{pgf,tikz}
\usetikzlibrary{arrows}
\usepackage{amsmath}%
\usepackage{amsfonts}%
\usepackage{amssymb}%
\usepackage{graphicx}%
\def\binom#1#2{{#1\choose#2}}

\newtheorem{lemma}{Lemma}

\begin{document}
\title{All Ramsey $(C_n,K_6)$ critical  graphs for large $n$}
\author{C. J. Jayawardene \\
Department of Mathematics\\
University of Colombo \\
Sri Lanka\\
email: c\_jayawardene@maths.cmb.ac.lk\\
\\
W. C. W. Navaratna  and J.N. Senadheera \\
Department of Mathematics \\
The Open University of Sri Lanka \\
Sri Lanka\\
email: wcper@ou.ac.lk, jnsen@ou.ac.lk \\
}

\maketitle
\begin{abstract} Let $G$ and $H$ be finite graphs without loops or multiple edges. If for any two-coloring of the edges of a complete graph $K_n$, there is a copy of $G$ in the first color, red, or a copy of $H$ in the second color, blue, we will say $K_n\rightarrow (G,H)$. The Ramsey number $r(G, H)$ is defined as the smallest positive integer $n$ such that $K_{n} \rightarrow (G, H)$. A two-coloring of $K_{r(G, H)-1}$ such that $K_{r(G, H)-1} \not \rightarrow (G,H)$ is called a critical coloring. A  Ramsey critical $r(G, H)$ graph is a graph induced by the first color of a critical coloring. In this paper, when $n \geq 15$, we show that there exist exactly  sixty eight non-isomorphic Ramsey critical $r(C_n, K_6)$ graphs.  
\end{abstract}


\section{Introduction}
\noindent For any two graphs $G$ and $H$, the \textit{Ramsey number} $r(G,H)$ is the smallest positive integer $n$ such that $K_n\rightarrow (G,H)$.   The \textit{classical Ramsey number} $r(s, t)$ is defined as $r(K_s,K_t)$ and the \textit{diagonal Ramsey number} is defined as  $r(K_n,K_n)$.  These numbers have been studied extensively in the last five decades. The difficulty in exact determination of \textit{diagonal Ramsey number} $r(n, n)$, swifts expeditiously from the apparent $r(3,3)=6$ to the unmanageable $r(5,5)$. Currently, the best known lower and upper bounds for $r(5,5)$ are  43 and 48 (see \cite{Ra}). Many interesting variations of the basic problem of finding  classical Ramsey numbers have emerged. One such variation  is the calculation of the number of \textit{Ramsey critical  $(G,H)$ graphs}, for any pair of graphs $(G,H)$. In this paper,  we show that there are exactly 68  Ramsey critical $r(C_n, K_6)$   graphs, for all $n$ exceeding fourteen.

\vspace{14pt}

\section{Notation}
\noindent The \textit{complete graph} on $n$ vertices, a \textit{cycle} on  $n$ vertices and a \textit{Star} on $n+1$ (see \cite{JaSa}) vertices are denoted by $K_n$, $C_n$  and $K_{1,n}$ respectively. Given a graph $G$ and a vertex $v \in V(G)$, we define the \textit{neighbourhood of $v$ in $G$},  $\Gamma(v)$, as the set of vertices adjacent to $v$ in $G$. The \textit{degree} of a vertex $v$, $d(v)$, is defined as the cardinality of  $ $ $\Gamma(v)$, $ $ i.e. $d(v)=| \Gamma(v)|$. The minimum degree of a graph $G(V,E)$ denoted by $\delta(G)$ is defined as  $ $ $ \min\{ d(v) | v\in V  \}$. 

\vspace{8pt}

\noindent Given a graph $G$, we say $I \subseteq V(G)$ is an \textit{independent set}, if no pair  of  vertices of $I$ is adjacent to each other in $G$. Equivalently, $I$ forms a clique in $G^c$.  Given a graph $G=(V,E)$, we define \textit{the independence number}, $\alpha(G)$, as the size of the largest independent set. Thus, $\alpha(G)=\max \{ \left| I\right| : I$ is an independent set of $G \}$. In the special case of $H=K_m$, alternatively $r(G,K_m)$ can be viewed as the smallest positive integer $n$ such that every graph  of order $n$ either contains $G$ as a subgraph or else satisfies $\alpha(G) \geq m$. For a non-empty subset $S$ of $V$, \textit{the induced subgraph of $S$ in $G$} denoted by $G[S]$ is defined as the  subgraph obtained by deleting all the vertices of $S^c$ from $G$. For two disjoint subgraphs $H$ and $K$ of $G$, we denote the set of edges between $H$ and $K$ by $E(H,K)$.

\vspace{10pt}

\section{Lemmas used to generate  Ramsey critical $(C_n,K_6)$  graphs for $n \geq 15$}

In an attempt to prove Bondy and Erd{\"o}s conjecture  $r(C_n,K_m) = (n-1)(m-1)+1$, for all $(n,m)\neq(3,3)$ satisfying $n \geq m \geq 3$ under certain restrictions, Schiermeyer has proved that  $r(C_n,K_6) = 5(n - 1)+1$, for $n \geq 6$  (see \cite{Ra,Sc}). Characterizing all  Ramsey critical $(C_n,K_6)$ graphs boils down to  finding all (red/blue) colorings of $K_{r(C_n,K_6)-1}$ such that there is no red $C_n$ or a blue $K_6$. This is achieved by finding all $C_n$-free graphs on $K_{r(C_n,K_6)-1}$ vertices such that $\alpha(G)<6$. We first  prove that any $C_n$ -free graph (where $n \geq 15$) of order $5(n-1)$ with $\alpha(G) \leq 5$  contains a $5K_{n-1}$. To prove this, we use seven lemmas of which the first three are already proven results. For ease of reference, we reiterate Lemma \ref{l1a} from \cite{JaRo1}, Lemma  \ref{l1b} from \cite{Ja} and Lemma  \ref{l2}  by Bollob{\'a}s et al from \cite{BoJaSh}.
\vspace{10pt}

\begin{lemma}
\label{l1a}
\noindent (\cite{JaRo1}, Lemma 2; \cite{JaRo1}). A $C_n$- free graph $G$ of order $N$  with independent number  less than or equal to $m$ has minimal degree  greater than or equal to $N-r(C_n,K_{m})$.
\end{lemma}

\vspace{2pt}

\begin{lemma}
\label{l1b}
\noindent (\cite{Ja}, Lemma 8).  A $C_n$- free graph (where $n \geq 7$) of order $4(n-1)$ with no independent set of 5 vertices  contains a  $4K_{n-1}$.
\end{lemma}

\vspace{2pt}

\begin{lemma} (\cite{BoJaSh}, Lemma 5)
\label{l2}
\noindent Suppose $G$ contains the cycle $(u_1, u_2, ... , u_{n-1}, u_1)$ of length $n - 1$ but
no cycle of length $n$. Let $Y = V(G) \setminus \{u_1, u_2, ... , u_{n-1}\}$. Then,

\vspace{8pt}
\noindent \textbf{(a)} No vertex $x \in Y$ is adjacent to two consecutive vertices on the cycle.

\vspace{8pt}
\noindent \textbf{(b)} If $x \in Y$ is adjacent to $u_i$ and $u_j$ then $u_{i+1}u_{j+1} \notin E(G)$.

\vspace{8pt}
\noindent \textbf{(c)} If $x \in Y$ is adjacent to $u_i$ and $u_j$ then no vertex $x' \in Y$ is adjacent to both
$u_{i+1}$ and $u_{j+2}$.

\vspace{8pt}
\noindent \textbf{(d)} Suppose $\alpha(G) = m - 1$ where $m \leq  \frac {n + 2}{2}$ and $\{x_1, x_2, ... ,x_{m -1} \} \subseteq Y$ is an $(m - 1)$-element independent set. Then, no member of this set is adjacent to $m -2$ or more vertices on the cycle (We have taken the liberty of making a slight correction to the inequality $m \leq \frac {n + 2}{2}$ of the original \cite{BoJaSh}, Lemma 5(d)).
\end{lemma}

\vspace{12pt}

\noindent The next lemma plays a pivotal role in proving the main results of this paper. 

\vspace{4pt}
\begin{lemma}
\label{l3}
\noindent A $C_n$ -free graph (where $n \geq 15$) of order $5(n-1)$ with no independent set of 6 vertices contains a $5K_{n-1}$.
\end{lemma}

\noindent {\bf Proof.} We shall assume that in each of the three cases  $n = 15$, $n = 16$ and $n \geq 17$ we consider,  $G$ as a  graph on $5(n-1)$ vertices  satisfying $C_n \not \subseteq G$ and $\alpha(G) \leq 5$. Since  $r(C_{n-1}, K_6)=5n-9 \leq 5(n-1)$ (see \cite{BoJaSh,Ra}), there exists a cycle $C=(u_1, u_2, ... , u_{n-1}, u_1)$ of length $n - 1$ in $G$. In consistent with the notation of \cite{BoJaSh}, define $H$  as the induced subgraph of $G$ not containing the vertices of the cycle $C$. Then, $|V(C)|=n-1$ and $|V(H)|=4(n-1)$. 

\vspace{6pt}

\noindent Suppose  there exists an independent set $Y=\{y_1, y_2,y_3, y_4, y_5 \}$ of size 5 in $H$, so that $\alpha(G)=5$. From  Lemma \ref{l2} (as $5\leq \frac{n+2}{2}$), it follows that no vertex of $Y$ is adjacent to four or more vertices of the $C_{n-1}$. Thus, $\left| E(Y,V(C)) \right| \le 15$.  For ease of reference, we define such a graph structure as a  \textbf{Standard Configuration ($n$)}.

\vspace{10pt}

\noindent \textbf{Case 1: $n \geq 17$}

\vspace{6pt}
\noindent  Now, $\left| E(Y,V(C)) \right| \le 15 < n-1$. Thus, there exists a vertex $x\in V(C)$ adjacent to no vertex of $Y$. This gives, an independent set $Y \cup \{x\}$ of size 6, a contradiction.

\vspace{10pt}
\noindent \textbf{Case 2: $n = 16$}

\vspace{6pt}  
\noindent In this case as $n-1=15$, in order to avoid an independent set of size 6, each vertex of $V(C)$ must be adjacent to at least one vertex of $Y$. Thus, we get that for each  $1\leq i\leq 5$,  $\left|\Gamma (y_i) \cap V(C) \right|=3$ and for each $1\leq j  < j' \leq 5$, $\Gamma (y_j) \cap \Gamma (y_{j'}) \cap V(C)=\phi$.

\vspace{10pt}

\noindent By Lemma \ref{l1a}, as $\delta(G) \geq 14$,  $\left|\Gamma (y_i) \cap V(H \setminus Y)\right| \geq 11$ for  $i=1,2$. Since $r(P_3,K_6)=11$ and $\alpha(G) < 6$, each of $G[\Gamma (y_1) \cap V(H \setminus Y)]$ and  $G[\Gamma (y_2) \cap V(H \setminus Y)]$ contains a copy of $P_3$. Thus, $P_3 \subseteq \Gamma (y_1) \cap V(H \setminus Y)$, where the $P_3$ is induced by  $\{x,y,z\}$ such that $(x,y), (y,z) \in E(G)$ and $P_3 \subseteq \Gamma (y_2) \cap V(H \setminus Y)$,  where this $P_3$ is induced by  $\{p,q,r\}$ such that $(p,q), (q,r) \in E(G)$.

\vspace{10pt}

\noindent Suppose that $x$ is not adjacent to any vertex of $\{y_2,y_3,y_4,y_5\}$ and $p$ is not adjacent to any vertex of $\{y_1,y_3,y_4,y_5\}$. Re-order the vertices of the cycle such that $y_1\in Y$ is adjacent to $u_1$. In this ordering, let $y_1$ be also adjacent to $u_t$ where $2 \leq t \leq 15$.

\begin{center}
\definecolor{uuuuuu}{rgb}{0.26666666666666666,0.26666666666666666,0.26666666666666666}
\definecolor{ffffff}{rgb}{1.0,1.0,1.0}
\begin{tikzpicture}[line cap=round,line join=round,>=triangle 45,x=1.3921568627450978cm,y=1.2600000000000002cm]
\clip(-5.699999999999999,5.0799999999999965) rectangle (4.5000000000000115,10.07999999999999);
\fill[color=ffffff,fill=ffffff,fill opacity=0.1] (-3.02,5.62) -- (-2.44,5.58) -- (-1.8938741688442602,5.77936543467826) -- (-1.4760526241370269,6.183624209200794) -- (-1.2587805067477513,6.722876308876984) -- (-1.279626139628259,7.303880146721288) -- (-1.5349851234768812,7.826174880916263) -- (-1.980703570122651,8.199451007020894) -- (-2.5397127109155475,8.359165691724547) -- (-3.1153547869738647,8.2777028150796) -- (-3.608096053803033,7.969148048524964) -- (-3.9327370699544106,7.4868533143299905) -- (-4.033144454464381,6.9142117535470975) -- (-3.891956858378718,6.350238294320912) -- (-3.533586899709677,5.892449071289668) -- cycle;
\fill[color=ffffff,fill=ffffff,fill opacity=0.1] (1.08,6.36) -- (2.14,6.34) -- (2.486579144363348,7.341939567385363) -- (1.6407768353717538,7.9811722747028835) -- (0.7714631162884599,7.374300247160362) -- cycle;
\draw (-4.8199999999999985,1.820000000000001) node[anchor=north west] {$v_{1,3}$};
\draw (-0.6399999999999939,1.6400000000000012) node[anchor=north west] {$v_{ 2,3}$};
\draw (6.8212102632969694E-15,0.24000000000000307) node[anchor=north west] {$v_{2,1}$};
\draw (0.3800000000000072,7.279999999999993) node[anchor=north west] {$y_{1}$};
\draw (1.6000000000000085,8.579999999999991) node[anchor=north west] {$y_{2}$};
\draw (2.5800000000000094,7.739999999999993) node[anchor=north west] {$y_{3}$};
\draw (2.200000000000009,6.559999999999994) node[anchor=north west] {$y_{4}$};
\draw (-2.599999999999996,8.93999999999999) node[anchor=north west] {$u_{6}$};
\draw (6.8212102632969694E-15,0.24000000000000307) node[anchor=north west] {$u_{4}$};
\draw (-1.4599999999999946,6.279999999999995) node[anchor=north west] {$u_{1}$};
\draw (-1.739999999999995,8.619999999999992) node[anchor=north west] {$u_t$};
\draw [color=ffffff] (-3.02,5.62)-- (-2.44,5.58);
\draw [color=ffffff] (-2.44,5.58)-- (-1.8938741688442602,5.77936543467826);
\draw [color=ffffff] (-1.8938741688442602,5.77936543467826)-- (-1.4760526241370269,6.183624209200794);
\draw [color=ffffff] (-1.4760526241370269,6.183624209200794)-- (-1.2587805067477513,6.722876308876984);
\draw [color=ffffff] (-1.2587805067477513,6.722876308876984)-- (-1.279626139628259,7.303880146721288);
\draw [color=ffffff] (-1.279626139628259,7.303880146721288)-- (-1.5349851234768812,7.826174880916263);
\draw [color=ffffff] (-1.5349851234768812,7.826174880916263)-- (-1.980703570122651,8.199451007020894);
\draw [color=ffffff] (-1.980703570122651,8.199451007020894)-- (-2.5397127109155475,8.359165691724547);
\draw [color=ffffff] (-2.5397127109155475,8.359165691724547)-- (-3.1153547869738647,8.2777028150796);
\draw [color=ffffff] (-3.1153547869738647,8.2777028150796)-- (-3.608096053803033,7.969148048524964);
\draw [color=ffffff] (-3.608096053803033,7.969148048524964)-- (-3.9327370699544106,7.4868533143299905);
\draw [color=ffffff] (-3.9327370699544106,7.4868533143299905)-- (-4.033144454464381,6.9142117535470975);
\draw [color=ffffff] (-4.033144454464381,6.9142117535470975)-- (-3.891956858378718,6.350238294320912);
\draw [color=ffffff] (-3.891956858378718,6.350238294320912)-- (-3.533586899709677,5.892449071289668);
\draw [color=ffffff] (-3.533586899709677,5.892449071289668)-- (-3.02,5.62);
\draw [color=ffffff] (1.08,6.36)-- (2.14,6.34);
\draw [color=ffffff] (2.14,6.34)-- (2.486579144363348,7.341939567385363);
\draw [color=ffffff] (2.486579144363348,7.341939567385363)-- (1.6407768353717538,7.9811722747028835);
\draw [color=ffffff] (1.6407768353717538,7.9811722747028835)-- (0.7714631162884599,7.374300247160362);
\draw [color=ffffff] (0.7714631162884599,7.374300247160362)-- (1.08,6.36);
\draw (6.8212102632969694E-15,0.24000000000000307) node[anchor=north west] {$y_{2}$};
\draw (0.6200000000000074,6.499999999999995) node[anchor=north west] {$y_{5}$};
\draw (-1.4760526241370269,6.183624209200794)-- (-1.2587805067477513,6.722876308876984);
\draw (-1.2587805067477513,6.722876308876984)-- (-1.279626139628259,7.303880146721288);
\draw (-1.279626139628259,7.303880146721288)-- (-1.5349851234768812,7.826174880916263);
\draw (-1.5349851234768812,7.826174880916263)-- (-1.980703570122651,8.199451007020894);
\draw (-1.980703570122651,8.199451007020894)-- (-2.5397127109155475,8.359165691724547);
\draw (-2.5397127109155475,8.359165691724547)-- (-3.1153547869738647,8.2777028150796);
\draw (-3.1153547869738647,8.2777028150796)-- (-3.608096053803033,7.969148048524964);
\draw (-3.608096053803033,7.969148048524964)-- (-3.9327370699544106,7.4868533143299905);
\draw (-3.9327370699544106,7.4868533143299905)-- (-4.033144454464381,6.9142117535470975);
\draw (-4.033144454464381,6.9142117535470975)-- (-3.891956858378718,6.350238294320912);
\draw (-3.891956858378718,6.350238294320912)-- (-3.533586899709677,5.892449071289668);
\draw (-3.533586899709677,5.892449071289668)-- (-3.02,5.62);
\draw (-3.02,5.62)-- (-2.44,5.58);
\draw (-2.44,5.58)-- (-1.8938741688442602,5.77936543467826);
\draw (-1.8938741688442602,5.77936543467826)-- (-1.4760526241370269,6.183624209200794);
\draw (0.7714631162884599,7.374300247160362)-- (-1.4760526241370269,6.183624209200794);
\draw [->] (0.7714631162884599,7.374300247160362) -- (-0.9800000000000001,7.96);
\draw (-4.039999999999997,5.859999999999995) node[anchor=north west] {$u_{12}$};
\draw (-4.039999999999997,8.519999999999992) node[anchor=north west] {$u_{8}$};
\draw (-4.719999999999998,7.259999999999994) node[anchor=north west] {$u_{10}$};
\draw (-2.399999999999996,5.499999999999996) node[anchor=north west] {$u_{14}$};
\draw (-0.1,9.12)-- (0.4,9.62);
\draw (0.4,9.62)-- (0.9,9.14);
\draw (0.42000000000000726,10.11999999999999) node[anchor=north west] {$y$};
\draw (-0.5199999999999937,9.59999999999999) node[anchor=north west] {$x$};
\draw (0.8800000000000078,9.61999999999999) node[anchor=north west] {$z$};
\draw (0.7714631162884599,7.374300247160362)-- (-0.1,9.12);
\draw (0.7714631162884599,7.374300247160362)-- (0.4,9.62);
\draw (0.7714631162884599,7.374300247160362)-- (0.9,9.14);
\draw (-2.9999999999999964,7.199999999999994) node[anchor=north west] {$C$};
\draw (1.3600000000000083,7.399999999999993) node[anchor=north west] {$Y$};
\draw [->] (-0.1,9.12) -- (-1.24,8.14);
\begin{scriptsize}
\draw [fill=black] (1.08,6.36) circle (1.5pt);
\draw [fill=black] (-3.02,5.62) circle (1.5pt);
\draw [fill=black] (-2.44,5.58) circle (1.5pt);
\draw [fill=uuuuuu] (-1.8938741688442602,5.77936543467826) circle (1.5pt);
\draw [fill=uuuuuu] (-1.4760526241370269,6.183624209200794) circle (1.5pt);
\draw [fill=uuuuuu] (-1.2587805067477513,6.722876308876984) circle (1.5pt);
\draw [fill=uuuuuu] (-1.279626139628259,7.303880146721288) circle (1.5pt);
\draw [fill=uuuuuu] (-1.5349851234768812,7.826174880916263) circle (1.5pt);
\draw [fill=uuuuuu] (-1.980703570122651,8.199451007020894) circle (1.5pt);
\draw [fill=uuuuuu] (-2.5397127109155475,8.359165691724547) circle (1.5pt);
\draw [fill=uuuuuu] (-3.1153547869738647,8.2777028150796) circle (1.5pt);
\draw [fill=uuuuuu] (-3.608096053803033,7.969148048524964) circle (1.5pt);
\draw [fill=uuuuuu] (-3.9327370699544106,7.4868533143299905) circle (1.5pt);
\draw [fill=uuuuuu] (-4.033144454464381,6.9142117535470975) circle (1.5pt);
\draw [fill=uuuuuu] (-3.891956858378718,6.350238294320912) circle (1.5pt);
\draw [fill=uuuuuu] (-3.533586899709677,5.892449071289668) circle (1.5pt);
\draw [fill=black] (2.14,6.34) circle (1.5pt);
\draw [fill=uuuuuu] (2.486579144363348,7.341939567385363) circle (1.5pt);
\draw [fill=uuuuuu] (1.6407768353717538,7.9811722747028835) circle (1.5pt);
\draw [fill=black] (0.7714631162884599,7.374300247160362) circle (1.5pt);
\draw [fill=black] (-0.1,9.12) circle (1.5pt);
\draw [fill=black] (0.9,9.14) circle (1.5pt);
\draw [fill=black] (0.4,9.62) circle (1.5pt);
\end{scriptsize}
\end{tikzpicture} 
\end{center}
\begin{center}
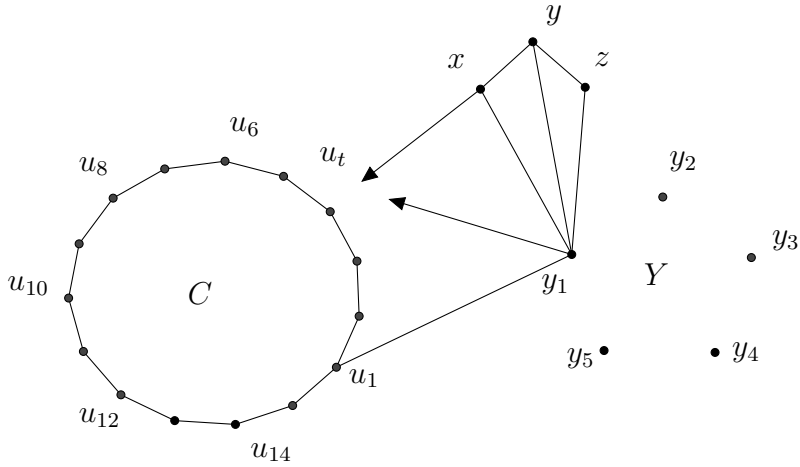
Figure 1. Configuration for $n=16$
\end{center}

\noindent By Lemma \ref{l2}(a), $t \neq 2$. In order to avoid an independent set of size 6, induced by $\{x,u_t,y_2,y_3,y_4,y_5\}$, we get that $(x,u_t) \in E(G)$. However, $t \neq 3$, in order to avoid a  $C_{16}$ comprising $(u_1,y_1,x, u_3,...,u_{15},u_1)$. Also, $t \neq 4$ in order to avoid a $C_{16}$ comprising  $(u_1,y_1,y,x,u_4,...,u_{15},u_1)$  and $t \neq 5$ in order to avoid a $C_{16}$ comprising  $(u_1,y_1,z,y,x, u_5,...,u_{15},u_1)$.

\vspace{8pt}
\noindent Thus, any pair of vertices adjacent to $y_1$ in $C$ cannot be separated by a path of length 1, 2, 3 or 4  along $C$. Thus, $\Gamma (y_1) \cap C =\{u_1,u_6,u_{11}\}$. In this scenario, we use the prerogative that $(y_2,u_2)\in E(G)$. Then,  by the previous argument $\Gamma (y_2) \cap C =\{u_2,u_7,u_{12}\}$. But by Lemma \ref{l2}(b),  $(u_2,u_7)\notin E(G)$. Henceforth, we will get that $\{u_2,u_7,y_1,y_3,y_4,y_6\}$  is an independent set of size 6, a contradiction. 
\vspace{10pt}

\noindent This implies that there is a vertex of $X=\{x,y,z\}$ adjacent to some vertex of $\{y_2,y_3,y_4,y_5\}$ or there is a vertex of $\{p,q,r\}$ adjacent to some vertex of $\{y_1,y_3,y_4,y_5\}$. Therefore, without loss of generality, we may assume that $y_1$ is adjacent to $X=\{x,y,z\} \subseteq V(H \setminus Y)$ and $y_1$ is adjacent to $X'=\{x',y',z'\} \subseteq V(C)$ where $X'$ induces a  $P_3$  and $y_2$ is adjacent to $x$. Next since $C$ has 15 points without loss of generality, $\{x',y',z'\}=\{u_1,u_3,u_5\}$ or $\{x',y',z'\}=\{u_1,u_3,u_6\}$ or $\{x',y',z'\}=\{u_1,u_3,u_7\}$  or  $\{x',y',z'\}=\{u_1,u_3,u_8\}$ or  $\{x',y',z'\}=\{u_1,u_3,u_9\}$ or $\{x',y',z'\}=\{u_1,u_4,u_7\}$ or $\{x',y',z'\}=\{u_1,u_4,u_8\}$ or $\{x',y',z'\}=\{u_1,u_4,u_9\}$ or $\{x',y',z'\}=\{u_1,u_4,u_{10}\}$ or $\{x',y',z'\}=\{u_1,u_5,u_9\}$ or $\{x',y',z'\}=\{u_1,u_5,u_{10}\}$ or $\{x',y',z'\}=\{u_1,u_6,u_{11}\}$. Moreover, as $y_1$ and $y_2$ are connected by paths of lengths 2, 3 and 4 in $H$, no pair of vertices selected from  $\Gamma (y_1) \cap V(C)$ and  $\Gamma (y_2) \cap V(C)$ can be separated by a path of length 3, 4 or 5 along the cycle $C$. Using this  we argue that when $\{x',y',z'\}=\{u_1,u_3,u_5\}$, $\Gamma (y_2) \cap V(C) = \phi$,  when $\{x',y',z'\}=\{u_1,u_3,u_6\}$, $\Gamma (y_2) \cap V(C) = \phi$ and when  $\{x',y',z'\}=\{u_1,u_3,u_7\}$, $\Gamma (y_2) \cap V(C) \subseteq \{u_9\}$.  Similarly, when $\{x',y',z'\}=\{u_1,u_3,u_8\}$, $\Gamma (y_2) \cap V(C) \subseteq \{u_2,u_9,u_{10}\}$, when $\{x',y',z'\}=\{u_1,u_3,u_9\}$, $\Gamma (y_2) \cap V(C) \subseteq \{u_2,u_{10}\}$, when $\{x',y',z'\}=\{u_1,u_4,u_7\}$, $\Gamma (y_2) \cap V(C) = \phi$ and when $\{x',y',z'\}=\{u_1,u_4,u_8\}$, $\Gamma (y_2) \cap V(C) \subseteq \{u_2,u_{10}\}$. When $\{x',y',z'\}=\{u_1,u_4,u_9\}$, $\Gamma (y_2) \cap V(C) \subseteq \{u_2,u_3,u_{10}\}$, when $\{x',y',z'\}=\{u_1,u_4,u_{10}\}$, $\Gamma (y_2) \cap V(C) \subseteq \{u_2,u_3\}$,  when $\{x',y',z'\}=\{u_1,u_5,u_9\}$, $\Gamma (y_2) \cap V(C) \subseteq \{u_3,u_7\}$ and when $\{x',y',z'\}=\{u_1,u_5,u_{10}\}$, $\Gamma (y_2) \cap V(C) \subseteq \{u_{3}\}$ and  when $\{x',y',z'\}=\{u_1,u_6,u_{11}\}$, $\Gamma (y_2) \cap V(C) =\phi$. Since none of these give a viable configuration, we get a contradiction.  

\vspace{8pt}
\noindent \textbf{Case 3: $n = 15$}

\vspace{8pt}

\noindent To deal with the case $n = 15$, we first prove three Lemmas. Lemma $\ref{l31}$, deals with the possible scenarios generated by the Standard Configuration ($15$). Lemmas $\ref{l32}$ and $\ref{l33}$ deal with showing that none of the scenarios generated by Lemma $\ref{l31}$ give viable configurations.

\vspace{4pt}  

\begin{lemma}
\label{l31}
\noindent In the Standard Configuration $(n=15)$, one of the following three scenarios $(a)$, $(b)$ and $(c)$ will occur:

\vspace{8pt}
\noindent \textbf{(a)} $y_1 \in Y$ is a vertex of the subgraph $K_4$ (see Figure 2(a)) in $H$.\vspace{8pt}

\vspace{8pt}
\noindent \textbf{(b)} $y_1,y_2\in Y$ are vertices of the subgraph $K$ (see Figure 2(b1)) or subgraph $K'$ (see Figure 2(b2)) in $H$.

\vspace{8pt}
\noindent \textbf{(c)} $y_1,y_2\in Y$ are vertices of the subgraph $L$ (see Figure 2(c)) in $H$.
\end{lemma}

\begin{center}
\begin{tikzpicture}[line cap=round,line join=round,>=triangle 45,x=1.0cm,y=1.4062499999999998cm]
\clip(-3.7600000000000025,4.999999999999994) rectangle (11.060000000000006,8.199999999999996);
\draw (-4.820000000000003,1.8199999999999918) node[anchor=north west] {$v_{1,3}$};
\draw (-0.6400000000000009,1.6399999999999917) node[anchor=north west] {$v_{ 2,3}$};
\draw (-5.684341886080805E-16,0.23999999999999103) node[anchor=north west] {$v_{2,1}$};
\draw (-5.684341886080805E-16,0.23999999999999103) node[anchor=north west] {$u_{4}$};
\draw (-5.684341886080805E-16,0.23999999999999103) node[anchor=north west] {$y_{2}$};
\draw (7.98,7.56)-- (8.4,7.52);
\draw (7.720000000000004,7.979999999999995) node[anchor=north west] {$x$};
\draw (7.98,7.56)-- (8.5,6.4);
\draw (8.5,6.4)-- (8.4,7.52);
\draw (8.180000000000005,6.3399999999999945) node[anchor=north west] {$y_1$};
\draw (0.54,7.6)-- (1.42,7.6);
\draw (1.42,7.6)-- (2.08,7.6);
\draw (0.23999999999999957,7.999999999999996) node[anchor=north west] {$x$};
\draw (0.54,7.6)-- (1.46,6.54);
\draw (1.46,6.54)-- (1.42,7.6);
\draw (1.08,6.459999999999995) node[anchor=north west] {$y_1$};
\draw (2.98,7.6)-- (2.08,7.6);
\draw (2.08,7.6)-- (2.12,6.52);
\draw (2.12,6.52)-- (2.98,7.6);
\draw (1.9800000000000004,6.439999999999994) node[anchor=north west] {$y_2$};
\draw (-3.7800000000000025,5.879999999999994) node[anchor=north west] {$(a): graph$};
\draw (0.23999999999999957,5.839999999999994) node[anchor=north west] {$(b1): graph$};
\draw (4.52,7.58)-- (5.58,7.58);
\draw (5.58,7.58)-- (6.62,7.56);
\draw (4.300000000000002,7.979999999999995) node[anchor=north west] {$x$};
\draw (4.52,7.58)-- (5.26,6.54);
\draw (5.26,6.54)-- (5.58,7.58);
\draw (4.880000000000002,6.459999999999995) node[anchor=north west] {$y_1$};
\draw (6.62,7.56)-- (5.92,6.52);
\draw (5.780000000000003,6.459999999999995) node[anchor=north west] {$y_2$};
\draw (4.2200000000000015,5.839999999999994) node[anchor=north west] {$(b2):  graph$};
\draw (5.92,6.52)-- (5.58,7.58);
\draw (1.9800000000000004,8.019999999999996) node[anchor=north west] {$z$};
\draw (2.760000000000001,8.019999999999996) node[anchor=north west] {$w$};
\draw (5.220000000000002,7.979999999999995) node[anchor=north west] {$y$};
\draw (6.520000000000003,7.9599999999999955) node[anchor=north west] {$z$};
\draw (1.06,7.999999999999996) node[anchor=north west] {$y$};
\draw (8.5,6.4)-- (8.92,7.56);
\draw (8.92,7.56)-- (9.3,7.56);
\draw (9.3,7.56)-- (8.5,6.4);
\draw (8.5,6.4)-- (9.68,7.54);
\draw (9.68,7.54)-- (10.22,7.54);
\draw (10.22,7.54)-- (8.5,6.4);
\draw (7.98,7.56)-- (9.1,6.44);
\draw (8.4,7.52)-- (9.1,6.44);
\draw (9.1,6.44)-- (8.92,7.56);
\draw (9.1,6.44)-- (9.3,7.56);
\draw (9.1,6.44)-- (9.68,7.54);
\draw (9.1,6.44)-- (10.22,7.54);
\draw (7.98,7.56)-- (9.74,6.36);
\draw (9.74,6.36)-- (8.4,7.52);
\draw (9.74,6.36)-- (8.92,7.56);
\draw (9.3,7.56)-- (9.74,6.36);
\draw (9.74,6.36)-- (9.68,7.54);
\draw (10.22,7.54)-- (9.74,6.36);
\draw (-3.3,7.6)-- (-2.24,7.6);
\draw (-2.24,7.6)-- (-1.18,7.6);
\draw (-3.5400000000000027,8.019999999999996) node[anchor=north west] {$x$};
\draw (-3.3,7.6)-- (-2.24,6.52);
\draw (-2.24,6.52)-- (-2.24,7.6);
\draw (-2.720000000000002,6.499999999999995) node[anchor=north west] {$y_1$};
\draw (-2.580000000000002,8.019999999999996) node[anchor=north west] {$y$};
\draw (-1.2800000000000014,8.019999999999996) node[anchor=north west] {$z$};
\draw (-1.18,7.6)-- (-2.24,6.52);
\draw (8.920000000000005,6.319999999999994) node[anchor=north west] {$y_2$};
\draw (9.680000000000005,6.2999999999999945) node[anchor=north west] {$y_3$};
\draw (8.300000000000004,7.9599999999999955) node[anchor=north west] {$y$};
\draw (8.780000000000005,7.9599999999999955) node[anchor=north west] {$z$};
\draw (9.220000000000004,7.9599999999999955) node[anchor=north west] {$u$};
\draw (9.640000000000006,7.9599999999999955) node[anchor=north west] {$v$};
\draw (10.140000000000006,7.939999999999995) node[anchor=north west] {$w$};
\draw (7.960000000000004,5.859999999999994) node[anchor=north west] {$(c): graph$};
\draw (-1.5400000000000014,5.839999999999994) node[anchor=north west] {$K_4$};
\draw (2.640000000000001,5.799999999999994) node[anchor=north west] {$K$};
\draw (6.620000000000004,5.839999999999994) node[anchor=north west] {$K'$};
\draw (10.200000000000006,5.819999999999994) node[anchor=north west] {$L$};
\begin{scriptsize}
\draw [fill=black] (7.98,7.56) circle (1.5pt);
\draw [fill=black] (8.4,7.52) circle (1.5pt);
\draw [fill=black] (8.5,6.4) circle (1.5pt);
\draw [fill=black] (0.54,7.6) circle (1.5pt);
\draw [fill=black] (2.08,7.6) circle (1.5pt);
\draw [fill=black] (1.42,7.6) circle (1.5pt);
\draw [fill=black] (1.46,6.54) circle (1.5pt);
\draw [fill=black] (2.98,7.6) circle (1.5pt);
\draw [fill=black] (2.12,6.52) circle (1.5pt);
\draw [fill=black] (4.52,7.58) circle (1.5pt);
\draw [fill=black] (6.62,7.56) circle (1.5pt);
\draw [fill=black] (5.58,7.58) circle (1.5pt);
\draw [fill=black] (5.26,6.54) circle (1.5pt);
\draw [fill=black] (5.92,6.52) circle (1.5pt);
\draw [fill=black] (8.92,7.56) circle (1.5pt);
\draw [fill=black] (9.3,7.56) circle (1.5pt);
\draw [fill=black] (9.68,7.54) circle (1.5pt);
\draw [fill=black] (10.22,7.54) circle (1.5pt);
\draw [fill=black] (9.1,6.44) circle (1.5pt);
\draw [fill=black] (9.74,6.36) circle (1.5pt);
\draw [fill=black] (-3.3,7.6) circle (1.5pt);
\draw [fill=black] (-1.18,7.6) circle (1.5pt);
\draw [fill=black] (-2.24,7.6) circle (1.5pt);
\draw [fill=black] (-2.24,6.52) circle (1.5pt);
\end{scriptsize}
\end{tikzpicture} 
\end{center}
\hspace{20pt} \textit{scenario (a)}                 \hspace{85pt}           \textit{scenario (b)}  \hspace{88pt}     \textit{scenario (c)}
\vspace{8pt}
\begin{center}
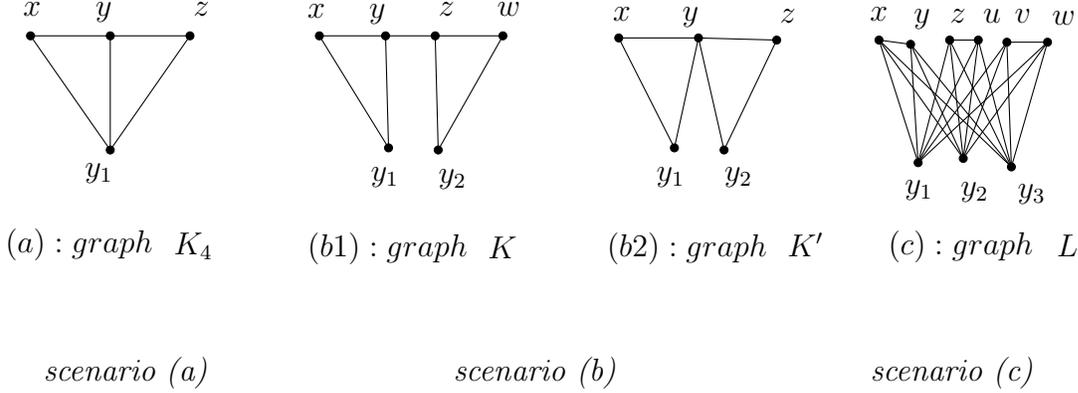
Figure 2: (a),(b1), (b2) and (c). The first three scenarios, \textit{scenario} (a),  \textit{scenario} (b) and \textit{scenario} (c) of the Standard Configuration $(n=15)$.
\end{center}

\noindent {\bf Proof.} 
\noindent As in the case of $n=16$, we get that without loss of generality $1\leq i\leq 4$, $\left|\Gamma (y_i) \cap V(C)\right| =3$ and for each $1\leq j  < j' \leq 4$, $\Gamma (y_j) \cap \Gamma (y_{j'}) \cap V(C)=\phi$.
Also $\left|\Gamma (y_5) \cap V(C)\right| \in \{ 2,3\}$. In particular, if $\left|\Gamma (y_5) \cap V(C)  \right|=2$ then, for each $1\leq j  < j' \leq 5$, $\Gamma (y_j) \cap \Gamma (y_{j'}) \cap V(C)=\phi$ and if  $\left|\Gamma (y_5) \cap V(C)\right| =3$ then, for each $1\leq j < j' \leq 5$, 
\[
\left|\Gamma (y_j) \cap \Gamma (y_{j'}) \cap V(C)\right|=
\begin{cases} 
$ $ 1 & \text{if }  j,j' \in \{4,5\}, \\
\hspace{20pt} & \\
$ $ 0 & \text{otherwise. }    \\
\end{cases}
\]

\vspace{8pt}

\noindent By Lemma \ref{l1a}, as $\delta(G) \geq 13$, we get that $\left|\Gamma (y_i) \cap V(H \setminus Y)\right| \geq 10$.  Suppose that there is some $1\leq i\leq 3$ (say $i=1$) such that $\left|\Gamma (y_i) \cap V(H \setminus Y)\right|\geq 11$. Then as $r(P_3,K_6)=11$ we get scenario $(a)$. Next, assume that for all $1\leq i\leq 3$,  $\left|\Gamma (y_i) \cap V(H \setminus Y)\right|=10$.  By the classification of the Ramsey critical $(P_3,K_6)$ graphs, we get that for all $1\leq i\leq 3$, $G[\Gamma (y_i) \cap V(H \setminus Y)] \supseteq 5K_2$. 

\vspace{8pt}

\noindent This gives  two possibilities. The first possibility is $\left|{\cup}_{i=1}^3 \Gamma (y_i) \cap V(H \setminus Y)\right|=10$. In this case, as for all $1\leq i\leq 3$, $5K_2 \subseteq G[\Gamma (y_i) \cap V(H \setminus Y)]$ we get scenario $(c)$.
The second possibility if $\left|{\cup}_{i=1}^3 \Gamma (y_i) \cup V(H \setminus Y)\right|\geq 11$. Without loss of generality, we may assume that $\left| \Gamma (y_1) \cup  \Gamma (y_2) \cup V(H \setminus Y)\right|\geq 11$.  Let $x_{11}$ be any vertex of $\Gamma (y_2)   \cap (\Gamma (y_1))^c \cup V(H \setminus Y)$. Since  $r(P_3,K_6)=11$, we get that $G[\Gamma (y_1) \cup V(H \setminus Y) \cup \{x_{11}\}]$ contains a subgraph $P$ isomorphic to a $P_3$. If  $P$ is contained in $G[\Gamma (y_1) \cup V(H \setminus Y)]$ we  get scenario $(a)$. Otherwise, $x_{11} \in P$. However, $x_{11}$ is an element of $5K_2 \subseteq G[\Gamma (y_2) \cap V(H \setminus Y)]$ and therefore, $x_{11}$ is adjacent to some other vertex say $w$ in $G[\Gamma (y_2) \cap V(H \setminus Y)]$. Depending on whether or not $w$ belongs to $V(P)$,  we get scenarios $(b2)$ or $(b1)$ respectively. Hence the Lemma.  
\vspace{4pt}

\begin{lemma}
\label{l32}
\noindent In the Standard Configuration $(n=15)$,  $y_1\in Y$ can not be a vertex of a  $K_4$  in $H$ (see Figure 2(a)).
\end{lemma}

\noindent {\bf Proof.} As indicated in Figure 2(a), let  $x \in \Gamma (y_1) \cap V(H \setminus Y)$. 
Then we get two possibilities depending on whether or not  $x$ is  adjacent to a vertex of $\{y_2,y_3,y_4,y_5\}$. In the first possibility,  $x$ is adjacent to some vertex of  $Y$ (say $y_2$). Then, as $y_1$ and $y_2$ are connected by paths of lengths 2, 3 and 4 in $H$, no pair of vertices selected from of $\Gamma (y_1) \cap V(C)$ and  $\Gamma (y_2) \cap V(C)$ can be separated by a path of length 3, 4 or 5 along the cycle $C$. However, as argued in $n=16$, we get that $\left|\Gamma (y_2) \cap V(C)\right| \leq 2$ and $\left|\Gamma (y_2) \cap V(C)\right| = 2$ only when $\{x',y',z'\}=\{u_1,u_3,u_{8}\}$, $\{x',y',z'\}=\{u_1,u_4,u_{9}\}$ or $\{x',y',z'\}=\{u_1,u_5,u_{9}\}$. This gives a contradiction. In the second possibility, re order the vertices of the cycle such that $y_1\in X$ is adjacent to $u_1$. In this ordering, suppose further that $y_1$ is also adjacent to $u_t$ where $2 \leq t \leq 14$.  By the argument used in $n=16$, we get that any pair of vertices adjacent to $y_1$ in $C$ cannot be separated by a path of length 1, 2, 3 or 4  along  $C$. However, this again leads to a contradiction.  
\vspace{14pt}

\begin{lemma}
\label{l33}
\noindent In the Standard Configuration $(n=15)$, the vertices $y_1,y_2\in Y$ can not be  vertices of the  subgraph $K$, $K'$ or $L$ in $H$ (see Figure 2).
\end{lemma}

\noindent {\bf Proof.} In the case $y_1,y_2 \in K$, since $y_1$ and $y_2$ are connected by paths of lengths 3, 4 and 5 in $H$, no pair of vertices selected from of $\Gamma (y_1) \cap V(C)$ and  $\Gamma (y_2) \cap V(C)$ can be separated by a path of length 4, 5 or 6 along the cycle $C$. The cardinality of the possible vertex sets of  $\Gamma (y_i) \cap V(C)$ ($i=1,2$), subject to this condition, are presented in Table 1. 

\begin{center}
\begin{tabular}{ ||p{1.89cm}||p{4.55cm}||p{6.8cm}|| } 
 \hline
\begin{center}$\{x',y',z'\}$ equals \end{center}
  & 
	\begin{center}$\Gamma (y_2) \cap V(C)$ is contained in \end{center}	 
	& 
	\begin{center}Cardinality of $\Gamma (y_2) \cap V(C)$ \end{center} \\

 \hline
$\{u_1,u_3,u_5\}$     &   $\{u_2,u_4\}$   &  $\left|\Gamma (y_2) \cap V(C)\right| \leq 2$  \\
 \hline 
$\{u_1,u_3,u_6\}$     &   $\{u_4\}$   &  $\left|\Gamma (y_2) \cap V(C)\right|  \leq  1$  \\
 \hline 
$\{u_1,u_3,u_7\}$     &   $\{u_4,u_{14}\}$  &  $\left|\Gamma (y_2) \cap V(C)\right|  \leq 2$  \\
 \hline 
$\{u_1,u_3,u_8\}$     &   $ $  $  $ $\phi$     &  $\left|\Gamma (y_2) \cap V(C)\right| = 0$  \\
 \hline 
$\{u_1,u_3,u_9\}$     &   $\{u_{2}\}$    &  $\left|\Gamma (y_2) \cap V(C)\right| \leq  1$  \\
 \hline 
$\{u_1,u_4,u_7\}$     &   $ $  $  $ $\phi$   &  $\left|\Gamma (y_2) \cap V(C)\right| =0$   \\
 \hline 
$\{u_1,u_4,u_8\}$     &   $ $  $  $ $\phi$  &  $\left|\Gamma (y_2) \cap V(C)\right| =0$  \\
 \hline 
$\{u_1,u_4,u_9\}$     &   $\{u_2\}$   &  $\left|\Gamma (y_2) \cap V(C)\right|  \leq 1$  \\
 \hline 
$\{u_1,u_5,u_9\}$     &   $\{u_2,u_8,u_{12}\}$  & $\left|\Gamma (y_2) \cap V(C)\right|  \leq  3$ \\
\hline 
$\{u_1,u_5,u_{10}\}$  &   $\{u_3,u_8,u_{12}\}$  & $\left|\Gamma (y_2) \cap V(C)\right|  \leq  3$ \\
 \hline 
\end{tabular}
\end{center}
\begin{center}
Table 1: Cardinality of $\Gamma (y_2) \cap V(C)$: Graph $K$.
\end{center}

\vspace{8pt}

\noindent Because  $\Gamma (y_2) \cap V(C)=3$, we are only left to deal with the last two possibilities of Table 1 for $\Gamma (y_2) \cap V(C)$. In both possibilities, $\Gamma (y_1) \cap V(C)$ will induce a $C_3$ by Lemma \ref{l2}(b). In the first possibility, $\Gamma (y_1) \cap V(C)=\{u_1,u_5,u_9\}$  gives rise to the 15-cycle given by  $(u_1,u_5,...,u_8,y_2,y,x,y_1, u_9,...,u_{14},u_1)$, a contradiction.  In the second possibility, $\Gamma (y_1) \cap V(C)=\{u_1,u_5,u_{10}\}$  gives rise to the 15-cycle given by $ (u_1,u_5,u_6,...,u_8,y_2,z,y,x,y_1,u_{10},...,u_{14},u_1)$, a contradiction.

\vspace{8pt}

\noindent In the case $y_1,y_2 \in K'$, since $y_1$ and $y_2$ are connected by paths of lengths 2, 3 or 4 in $H$, no pair of vertices selected from of $\Gamma (y_1) \cap V(C)$ and  $\Gamma (y_2) \cap V(C)$ can be separated by a path of length 3, 4 or 5 along the cycle $C$. The cardinality of the possible vertex set of  $\Gamma (y_1) \cap V(C)$ is presented in Table 2 and each of these leads to a contradiction as $\Gamma (y_2) \cap V(C)<3$

\vspace{10pt}

\begin{center}
\begin{tabular}{ ||p{1.9cm}||p{2.7cm}||p{4.7cm}|| }
 \hline
\begin{center}$\{x',y',z'\}$ equals \end{center}
  & 
	\begin{center}$\Gamma (y_2) \cap V(C)$ is contained in \end{center}
	  & 
	\begin{center}Cardinality of $\Gamma (y_2) \cap V(C)$ \end{center}
		\\
 \hline
$\{u_1,u_3,u_5\}$     &   $ $  $  $ $\phi$   &  $\left|\Gamma (y_2) \cap V(C)\right| =0$    \\
$\{u_1,u_3,u_6\}$     &   $ $  $  $ $\phi$   &  $\left|\Gamma (y_2) \cap V(C)\right| =0$   \\
$\{u_1,u_3,u_7\}$     &   $\{u_9\}$ &  $\left|\Gamma (y_2) \cap V(C)\right| \leq 1$   \\
$\{u_1,u_3,u_8\}$     &   $\{u_2,u_9\}$  &  $\left|\Gamma (y_2) \cap V(C)\right| \leq 2$    \\
$\{u_1,u_3,u_9\}$     &   $\{u_2\}$  &  $\left|\Gamma (y_2) \cap V(C)\right| \leq 1$    \\
$\{u_1,u_4,u_7\}$     &   $ $  $  $ $\phi$  &  $\left|\Gamma (y_2) \cap V(C)\right| =0$    \\
$\{u_1,u_4,u_8\}$     &   $\{u_2\}$ &  $\left|\Gamma (y_2) \cap V(C)\right| \leq 1$    \\
$\{u_1,u_4,u_9\}$     &   $\{u_2,u_3\}$ &  $\left|\Gamma (y_2) \cap V(C)\right| \leq 2$   \\
$\{u_1,u_5,u_9\}$     &   $\{u_3,u_7\}$  &  $\left|\Gamma (y_2) \cap V(C)\right| \leq 2$   \\
$\{u_1,u_5,u_{10}\}$  &   $\{u_3\}$  &  $\left|\Gamma (y_2) \cap V(C)\right| \leq 1$     \\
 \hline
\end{tabular}
\end{center}
\begin{center}
Table 2: Cardinality of $\Gamma (y_2) \cap V(C)$: Graph $K'$.
\end{center}
\vspace{4pt}

\noindent In the case $y_1,y_2 \in K'$, since  $y_1$ and $y_2$ are connected by paths of lengths 2, 3, 4 and 5 in $H$, no pair of vertices selected from of $\Gamma (y_1) \cap V(C)$ and  $\Gamma (y_2) \cap V(C)$ can be separated by  paths of length 3, 4, 5 or 6 along the cycle $C$. Thus, Table 2 will give us the required contradiction. Similarly,  in the case $y_1,y_2 \in L$, since  $y_1$ and $y_2$ are connected by paths of lengths 2, 3, 4, 5 and 6 in $H$, no pair of vertices selected from of $\Gamma (y_1) \cap V(C)$ and  $\Gamma (y_2) \cap V(C)$ can be separated by  paths of length 3, 4, 5, 6 or 7 along the cycle $C$. As before, for all possibilities  $\left|\Gamma (y_2) \cap V(C)\right| <3$, a contradiction. Thus, lemmas 5, 6 and 7 imply that $H$ cannot have an independent set of size 5.

\vspace{8pt}
\noindent Having proved that  $H$ cannot have an independent set of size 5 in all three cases $n=15,16$ and 17,  we next continue with the main proof. Since, $H$  satisfies all conditions of Lemma  \ref{l1b}, $H$ contains a $4K_{n-1}$.

\vspace{8pt}
\noindent Next we show that $V(C_{n-1})$ induced a $K_{n-1}$. Suppose that there exists two vertices  of $V(C)$, say $v$ and $w$, such that $(v,w) \not \in E(G)$. In order to avoid a $C_n$ both $v$ and $w$ will have to be adjacent to at most one vertex of each of the four copies of $K_{n-1}$ in $H$. Moreover, any vertex of any copy of $K_{n-1}$ in $H$ will have to be adjacent to at most one vertex of another copy of a $K_{n-1}$ in $H$. Thus, each copy of a $K_{n-1}$  will have at most 5 vertices adjacent to some vertex outside that of $K_{n-1}$, in $V(H) \cup \{v,w\}$. Since $(n-1)-5 \geq 1$, we can select $x_1$ in the first $K_{n-1}$, $x_2$ in the second $K_{n-1}$, $x_3$ in the third $K_{n-1}$ and $x_4$ in the fourth $K_{n-1}$ such that $\{x_1, x_2, x_3, x_4\}$ is an independent set of size four and no vertex of $\{x_1, x_2, x_3, x_4\}$ is adjacent to any vertex of $\{v, w\}$. Hence $\{x_1, x_2, x_3, x_4, v, w\}$ is an independent set of size 6, a contradiction. Therefore, we get that any two pair of vertices of $V(C)$ are  connected by an edge. Hence, $G[V(C_{n-1})]=K_{n-1}$ as required. This $K_{n-1}$ along with the $4K_{n-1}$ contained in $H$  gives the required result. 

\section{ All Ramsey $(C_n,K_6)$ critical graphs for $n \geq 15$}

\noindent We have already observed that any Ramsey $(C_n,K_6)$ critical graph will consist of a red graph containing $5K_{n-1}$, with respect to the red/blue coloring. Let $\{V_i : $ $ i \in \{1,2,...,5\} \}$ be the vertex set of  the five $K_{n-1}$ graphs.  We notice that there are two types of Ramsey $(C_n,K_6)$ critical graphs. The first type (Type1) of Ramsey $(C_n,K_6)$ critical graphs will satisfy the condition that at most one vertex of each $V_i$ is adjacent to any other vertex in $V_i^c$. The second type (Type2) of Ramsey $(C_n,K_6)$ critical graphs  will satisfy the condition that there exists a  $V_k$ for some $1 \leq k \leq 5$ such that at least two vertices of $V_k$ have neighbors in  $V_k^c$. Moreover, it is  worth noting that a  Type1 critical graph is completely determined by the structure of the external edges between $V_i$'s and not by the $\binom{n}{2}$ edges inside each of the five $V_i$'s. This fact is taken into consideration when representing the Ramsey $(C_n,K_6)$ critical graphs.
\vspace{10pt}

\noindent Each subgraph of $K_5$ generates a unique Ramsey $(C_n,K_6)$ critical graph of Type1. Thus, as illustrated in the following figure, there are 34 critical graphs ($R_i$,$ $  $1 \leq i \leq 34$) of Type1 generated by the 34 subgraphs of $K_5$. 
\vspace{2pt}

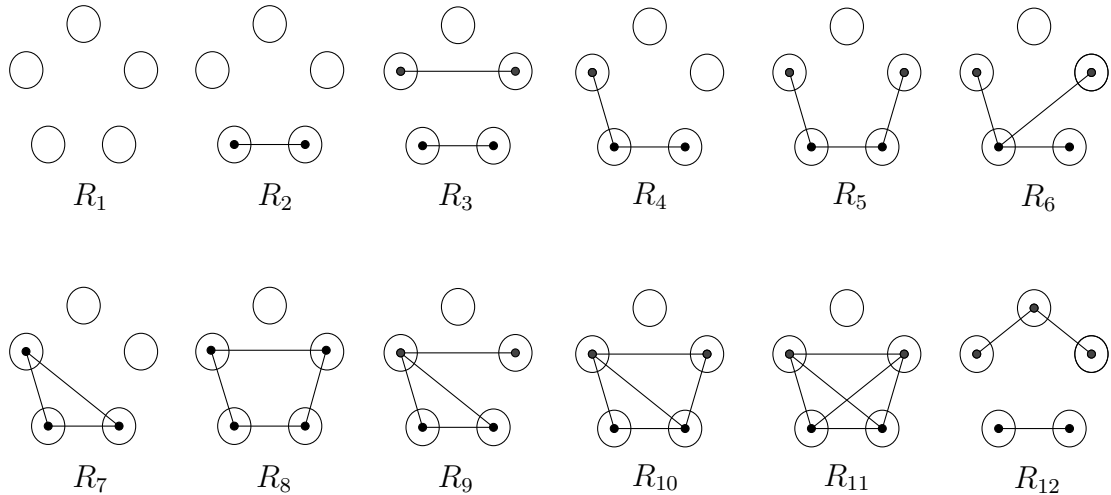
\begin{center}
\definecolor{uuuuuu}{rgb}{0.26666666666666666,0.26666666666666666,0.26666666666666666}
\begin{tikzpicture}[line cap=round,line join=round,>=triangle 45,x=0.7363281250000001cm,y=0.8138297872340424cm]
\clip(-5.600000000000002,-2.8599999999999888) rectangle (14.880000000000006,0.9000000000000069);
\draw(-4.2,0.42971746379216136) ellipse (0.22089843750000002cm and 0.24414893617021272cm);
\draw(13.995541752799934,-0.362647659142205) ellipse (0.22089843750000002cm and 0.24414893617021272cm);
\draw(-5.235541752799932,-0.3226476591422037) ellipse (0.22089843750000002cm and 0.24414893617021272cm);
\draw(-3.1644582472000677,-0.3226476591422038) ellipse (0.22089843750000002cm and 0.24414893617021272cm);
\draw(-3.56,-1.54) ellipse (0.22089843750000002cm and 0.24414893617021272cm);
\draw(-4.84,-1.54) ellipse (0.22089843750000002cm and 0.24414893617021272cm);
\draw(-0.8399999999999999,0.4297174637921619) ellipse (0.22089843750000002cm and 0.24414893617021272cm);
\draw(-1.8755417527999323,-0.32264765914220334) ellipse (0.22089843750000002cm and 0.24414893617021272cm);
\draw(0.1955417527999329,-0.32264765914220384) ellipse (0.22089843750000002cm and 0.24414893617021272cm);
\draw(-0.2,-1.54) ellipse (0.22089843750000002cm and 0.24414893617021272cm);
\draw(-1.48,-1.54) ellipse (0.22089843750000002cm and 0.24414893617021272cm);
\draw (1.92,-1.56)-- (3.2,-1.56);
\draw(2.5600000000000005,0.40971746379216223) ellipse (0.22089843750000002cm and 0.24414893617021272cm);
\draw(1.5244582472000674,-0.34264765914220313) ellipse (0.22089843750000002cm and 0.24414893617021272cm);
\draw(3.595541752799933,-0.3426476591422037) ellipse (0.22089843750000002cm and 0.24414893617021272cm);
\draw(3.2,-1.56) ellipse (0.22089843750000002cm and 0.24414893617021272cm);
\draw(1.92,-1.56) ellipse (0.22089843750000002cm and 0.24414893617021272cm);
\draw (5.38,-1.58)-- (6.66,-1.58);
\draw(6.02,0.389717463792162) ellipse (0.22089843750000002cm and 0.24414893617021272cm);
\draw(4.984458247200067,-0.3626476591422036) ellipse (0.22089843750000002cm and 0.24414893617021272cm);
\draw(7.0555417527999325,-0.3626476591422035) ellipse (0.22089843750000002cm and 0.24414893617021272cm);
\draw(6.66,-1.58) ellipse (0.22089843750000002cm and 0.24414893617021272cm);
\draw(5.38,-1.58) ellipse (0.22089843750000002cm and 0.24414893617021272cm);
\draw (4.984458247200067,-0.3626476591422036)-- (5.38,-1.58);
\draw (8.94,-1.58)-- (10.22,-1.58);
\draw(9.58,0.38971746379216354) ellipse (0.22089843750000002cm and 0.24414893617021272cm);
\draw(8.544458247200067,-0.36264765914220226) ellipse (0.22089843750000002cm and 0.24414893617021272cm);
\draw(10.615541752799935,-0.36264765914220276) ellipse (0.22089843750000002cm and 0.24414893617021272cm);
\draw(10.22,-1.58) ellipse (0.22089843750000002cm and 0.24414893617021272cm);
\draw(8.94,-1.58) ellipse (0.22089843750000002cm and 0.24414893617021272cm);
\draw (10.22,-1.58)-- (10.615541752799935,-0.36264765914220276);
\draw (8.544458247200067,-0.36264765914220226)-- (8.94,-1.58);
\draw (12.32,-1.58)-- (13.6,-1.58);
\draw(12.960000000000003,0.38971746379216143) ellipse (0.22089843750000002cm and 0.24414893617021272cm);
\draw(11.924458247200068,-0.3626476591422031) ellipse (0.22089843750000002cm and 0.24414893617021272cm);
\draw(13.995541752799934,-0.362647659142205) ellipse (0.22089843750000002cm and 0.24414893617021272cm);
\draw(13.6,-1.58) ellipse (0.22089843750000002cm and 0.24414893617021272cm);
\draw(12.32,-1.58) ellipse (0.22089843750000002cm and 0.24414893617021272cm);
\draw (11.924458247200068,-0.3626476591422031)-- (12.32,-1.58);
\draw (12.32,-1.58)-- (13.995541752799934,-0.362647659142205);
\draw (-4.600000000000001,-1.9999999999999896) node[anchor=north west] {$R_1$};
\draw (-1.2600000000000005,-1.9999999999999896) node[anchor=north west] {$R_2$};
\draw (2.000000000000001,-1.9999999999999896) node[anchor=north west] {$R_3$};
\draw (5.500000000000002,-1.9799999999999898) node[anchor=north west] {$R_4$};
\draw (9.120000000000003,-1.9799999999999898) node[anchor=north west] {$R_5$};
\draw (12.520000000000005,-2.01999999999999) node[anchor=north west] {$R_6$};
\draw (-1.48,-1.54)-- (-0.2,-1.54);
\draw (3.595541752799933,-0.3426476591422037)-- (1.5244582472000674,-0.34264765914220313);
\begin{scriptsize}
\draw [fill=black] (-1.48,-1.54) circle (1.5pt);
\draw [fill=black] (-0.2,-1.54) circle (1.5pt);
\draw [fill=black] (1.92,-1.56) circle (1.5pt);
\draw [fill=black] (3.2,-1.56) circle (1.5pt);
\draw [fill=uuuuuu] (1.5244582472000674,-0.34264765914220313) circle (1.5pt);
\draw [fill=black] (5.38,-1.58) circle (1.5pt);
\draw [fill=black] (6.66,-1.58) circle (1.5pt);
\draw [fill=uuuuuu] (4.984458247200067,-0.3626476591422036) circle (1.5pt);
\draw [fill=black] (8.94,-1.58) circle (1.5pt);
\draw [fill=black] (10.22,-1.58) circle (1.5pt);
\draw [fill=uuuuuu] (10.615541752799935,-0.36264765914220276) circle (1.5pt);
\draw [fill=uuuuuu] (8.544458247200067,-0.36264765914220226) circle (1.5pt);
\draw [fill=black] (12.32,-1.58) circle (1.5pt);
\draw [fill=black] (13.6,-1.58) circle (1.5pt);
\draw [fill=uuuuuu] (13.995541752799934,-0.362647659142205) circle (1.5pt);
\draw [fill=uuuuuu] (11.924458247200068,-0.3626476591422031) circle (1.5pt);
\draw [fill=uuuuuu] (3.595541752799933,-0.3426476591422037) circle (1.5pt);
\draw [fill=uuuuuu] (1.5244582472000674,-0.34264765914220313) circle (1.5pt);
\draw [fill=uuuuuu] (4.984458247200067,-0.3626476591422036) circle (1.5pt);
\draw [fill=uuuuuu] (10.615541752799935,-0.36264765914220276) circle (1.5pt);
\draw [fill=uuuuuu] (8.544458247200067,-0.36264765914220226) circle (1.5pt);
\draw [fill=uuuuuu] (13.995541752799934,-0.362647659142205) circle (1.5pt);
\draw [fill=uuuuuu] (11.924458247200068,-0.3626476591422031) circle (1.5pt);
\end{scriptsize}
\end{tikzpicture} 
\vspace{4pt}
\definecolor{uuuuuu}{rgb}{0.26666666666666666,0.26666666666666666,0.26666666666666666}
\begin{tikzpicture}[line cap=round,line join=round,>=triangle 45,x=0.7363281250000001cm,y=0.8138297872340424cm]
\clip(-5.600000000000002,-2.8599999999999888) rectangle (14.880000000000006,0.9000000000000069);
\draw(-4.2,0.42971746379216136) ellipse (0.22089843750000002cm and 0.24414893617021272cm);
\draw(13.995541752799934,-0.362647659142205) ellipse (0.22089843750000002cm and 0.24414893617021272cm);
\draw(-5.235541752799932,-0.3226476591422037) ellipse (0.22089843750000002cm and 0.24414893617021272cm);
\draw(-3.1644582472000677,-0.3226476591422038) ellipse (0.22089843750000002cm and 0.24414893617021272cm);
\draw(-3.56,-1.54) ellipse (0.22089843750000002cm and 0.24414893617021272cm);
\draw(-4.84,-1.54) ellipse (0.22089843750000002cm and 0.24414893617021272cm);
\draw(-0.8399999999999999,0.4297174637921619) ellipse (0.22089843750000002cm and 0.24414893617021272cm);
\draw(-1.8755417527999323,-0.32264765914220334) ellipse (0.22089843750000002cm and 0.24414893617021272cm);
\draw(0.1955417527999329,-0.32264765914220384) ellipse (0.22089843750000002cm and 0.24414893617021272cm);
\draw(-0.2,-1.54) ellipse (0.22089843750000002cm and 0.24414893617021272cm);
\draw(-1.48,-1.54) ellipse (0.22089843750000002cm and 0.24414893617021272cm);
\draw (1.92,-1.56)-- (3.2,-1.56);
\draw(2.5600000000000005,0.40971746379216223) ellipse (0.22089843750000002cm and 0.24414893617021272cm);
\draw(1.5244582472000674,-0.34264765914220313) ellipse (0.22089843750000002cm and 0.24414893617021272cm);
\draw(3.595541752799933,-0.3426476591422037) ellipse (0.22089843750000002cm and 0.24414893617021272cm);
\draw(3.2,-1.56) ellipse (0.22089843750000002cm and 0.24414893617021272cm);
\draw(1.92,-1.56) ellipse (0.22089843750000002cm and 0.24414893617021272cm);
\draw (5.38,-1.58)-- (6.66,-1.58);
\draw(6.02,0.389717463792162) ellipse (0.22089843750000002cm and 0.24414893617021272cm);
\draw(4.984458247200067,-0.3626476591422036) ellipse (0.22089843750000002cm and 0.24414893617021272cm);
\draw(7.0555417527999325,-0.3626476591422035) ellipse (0.22089843750000002cm and 0.24414893617021272cm);
\draw(6.66,-1.58) ellipse (0.22089843750000002cm and 0.24414893617021272cm);
\draw(5.38,-1.58) ellipse (0.22089843750000002cm and 0.24414893617021272cm);
\draw (4.984458247200067,-0.3626476591422036)-- (5.38,-1.58);
\draw (8.94,-1.58)-- (10.22,-1.58);
\draw(9.58,0.38971746379216354) ellipse (0.22089843750000002cm and 0.24414893617021272cm);
\draw(8.544458247200067,-0.36264765914220226) ellipse (0.22089843750000002cm and 0.24414893617021272cm);
\draw(10.615541752799935,-0.36264765914220276) ellipse (0.22089843750000002cm and 0.24414893617021272cm);
\draw(10.22,-1.58) ellipse (0.22089843750000002cm and 0.24414893617021272cm);
\draw(8.94,-1.58) ellipse (0.22089843750000002cm and 0.24414893617021272cm);
\draw (10.22,-1.58)-- (10.615541752799935,-0.36264765914220276);
\draw (8.544458247200067,-0.36264765914220226)-- (8.94,-1.58);
\draw (12.32,-1.58)-- (13.6,-1.58);
\draw(12.960000000000003,0.38971746379216143) ellipse (0.22089843750000002cm and 0.24414893617021272cm);
\draw(11.924458247200068,-0.3626476591422031) ellipse (0.22089843750000002cm and 0.24414893617021272cm);
\draw(13.995541752799934,-0.362647659142205) ellipse (0.22089843750000002cm and 0.24414893617021272cm);
\draw(13.6,-1.58) ellipse (0.22089843750000002cm and 0.24414893617021272cm);
\draw(12.32,-1.58) ellipse (0.22089843750000002cm and 0.24414893617021272cm);
\draw (-4.600000000000001,-1.9999999999999896) node[anchor=north west] {$R_7$};
\draw (-1.2800000000000005,-1.9999999999999896) node[anchor=north west] {$R_8$};
\draw (2.000000000000001,-1.9999999999999896) node[anchor=north west] {$R_9$};
\draw (5.500000000000002,-1.9799999999999898) node[anchor=north west] {$R_{10}$};
\draw (8.960000000000003,-1.9799999999999898) node[anchor=north west] {$R_{11}$};
\draw (12.400000000000004,-2.01999999999999) node[anchor=north west] {$R_{12}$};
\draw (-1.48,-1.54)-- (-0.2,-1.54);
\draw (3.595541752799933,-0.3426476591422037)-- (1.5244582472000674,-0.34264765914220313);
\draw (-5.24,-0.32)-- (-4.84,-1.54);
\draw (-3.56,-1.54)-- (-4.84,-1.54);
\draw (-3.56,-1.54)-- (-5.24,-0.32);
\draw (-1.48,-1.54)-- (-1.9,-0.3);
\draw (-1.9,-0.3)-- (0.18,-0.3);
\draw (0.18,-0.3)-- (-0.2,-1.54);
\draw (1.5244582472000674,-0.34264765914220313)-- (3.2,-1.56);
\draw (1.92,-1.56)-- (1.5244582472000674,-0.34264765914220313);
\draw (4.984458247200067,-0.3626476591422036)-- (6.66,-1.58);
\draw (6.66,-1.58)-- (7.0555417527999325,-0.3626476591422035);
\draw (7.0555417527999325,-0.3626476591422035)-- (4.984458247200067,-0.3626476591422036);
\draw (10.615541752799935,-0.36264765914220276)-- (8.544458247200067,-0.36264765914220226);
\draw (8.544458247200067,-0.36264765914220226)-- (10.22,-1.58);
\draw (10.615541752799935,-0.36264765914220276)-- (8.94,-1.58);
\draw (11.924458247200068,-0.3626476591422031)-- (12.960000000000003,0.38971746379216143);
\draw (12.960000000000003,0.38971746379216143)-- (13.995541752799934,-0.362647659142205);
\begin{scriptsize}
\draw [fill=black] (-4.84,-1.54) circle (1.5pt);
\draw [fill=black] (-3.56,-1.54) circle (1.5pt);
\draw [fill=black] (-1.48,-1.54) circle (1.5pt);
\draw [fill=black] (-0.2,-1.54) circle (1.5pt);
\draw [fill=black] (1.92,-1.56) circle (1.5pt);
\draw [fill=black] (3.2,-1.56) circle (1.5pt);
\draw [fill=uuuuuu] (1.5244582472000674,-0.34264765914220313) circle (1.5pt);
\draw [fill=black] (5.38,-1.58) circle (1.5pt);
\draw [fill=black] (6.66,-1.58) circle (1.5pt);
\draw [fill=uuuuuu] (7.0555417527999325,-0.3626476591422035) circle (1.5pt);
\draw [fill=uuuuuu] (4.984458247200067,-0.3626476591422036) circle (1.5pt);
\draw [fill=black] (8.94,-1.58) circle (1.5pt);
\draw [fill=black] (10.22,-1.58) circle (1.5pt);
\draw [fill=uuuuuu] (10.615541752799935,-0.36264765914220276) circle (1.5pt);
\draw [fill=uuuuuu] (8.544458247200067,-0.36264765914220226) circle (1.5pt);
\draw [fill=black] (12.32,-1.58) circle (1.5pt);
\draw [fill=black] (13.6,-1.58) circle (1.5pt);
\draw [fill=uuuuuu] (13.995541752799934,-0.362647659142205) circle (1.5pt);
\draw [fill=uuuuuu] (12.960000000000003,0.38971746379216143) circle (1.5pt);
\draw [fill=uuuuuu] (11.924458247200068,-0.3626476591422031) circle (1.5pt);
\draw [fill=uuuuuu] (3.595541752799933,-0.3426476591422037) circle (1.5pt);
\draw [fill=uuuuuu] (1.5244582472000674,-0.34264765914220313) circle (1.5pt);
\draw [fill=uuuuuu] (7.0555417527999325,-0.3626476591422035) circle (1.5pt);
\draw [fill=uuuuuu] (4.984458247200067,-0.3626476591422036) circle (1.5pt);
\draw [fill=uuuuuu] (10.615541752799935,-0.36264765914220276) circle (1.5pt);
\draw [fill=uuuuuu] (8.544458247200067,-0.36264765914220226) circle (1.5pt);
\draw [fill=uuuuuu] (13.995541752799934,-0.362647659142205) circle (1.5pt);
\draw [fill=uuuuuu] (12.960000000000003,0.38971746379216143) circle (1.5pt);
\draw [fill=uuuuuu] (11.924458247200068,-0.3626476591422031) circle (1.5pt);
\draw [fill=black] (-5.24,-0.32) circle (1.5pt);
\draw [fill=black] (-1.9,-0.3) circle (1.5pt);
\draw [fill=black] (0.18,-0.3) circle (1.5pt);
\end{scriptsize}
\end{tikzpicture} 
\vspace{4pt}
\definecolor{uuuuuu}{rgb}{0.26666666666666666,0.26666666666666666,0.26666666666666666}
\begin{tikzpicture}[line cap=round,line join=round,>=triangle 45,x=0.7363281250000001cm,y=0.8138297872340424cm]
\clip(-5.600000000000002,-2.8599999999999888) rectangle (14.880000000000006,0.9000000000000069);
\draw(-4.2,0.42971746379216136) ellipse (0.22089843750000002cm and 0.24414893617021272cm);
\draw(13.995541752799934,-0.362647659142205) ellipse (0.22089843750000002cm and 0.24414893617021272cm);
\draw(-5.235541752799932,-0.3226476591422037) ellipse (0.22089843750000002cm and 0.24414893617021272cm);
\draw(-3.1644582472000677,-0.3226476591422038) ellipse (0.22089843750000002cm and 0.24414893617021272cm);
\draw(-3.56,-1.54) ellipse (0.22089843750000002cm and 0.24414893617021272cm);
\draw(-4.84,-1.54) ellipse (0.22089843750000002cm and 0.24414893617021272cm);
\draw(-0.8399999999999999,0.4097174637921619) ellipse (0.22089843750000002cm and 0.24414893617021272cm);
\draw(-1.8755417527999323,-0.34264765914220335) ellipse (0.22089843750000002cm and 0.24414893617021272cm);
\draw(0.1955417527999329,-0.34264765914220385) ellipse (0.22089843750000002cm and 0.24414893617021272cm);
\draw(-0.2,-1.56) ellipse (0.22089843750000002cm and 0.24414893617021272cm);
\draw(-1.48,-1.56) ellipse (0.22089843750000002cm and 0.24414893617021272cm);
\draw(2.5600000000000005,0.40971746379216223) ellipse (0.22089843750000002cm and 0.24414893617021272cm);
\draw(1.5244582472000674,-0.34264765914220313) ellipse (0.22089843750000002cm and 0.24414893617021272cm);
\draw(3.595541752799933,-0.3426476591422037) ellipse (0.22089843750000002cm and 0.24414893617021272cm);
\draw(3.2,-1.56) ellipse (0.22089843750000002cm and 0.24414893617021272cm);
\draw(1.92,-1.56) ellipse (0.22089843750000002cm and 0.24414893617021272cm);
\draw (5.38,-1.58)-- (6.66,-1.58);
\draw(6.02,0.389717463792162) ellipse (0.22089843750000002cm and 0.24414893617021272cm);
\draw(4.984458247200067,-0.3626476591422036) ellipse (0.22089843750000002cm and 0.24414893617021272cm);
\draw(7.0555417527999325,-0.3626476591422035) ellipse (0.22089843750000002cm and 0.24414893617021272cm);
\draw(6.66,-1.58) ellipse (0.22089843750000002cm and 0.24414893617021272cm);
\draw(5.38,-1.58) ellipse (0.22089843750000002cm and 0.24414893617021272cm);
\draw (4.984458247200067,-0.3626476591422036)-- (5.38,-1.58);
\draw (8.94,-1.58)-- (10.22,-1.58);
\draw(9.58,0.38971746379216354) ellipse (0.22089843750000002cm and 0.24414893617021272cm);
\draw(8.544458247200067,-0.36264765914220226) ellipse (0.22089843750000002cm and 0.24414893617021272cm);
\draw(10.615541752799935,-0.36264765914220276) ellipse (0.22089843750000002cm and 0.24414893617021272cm);
\draw(10.22,-1.58) ellipse (0.22089843750000002cm and 0.24414893617021272cm);
\draw(8.94,-1.58) ellipse (0.22089843750000002cm and 0.24414893617021272cm);
\draw (10.22,-1.58)-- (10.615541752799935,-0.36264765914220276);
\draw (8.544458247200067,-0.36264765914220226)-- (8.94,-1.58);
\draw(12.960000000000003,0.38971746379216143) ellipse (0.22089843750000002cm and 0.24414893617021272cm);
\draw(11.924458247200068,-0.3626476591422031) ellipse (0.22089843750000002cm and 0.24414893617021272cm);
\draw(13.995541752799934,-0.362647659142205) ellipse (0.22089843750000002cm and 0.24414893617021272cm);
\draw(13.6,-1.58) ellipse (0.22089843750000002cm and 0.24414893617021272cm);
\draw(12.32,-1.58) ellipse (0.22089843750000002cm and 0.24414893617021272cm);
\draw (-4.780000000000002,-1.9799999999999898) node[anchor=north west] {$R_{13}$};
\draw (-1.4200000000000006,-1.9999999999999896) node[anchor=north west] {$R_{14}$};
\draw (2.000000000000001,-2.03999999999999) node[anchor=north west] {$R_{15}$};
\draw (5.500000000000002,-1.9799999999999898) node[anchor=north west] {$R_{16}$};
\draw (8.980000000000004,-2.01999999999999) node[anchor=north west] {$R_{17}$};
\draw (12.400000000000004,-2.01999999999999) node[anchor=north west] {$R_{18}$};
\draw (-3.56,-1.54)-- (-4.84,-1.54);
\draw (-1.48,-1.56)-- (-1.9,-0.3);
\draw (0.18,-0.3)-- (-0.2,-1.56);
\draw (6.66,-1.58)-- (7.0555417527999325,-0.3626476591422035);
\draw (10.615541752799935,-0.36264765914220276)-- (8.544458247200067,-0.36264765914220226);
\draw (11.924458247200068,-0.3626476591422031)-- (12.960000000000003,0.38971746379216143);
\draw (12.960000000000003,0.38971746379216143)-- (13.995541752799934,-0.362647659142205);
\draw (-5.24,-0.32)-- (-3.18,-0.32);
\draw (-3.18,-0.32)-- (-4.2,0.42);
\draw (-4.2,0.42)-- (-5.24,-0.32);
\draw (-0.8,0.42)-- (-1.9,-0.3);
\draw (-0.8,0.42)-- (0.18,-0.3);
\draw (1.5244582472000674,-0.34264765914220313)-- (2.54,0.42);
\draw (2.54,0.42)-- (3.595541752799933,-0.3426476591422037);
\draw (2.54,0.42)-- (3.2,-1.56);
\draw (2.54,0.42)-- (1.92,-1.56);
\draw (6.02,0.389717463792162)-- (5.38,-1.58);
\draw (9.58,0.38971746379216354)-- (8.544458247200067,-0.36264765914220226);
\draw (11.924458247200068,-0.3626476591422031)-- (13.6,-1.58);
\draw (13.995541752799934,-0.362647659142205)-- (13.6,-1.58);
\draw (13.995541752799934,-0.362647659142205)-- (12.32,-1.58);
\draw (12.32,-1.58)-- (11.924458247200068,-0.3626476591422031);
\begin{scriptsize}
\draw [fill=black] (-4.84,-1.54) circle (1.5pt);
\draw [fill=black] (-3.56,-1.54) circle (1.5pt);
\draw [fill=black] (-1.48,-1.56) circle (1.5pt);
\draw [fill=black] (-0.2,-1.56) circle (1.5pt);
\draw [fill=black] (1.92,-1.56) circle (1.5pt);
\draw [fill=black] (3.2,-1.56) circle (1.5pt);
\draw [fill=uuuuuu] (1.5244582472000674,-0.34264765914220313) circle (1.5pt);
\draw [fill=black] (5.38,-1.58) circle (1.5pt);
\draw [fill=black] (6.66,-1.58) circle (1.5pt);
\draw [fill=uuuuuu] (7.0555417527999325,-0.3626476591422035) circle (1.5pt);
\draw [fill=uuuuuu] (4.984458247200067,-0.3626476591422036) circle (1.5pt);
\draw [fill=black] (8.94,-1.58) circle (1.5pt);
\draw [fill=black] (10.22,-1.58) circle (1.5pt);
\draw [fill=uuuuuu] (10.615541752799935,-0.36264765914220276) circle (1.5pt);
\draw [fill=uuuuuu] (8.544458247200067,-0.36264765914220226) circle (1.5pt);
\draw [fill=black] (12.32,-1.58) circle (1.5pt);
\draw [fill=black] (13.6,-1.58) circle (1.5pt);
\draw [fill=uuuuuu] (13.995541752799934,-0.362647659142205) circle (1.5pt);
\draw [fill=uuuuuu] (12.960000000000003,0.38971746379216143) circle (1.5pt);
\draw [fill=uuuuuu] (11.924458247200068,-0.3626476591422031) circle (1.5pt);
\draw [fill=uuuuuu] (3.595541752799933,-0.3426476591422037) circle (1.5pt);
\draw [fill=uuuuuu] (1.5244582472000674,-0.34264765914220313) circle (1.5pt);
\draw [fill=uuuuuu] (7.0555417527999325,-0.3626476591422035) circle (1.5pt);
\draw [fill=uuuuuu] (6.02,0.389717463792162) circle (1.5pt);
\draw [fill=uuuuuu] (4.984458247200067,-0.3626476591422036) circle (1.5pt);
\draw [fill=uuuuuu] (10.615541752799935,-0.36264765914220276) circle (1.5pt);
\draw [fill=uuuuuu] (9.58,0.38971746379216354) circle (1.5pt);
\draw [fill=uuuuuu] (8.544458247200067,-0.36264765914220226) circle (1.5pt);
\draw [fill=uuuuuu] (13.995541752799934,-0.362647659142205) circle (1.5pt);
\draw [fill=uuuuuu] (12.960000000000003,0.38971746379216143) circle (1.5pt);
\draw [fill=uuuuuu] (11.924458247200068,-0.3626476591422031) circle (1.5pt);
\draw [fill=black] (-5.24,-0.32) circle (1.5pt);
\draw [fill=black] (-1.9,-0.3) circle (1.5pt);
\draw [fill=black] (0.18,-0.3) circle (1.5pt);
\draw [fill=black] (-3.18,-0.32) circle (1.5pt);
\draw [fill=black] (-4.2,0.42) circle (1.5pt);
\draw [fill=black] (-0.8,0.42) circle (1.5pt);
\draw [fill=black] (2.54,0.42) circle (1.5pt);
\end{scriptsize}
\end{tikzpicture} 
\vspace{4pt}
\definecolor{uuuuuu}{rgb}{0.26666666666666666,0.26666666666666666,0.26666666666666666}
\begin{tikzpicture}[line cap=round,line join=round,>=triangle 45,x=0.7363281250000001cm,y=0.8138297872340424cm]
\clip(-5.600000000000002,-2.8599999999999888) rectangle (14.880000000000006,0.9000000000000069);
\draw(-4.2,0.42971746379216136) ellipse (0.22089843750000002cm and 0.24414893617021272cm);
\draw(13.995541752799934,-0.362647659142205) ellipse (0.22089843750000002cm and 0.24414893617021272cm);
\draw(-5.235541752799932,-0.3226476591422037) ellipse (0.22089843750000002cm and 0.24414893617021272cm);
\draw(-3.1644582472000677,-0.3226476591422038) ellipse (0.22089843750000002cm and 0.24414893617021272cm);
\draw(-3.56,-1.54) ellipse (0.22089843750000002cm and 0.24414893617021272cm);
\draw(-4.84,-1.54) ellipse (0.22089843750000002cm and 0.24414893617021272cm);
\draw(-0.8399999999999999,0.4097174637921619) ellipse (0.22089843750000002cm and 0.24414893617021272cm);
\draw(-1.8755417527999323,-0.34264765914220335) ellipse (0.22089843750000002cm and 0.24414893617021272cm);
\draw(0.1955417527999329,-0.34264765914220385) ellipse (0.22089843750000002cm and 0.24414893617021272cm);
\draw(-0.2,-1.56) ellipse (0.22089843750000002cm and 0.24414893617021272cm);
\draw(-1.48,-1.56) ellipse (0.22089843750000002cm and 0.24414893617021272cm);
\draw(2.5600000000000005,0.40971746379216223) ellipse (0.22089843750000002cm and 0.24414893617021272cm);
\draw(1.5244582472000674,-0.34264765914220313) ellipse (0.22089843750000002cm and 0.24414893617021272cm);
\draw(3.595541752799933,-0.3426476591422037) ellipse (0.22089843750000002cm and 0.24414893617021272cm);
\draw(3.2,-1.56) ellipse (0.22089843750000002cm and 0.24414893617021272cm);
\draw(1.92,-1.56) ellipse (0.22089843750000002cm and 0.24414893617021272cm);
\draw (5.38,-1.58)-- (6.66,-1.58);
\draw(6.02,0.389717463792162) ellipse (0.22089843750000002cm and 0.24414893617021272cm);
\draw(4.984458247200067,-0.3626476591422036) ellipse (0.22089843750000002cm and 0.24414893617021272cm);
\draw(7.0555417527999325,-0.3626476591422035) ellipse (0.22089843750000002cm and 0.24414893617021272cm);
\draw(6.66,-1.58) ellipse (0.22089843750000002cm and 0.24414893617021272cm);
\draw(5.38,-1.58) ellipse (0.22089843750000002cm and 0.24414893617021272cm);
\draw (4.984458247200067,-0.3626476591422036)-- (5.38,-1.58);
\draw(9.58,0.38971746379216354) ellipse (0.22089843750000002cm and 0.24414893617021272cm);
\draw(8.544458247200067,-0.36264765914220226) ellipse (0.22089843750000002cm and 0.24414893617021272cm);
\draw(10.615541752799935,-0.36264765914220276) ellipse (0.22089843750000002cm and 0.24414893617021272cm);
\draw(10.22,-1.58) ellipse (0.22089843750000002cm and 0.24414893617021272cm);
\draw(8.94,-1.58) ellipse (0.22089843750000002cm and 0.24414893617021272cm);
\draw (10.22,-1.58)-- (10.615541752799935,-0.36264765914220276);
\draw (8.544458247200067,-0.36264765914220226)-- (8.94,-1.58);
\draw(12.960000000000003,0.38971746379216143) ellipse (0.22089843750000002cm and 0.24414893617021272cm);
\draw(11.924458247200068,-0.3626476591422031) ellipse (0.22089843750000002cm and 0.24414893617021272cm);
\draw(13.995541752799934,-0.362647659142205) ellipse (0.22089843750000002cm and 0.24414893617021272cm);
\draw(13.6,-1.58) ellipse (0.22089843750000002cm and 0.24414893617021272cm);
\draw(12.32,-1.58) ellipse (0.22089843750000002cm and 0.24414893617021272cm);
\draw (-4.780000000000002,-1.9799999999999898) node[anchor=north west] {$R_{19}$};
\draw (-1.4200000000000006,-1.9999999999999896) node[anchor=north west] {$R_{20}$};
\draw (2.000000000000001,-2.03999999999999) node[anchor=north west] {$R_{21}$};
\draw (5.500000000000002,-1.9799999999999898) node[anchor=north west] {$R_{22}$};
\draw (8.980000000000004,-1.9999999999999896) node[anchor=north west] {$R_{23}$};
\draw (12.400000000000004,-2.01999999999999) node[anchor=north west] {$R_{24}$};
\draw (-3.56,-1.54)-- (-4.84,-1.54);
\draw (-1.48,-1.56)-- (-1.9,-0.3);
\draw (0.18,-0.3)-- (-0.2,-1.56);
\draw (6.66,-1.58)-- (7.0555417527999325,-0.3626476591422035);
\draw (12.960000000000003,0.38971746379216143)-- (13.995541752799934,-0.362647659142205);
\draw (-5.24,-0.32)-- (-3.18,-0.32);
\draw (-3.18,-0.32)-- (-4.2,0.42);
\draw (-4.2,0.42)-- (-5.24,-0.32);
\draw (-0.8,0.42)-- (-1.9,-0.3);
\draw (-0.8,0.42)-- (0.18,-0.3);
\draw (1.5244582472000674,-0.34264765914220313)-- (2.54,0.42);
\draw (2.54,0.42)-- (3.595541752799933,-0.3426476591422037);
\draw (9.58,0.38971746379216354)-- (8.544458247200067,-0.36264765914220226);
\draw (13.995541752799934,-0.362647659142205)-- (13.6,-1.58);
\draw (13.995541752799934,-0.362647659142205)-- (12.32,-1.58);
\draw (12.32,-1.58)-- (11.924458247200068,-0.3626476591422031);
\draw (-5.24,-0.32)-- (-4.84,-1.54);
\draw (-1.9,-0.3)-- (0.18,-0.3);
\draw (1.5244582472000674,-0.34264765914220313)-- (3.595541752799933,-0.3426476591422037);
\draw (1.5244582472000674,-0.34264765914220313)-- (3.2,-1.56);
\draw (1.5244582472000674,-0.34264765914220313)-- (1.92,-1.56);
\draw (4.984458247200067,-0.3626476591422036)-- (6.02,0.389717463792162);
\draw (6.02,0.389717463792162)-- (7.0555417527999325,-0.3626476591422035);
\draw (9.58,0.38971746379216354)-- (10.615541752799935,-0.36264765914220276);
\draw (10.22,-1.58)-- (9.58,0.38971746379216354);
\draw (9.58,0.38971746379216354)-- (8.94,-1.58);
\draw (13.995541752799934,-0.362647659142205)-- (11.924458247200068,-0.3626476591422031);
\draw (13.6,-1.58)-- (12.32,-1.58);
\begin{scriptsize}
\draw [fill=black] (-4.84,-1.54) circle (1.5pt);
\draw [fill=black] (-3.56,-1.54) circle (1.5pt);
\draw [fill=black] (-1.48,-1.56) circle (1.5pt);
\draw [fill=black] (-0.2,-1.56) circle (1.5pt);
\draw [fill=black] (1.92,-1.56) circle (1.5pt);
\draw [fill=black] (3.2,-1.56) circle (1.5pt);
\draw [fill=uuuuuu] (1.5244582472000674,-0.34264765914220313) circle (1.5pt);
\draw [fill=black] (5.38,-1.58) circle (1.5pt);
\draw [fill=black] (6.66,-1.58) circle (1.5pt);
\draw [fill=uuuuuu] (7.0555417527999325,-0.3626476591422035) circle (1.5pt);
\draw [fill=uuuuuu] (4.984458247200067,-0.3626476591422036) circle (1.5pt);
\draw [fill=black] (8.94,-1.58) circle (1.5pt);
\draw [fill=black] (10.22,-1.58) circle (1.5pt);
\draw [fill=uuuuuu] (10.615541752799935,-0.36264765914220276) circle (1.5pt);
\draw [fill=uuuuuu] (8.544458247200067,-0.36264765914220226) circle (1.5pt);
\draw [fill=black] (12.32,-1.58) circle (1.5pt);
\draw [fill=black] (13.6,-1.58) circle (1.5pt);
\draw [fill=uuuuuu] (13.995541752799934,-0.362647659142205) circle (1.5pt);
\draw [fill=uuuuuu] (12.960000000000003,0.38971746379216143) circle (1.5pt);
\draw [fill=uuuuuu] (11.924458247200068,-0.3626476591422031) circle (1.5pt);
\draw [fill=uuuuuu] (3.595541752799933,-0.3426476591422037) circle (1.5pt);
\draw [fill=uuuuuu] (1.5244582472000674,-0.34264765914220313) circle (1.5pt);
\draw [fill=uuuuuu] (7.0555417527999325,-0.3626476591422035) circle (1.5pt);
\draw [fill=uuuuuu] (6.02,0.389717463792162) circle (1.5pt);
\draw [fill=uuuuuu] (4.984458247200067,-0.3626476591422036) circle (1.5pt);
\draw [fill=uuuuuu] (10.615541752799935,-0.36264765914220276) circle (1.5pt);
\draw [fill=uuuuuu] (9.58,0.38971746379216354) circle (1.5pt);
\draw [fill=uuuuuu] (8.544458247200067,-0.36264765914220226) circle (1.5pt);
\draw [fill=uuuuuu] (13.995541752799934,-0.362647659142205) circle (1.5pt);
\draw [fill=uuuuuu] (12.960000000000003,0.38971746379216143) circle (1.5pt);
\draw [fill=uuuuuu] (11.924458247200068,-0.3626476591422031) circle (1.5pt);
\draw [fill=black] (-5.24,-0.32) circle (1.5pt);
\draw [fill=black] (-1.9,-0.3) circle (1.5pt);
\draw [fill=black] (0.18,-0.3) circle (1.5pt);
\draw [fill=black] (-3.18,-0.32) circle (1.5pt);
\draw [fill=black] (-4.2,0.42) circle (1.5pt);
\draw [fill=black] (-0.8,0.42) circle (1.5pt);
\draw [fill=black] (2.54,0.42) circle (1.5pt);
\end{scriptsize}
\end{tikzpicture} 
\vspace{4pt}
\definecolor{uuuuuu}{rgb}{0.26666666666666666,0.26666666666666666,0.26666666666666666}
\begin{tikzpicture}[line cap=round,line join=round,>=triangle 45,x=0.7363281250000001cm,y=0.8138297872340424cm]
\clip(-5.600000000000002,-2.8599999999999888) rectangle (14.880000000000006,0.9000000000000069);
\draw(-4.2,0.42971746379216136) ellipse (0.22089843750000002cm and 0.24414893617021272cm);
\draw(13.995541752799934,-0.362647659142205) ellipse (0.22089843750000002cm and 0.24414893617021272cm);
\draw(-5.235541752799932,-0.3226476591422037) ellipse (0.22089843750000002cm and 0.24414893617021272cm);
\draw(-3.1644582472000677,-0.3226476591422038) ellipse (0.22089843750000002cm and 0.24414893617021272cm);
\draw(-3.56,-1.54) ellipse (0.22089843750000002cm and 0.24414893617021272cm);
\draw(-4.84,-1.54) ellipse (0.22089843750000002cm and 0.24414893617021272cm);
\draw(-0.8399999999999999,0.4097174637921619) ellipse (0.22089843750000002cm and 0.24414893617021272cm);
\draw(-1.8755417527999323,-0.34264765914220335) ellipse (0.22089843750000002cm and 0.24414893617021272cm);
\draw(0.1955417527999329,-0.34264765914220385) ellipse (0.22089843750000002cm and 0.24414893617021272cm);
\draw(-0.2,-1.56) ellipse (0.22089843750000002cm and 0.24414893617021272cm);
\draw(-1.48,-1.56) ellipse (0.22089843750000002cm and 0.24414893617021272cm);
\draw(2.5600000000000005,0.40971746379216223) ellipse (0.22089843750000002cm and 0.24414893617021272cm);
\draw(1.5244582472000674,-0.34264765914220313) ellipse (0.22089843750000002cm and 0.24414893617021272cm);
\draw(3.595541752799933,-0.3426476591422037) ellipse (0.22089843750000002cm and 0.24414893617021272cm);
\draw(3.2,-1.56) ellipse (0.22089843750000002cm and 0.24414893617021272cm);
\draw(1.92,-1.56) ellipse (0.22089843750000002cm and 0.24414893617021272cm);
\draw (5.38,-1.58)-- (6.66,-1.58);
\draw(6.02,0.389717463792162) ellipse (0.22089843750000002cm and 0.24414893617021272cm);
\draw(4.984458247200067,-0.3626476591422036) ellipse (0.22089843750000002cm and 0.24414893617021272cm);
\draw(7.0555417527999325,-0.3626476591422035) ellipse (0.22089843750000002cm and 0.24414893617021272cm);
\draw(6.66,-1.58) ellipse (0.22089843750000002cm and 0.24414893617021272cm);
\draw(5.38,-1.58) ellipse (0.22089843750000002cm and 0.24414893617021272cm);
\draw (4.984458247200067,-0.3626476591422036)-- (5.38,-1.58);
\draw(9.58,0.38971746379216354) ellipse (0.22089843750000002cm and 0.24414893617021272cm);
\draw(8.544458247200067,-0.36264765914220226) ellipse (0.22089843750000002cm and 0.24414893617021272cm);
\draw(10.615541752799935,-0.36264765914220276) ellipse (0.22089843750000002cm and 0.24414893617021272cm);
\draw(10.22,-1.58) ellipse (0.22089843750000002cm and 0.24414893617021272cm);
\draw(8.94,-1.58) ellipse (0.22089843750000002cm and 0.24414893617021272cm);
\draw (10.22,-1.58)-- (10.615541752799935,-0.36264765914220276);
\draw (8.544458247200067,-0.36264765914220226)-- (8.94,-1.58);
\draw(12.960000000000003,0.38971746379216143) ellipse (0.22089843750000002cm and 0.24414893617021272cm);
\draw(11.924458247200068,-0.3626476591422031) ellipse (0.22089843750000002cm and 0.24414893617021272cm);
\draw(13.995541752799934,-0.362647659142205) ellipse (0.22089843750000002cm and 0.24414893617021272cm);
\draw(13.6,-1.58) ellipse (0.22089843750000002cm and 0.24414893617021272cm);
\draw(12.32,-1.58) ellipse (0.22089843750000002cm and 0.24414893617021272cm);
\draw (-4.760000000000002,-1.9799999999999898) node[anchor=north west] {$R_{25}$};
\draw (-1.4200000000000006,-1.9999999999999896) node[anchor=north west] {$R_{26}$};
\draw (2.000000000000001,-2.03999999999999) node[anchor=north west] {$R_{27}$};
\draw (5.500000000000002,-1.9799999999999898) node[anchor=north west] {$R_{28}$};
\draw (9.000000000000004,-1.9999999999999896) node[anchor=north west] {$R_{29}$};
\draw (12.380000000000004,-2.01999999999999) node[anchor=north west] {$R_{30}$};
\draw (-3.56,-1.54)-- (-4.84,-1.54);
\draw (-1.48,-1.56)-- (-1.9,-0.3);
\draw (0.18,-0.3)-- (-0.2,-1.56);
\draw (6.66,-1.58)-- (7.0555417527999325,-0.3626476591422035);
\draw (-5.24,-0.32)-- (-3.18,-0.32);
\draw (-4.2,0.42)-- (-5.24,-0.32);
\draw (-0.8,0.42)-- (-1.9,-0.3);
\draw (-0.8,0.42)-- (0.18,-0.3);
\draw (1.5244582472000674,-0.34264765914220313)-- (2.54,0.42);
\draw (2.54,0.42)-- (3.595541752799933,-0.3426476591422037);
\draw (9.58,0.38971746379216354)-- (8.544458247200067,-0.36264765914220226);
\draw (13.995541752799934,-0.362647659142205)-- (13.6,-1.58);
\draw (13.995541752799934,-0.362647659142205)-- (12.32,-1.58);
\draw (12.32,-1.58)-- (11.924458247200068,-0.3626476591422031);
\draw (-5.24,-0.32)-- (-4.84,-1.54);
\draw (-1.9,-0.3)-- (0.18,-0.3);
\draw (1.5244582472000674,-0.34264765914220313)-- (1.92,-1.56);
\draw (4.984458247200067,-0.3626476591422036)-- (6.02,0.389717463792162);
\draw (6.02,0.389717463792162)-- (7.0555417527999325,-0.3626476591422035);
\draw (9.58,0.38971746379216354)-- (10.615541752799935,-0.36264765914220276);
\draw (13.995541752799934,-0.362647659142205)-- (11.924458247200068,-0.3626476591422031);
\draw (13.6,-1.58)-- (12.32,-1.58);
\draw (-3.18,-0.32)-- (-4.84,-1.54);
\draw (-3.18,-0.32)-- (-3.56,-1.54);
\draw (-1.48,-1.56)-- (-0.2,-1.56);
\draw (1.92,-1.56)-- (3.2,-1.56);
\draw (3.2,-1.56)-- (3.595541752799933,-0.3426476591422037);
\draw (2.54,0.42)-- (1.92,-1.56);
\draw (2.54,0.42)-- (3.2,-1.56);
\draw (6.66,-1.58)-- (4.984458247200067,-0.3626476591422036);
\draw (5.38,-1.58)-- (7.0555417527999325,-0.3626476591422035);
\draw (8.94,-1.58)-- (10.615541752799935,-0.36264765914220276);
\draw (10.22,-1.58)-- (8.544458247200067,-0.36264765914220226);
\draw (8.544458247200067,-0.36264765914220226)-- (10.615541752799935,-0.36264765914220276);
\draw (11.924458247200068,-0.3626476591422031)-- (13.6,-1.58);
\draw (12.960000000000003,0.38971746379216143)-- (11.924458247200068,-0.3626476591422031);
\begin{scriptsize}
\draw [fill=black] (-4.84,-1.54) circle (1.5pt);
\draw [fill=black] (-3.56,-1.54) circle (1.5pt);
\draw [fill=black] (-1.48,-1.56) circle (1.5pt);
\draw [fill=black] (-0.2,-1.56) circle (1.5pt);
\draw [fill=black] (1.92,-1.56) circle (1.5pt);
\draw [fill=black] (3.2,-1.56) circle (1.5pt);
\draw [fill=uuuuuu] (1.5244582472000674,-0.34264765914220313) circle (1.5pt);
\draw [fill=black] (5.38,-1.58) circle (1.5pt);
\draw [fill=black] (6.66,-1.58) circle (1.5pt);
\draw [fill=uuuuuu] (7.0555417527999325,-0.3626476591422035) circle (1.5pt);
\draw [fill=uuuuuu] (4.984458247200067,-0.3626476591422036) circle (1.5pt);
\draw [fill=black] (8.94,-1.58) circle (1.5pt);
\draw [fill=black] (10.22,-1.58) circle (1.5pt);
\draw [fill=uuuuuu] (10.615541752799935,-0.36264765914220276) circle (1.5pt);
\draw [fill=uuuuuu] (8.544458247200067,-0.36264765914220226) circle (1.5pt);
\draw [fill=black] (12.32,-1.58) circle (1.5pt);
\draw [fill=black] (13.6,-1.58) circle (1.5pt);
\draw [fill=uuuuuu] (13.995541752799934,-0.362647659142205) circle (1.5pt);
\draw [fill=uuuuuu] (12.960000000000003,0.38971746379216143) circle (1.5pt);
\draw [fill=uuuuuu] (11.924458247200068,-0.3626476591422031) circle (1.5pt);
\draw [fill=uuuuuu] (3.595541752799933,-0.3426476591422037) circle (1.5pt);
\draw [fill=uuuuuu] (1.5244582472000674,-0.34264765914220313) circle (1.5pt);
\draw [fill=uuuuuu] (7.0555417527999325,-0.3626476591422035) circle (1.5pt);
\draw [fill=uuuuuu] (6.02,0.389717463792162) circle (1.5pt);
\draw [fill=uuuuuu] (4.984458247200067,-0.3626476591422036) circle (1.5pt);
\draw [fill=uuuuuu] (10.615541752799935,-0.36264765914220276) circle (1.5pt);
\draw [fill=uuuuuu] (9.58,0.38971746379216354) circle (1.5pt);
\draw [fill=uuuuuu] (8.544458247200067,-0.36264765914220226) circle (1.5pt);
\draw [fill=uuuuuu] (13.995541752799934,-0.362647659142205) circle (1.5pt);
\draw [fill=uuuuuu] (12.960000000000003,0.38971746379216143) circle (1.5pt);
\draw [fill=uuuuuu] (11.924458247200068,-0.3626476591422031) circle (1.5pt);
\draw [fill=black] (-5.24,-0.32) circle (1.5pt);
\draw [fill=black] (-1.9,-0.3) circle (1.5pt);
\draw [fill=black] (0.18,-0.3) circle (1.5pt);
\draw [fill=black] (-3.18,-0.32) circle (1.5pt);
\draw [fill=black] (-4.2,0.42) circle (1.5pt);
\draw [fill=black] (-0.8,0.42) circle (1.5pt);
\draw [fill=black] (2.54,0.42) circle (1.5pt);
\end{scriptsize}
\end{tikzpicture} 
\vspace{4pt}
\definecolor{uuuuuu}{rgb}{0.26666666666666666,0.26666666666666666,0.26666666666666666}
\begin{tikzpicture}[line cap=round,line join=round,>=triangle 45,x=0.7360703812316717cm,y=0.8138297872340424cm]
\clip(-5.600000000000002,-2.8599999999999888) rectangle (14.860000000000005,0.9000000000000069);
\draw(-0.7399999999999998,0.3297174637921617) ellipse (0.2208211143695015cm and 0.24414893617021272cm);
\draw(-1.7755417527999324,-0.42264765914220337) ellipse (0.2208211143695015cm and 0.24414893617021272cm);
\draw(0.29554175279993244,-0.42264765914220376) ellipse (0.2208211143695015cm and 0.24414893617021272cm);
\draw(-0.1,-1.64) ellipse (0.2208211143695015cm and 0.24414893617021272cm);
\draw(-1.38,-1.64) ellipse (0.2208211143695015cm and 0.24414893617021272cm);
\draw(2.6200000000000006,0.3097174637921619) ellipse (0.2208211143695015cm and 0.24414893617021272cm);
\draw(1.5844582472000677,-0.44264765914220316) ellipse (0.2208211143695015cm and 0.24414893617021272cm);
\draw(3.655541752799933,-0.4426476591422039) ellipse (0.2208211143695015cm and 0.24414893617021272cm);
\draw(3.26,-1.66) ellipse (0.2208211143695015cm and 0.24414893617021272cm);
\draw(1.98,-1.66) ellipse (0.2208211143695015cm and 0.24414893617021272cm);
\draw(6.02,0.30971746379216214) ellipse (0.2208211143695015cm and 0.24414893617021272cm);
\draw(4.984458247200067,-0.44264765914220344) ellipse (0.2208211143695015cm and 0.24414893617021272cm);
\draw(7.0555417527999325,-0.44264765914220333) ellipse (0.2208211143695015cm and 0.24414893617021272cm);
\draw(6.66,-1.66) ellipse (0.2208211143695015cm and 0.24414893617021272cm);
\draw(5.38,-1.66) ellipse (0.2208211143695015cm and 0.24414893617021272cm);
\draw (8.84,-1.68)-- (10.12,-1.68);
\draw(9.480000000000002,0.28971746379216157) ellipse (0.2208211143695015cm and 0.24414893617021272cm);
\draw(8.444458247200068,-0.46264765914220296) ellipse (0.2208211143695015cm and 0.24414893617021272cm);
\draw(10.515541752799933,-0.46264765914220485) ellipse (0.2208211143695015cm and 0.24414893617021272cm);
\draw(10.12,-1.68) ellipse (0.2208211143695015cm and 0.24414893617021272cm);
\draw(8.84,-1.68) ellipse (0.2208211143695015cm and 0.24414893617021272cm);
\draw (8.444458247200068,-0.46264765914220296)-- (8.84,-1.68);
\draw (-1.3000000000000005,-2.0799999999999894) node[anchor=north west] {$R_{31}$};
\draw (2.040000000000001,-2.0999999999999894) node[anchor=north west] {$R_{32}$};
\draw (5.480000000000002,-2.1199999999999894) node[anchor=north west] {$R_{33}$};
\draw (8.960000000000003,-2.0799999999999894) node[anchor=north west] {$R_{34}$};
\draw (-0.1,-1.64)-- (-1.38,-1.64);
\draw (1.98,-1.66)-- (1.56,-0.4);
\draw (3.64,-0.4)-- (3.26,-1.66);
\draw (10.12,-1.68)-- (10.515541752799933,-0.46264765914220485);
\draw (-1.78,-0.42)-- (0.28,-0.42);
\draw (-0.74,0.32)-- (-1.78,-0.42);
\draw (2.66,0.32)-- (1.56,-0.4);
\draw (2.66,0.32)-- (3.64,-0.4);
\draw (4.984458247200067,-0.44264765914220344)-- (6.0,0.32);
\draw (6.0,0.32)-- (7.0555417527999325,-0.44264765914220333);
\draw (-1.78,-0.42)-- (-1.38,-1.64);
\draw (1.56,-0.4)-- (3.64,-0.4);
\draw (4.984458247200067,-0.44264765914220344)-- (5.38,-1.66);
\draw (8.444458247200068,-0.46264765914220296)-- (9.480000000000002,0.28971746379216157);
\draw (9.480000000000002,0.28971746379216157)-- (10.515541752799933,-0.46264765914220485);
\draw (0.28,-0.42)-- (-1.38,-1.64);
\draw (0.28,-0.42)-- (-0.1,-1.64);
\draw (1.98,-1.66)-- (3.26,-1.66);
\draw (6.66,-1.66)-- (7.0555417527999325,-0.44264765914220333);
\draw (6.0,0.32)-- (5.38,-1.66);
\draw (6.0,0.32)-- (6.66,-1.66);
\draw (10.12,-1.68)-- (8.444458247200068,-0.46264765914220296);
\draw (8.84,-1.68)-- (10.515541752799933,-0.46264765914220485);
\draw (-1.78,-0.42)-- (-0.1,-1.64);
\draw (0.28,-0.42)-- (-0.74,0.32);
\draw (1.98,-1.66)-- (2.66,0.32);
\draw (2.66,0.32)-- (3.26,-1.66);
\draw (4.984458247200067,-0.44264765914220344)-- (7.0555417527999325,-0.44264765914220333);
\draw (7.0555417527999325,-0.44264765914220333)-- (5.38,-1.66);
\draw (4.984458247200067,-0.44264765914220344)-- (6.66,-1.66);
\draw (8.444458247200068,-0.46264765914220296)-- (10.515541752799933,-0.46264765914220485);
\draw (8.84,-1.68)-- (9.480000000000002,0.28971746379216157);
\draw (9.480000000000002,0.28971746379216157)-- (10.12,-1.68);
\begin{scriptsize}
\draw [fill=black] (-1.38,-1.64) circle (1.5pt);
\draw [fill=black] (-0.1,-1.64) circle (1.5pt);
\draw [fill=black] (1.98,-1.66) circle (1.5pt);
\draw [fill=black] (3.26,-1.66) circle (1.5pt);
\draw [fill=black] (5.38,-1.66) circle (1.5pt);
\draw [fill=black] (6.66,-1.66) circle (1.5pt);
\draw [fill=uuuuuu] (4.984458247200067,-0.44264765914220344) circle (1.5pt);
\draw [fill=black] (8.84,-1.68) circle (1.5pt);
\draw [fill=black] (10.12,-1.68) circle (1.5pt);
\draw [fill=uuuuuu] (10.515541752799933,-0.46264765914220485) circle (1.5pt);
\draw [fill=uuuuuu] (8.444458247200068,-0.46264765914220296) circle (1.5pt);
\draw [fill=uuuuuu] (7.0555417527999325,-0.44264765914220333) circle (1.5pt);
\draw [fill=uuuuuu] (4.984458247200067,-0.44264765914220344) circle (1.5pt);
\draw [fill=uuuuuu] (10.515541752799933,-0.46264765914220485) circle (1.5pt);
\draw [fill=uuuuuu] (9.480000000000002,0.28971746379216157) circle (1.5pt);
\draw [fill=uuuuuu] (8.444458247200068,-0.46264765914220296) circle (1.5pt);
\draw [fill=black] (-1.78,-0.42) circle (1.5pt);
\draw [fill=black] (1.56,-0.4) circle (1.5pt);
\draw [fill=black] (3.64,-0.4) circle (1.5pt);
\draw [fill=black] (0.28,-0.42) circle (1.5pt);
\draw [fill=black] (-0.74,0.32) circle (1.5pt);
\draw [fill=black] (2.66,0.32) circle (1.5pt);
\draw [fill=black] (6.0,0.32) circle (1.5pt);
\end{scriptsize}
\end{tikzpicture} 
\vspace{12pt}

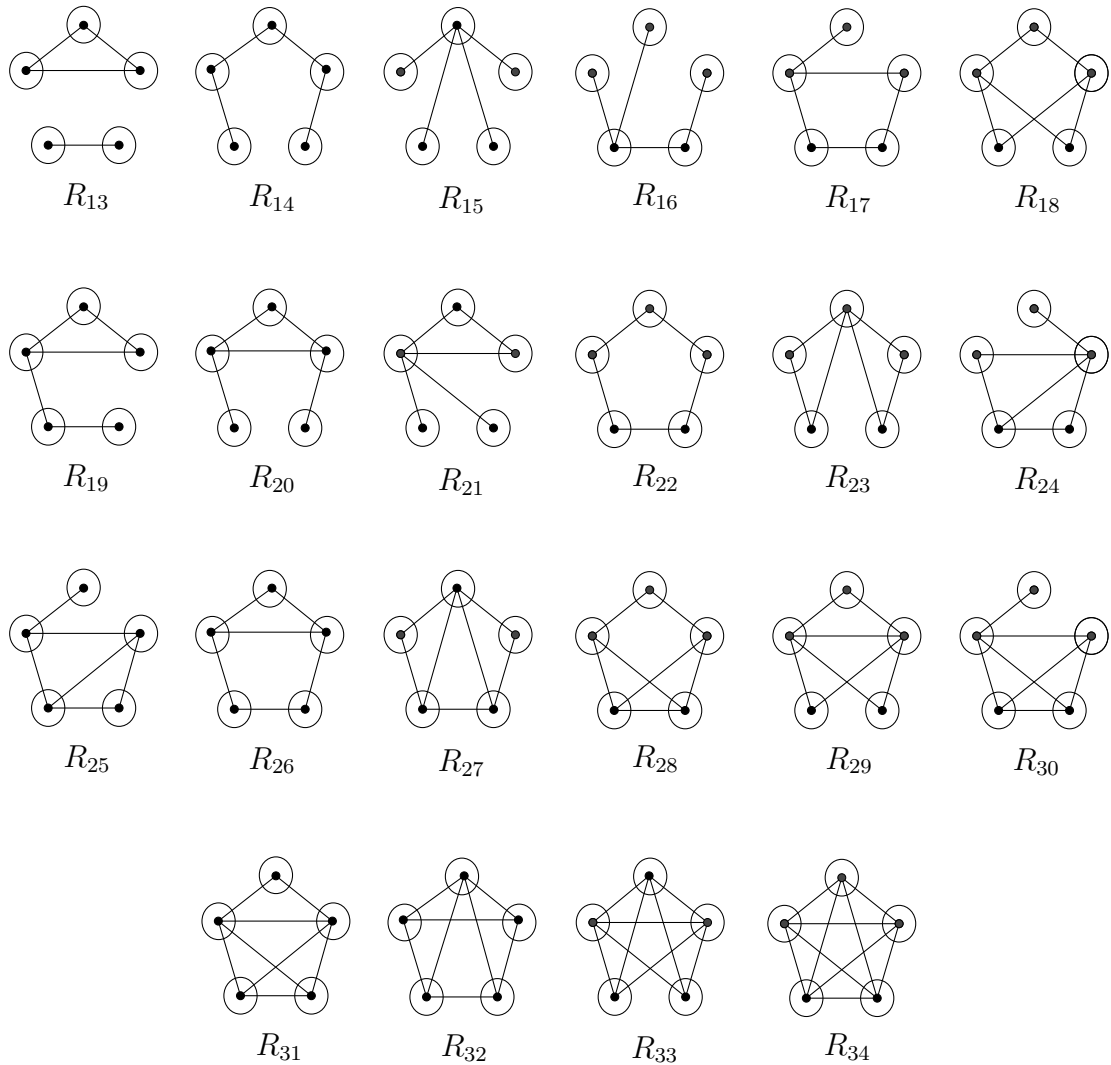
Figure 3. Ramsey $(C_n,K_6)$ critical graphs of Type1, $R_i$ ($1 \leq i \leq 34$) 
\end{center}
\vspace{8pt}

\noindent First note that each and every Type2 critical graph is obtained by an appropriate vertex splitting of some Type1 critical graph. As illustrated in the Figure 4, there are exactly 34 Type2 critical graphs (labeled $S_i$ where $1 \leq i \leq 34$) generated by 18 critical graphs of Type1, since exactly sixteen Type1 critical graphs do not generate  Type2 critical graphs. 

\vspace{2pt}

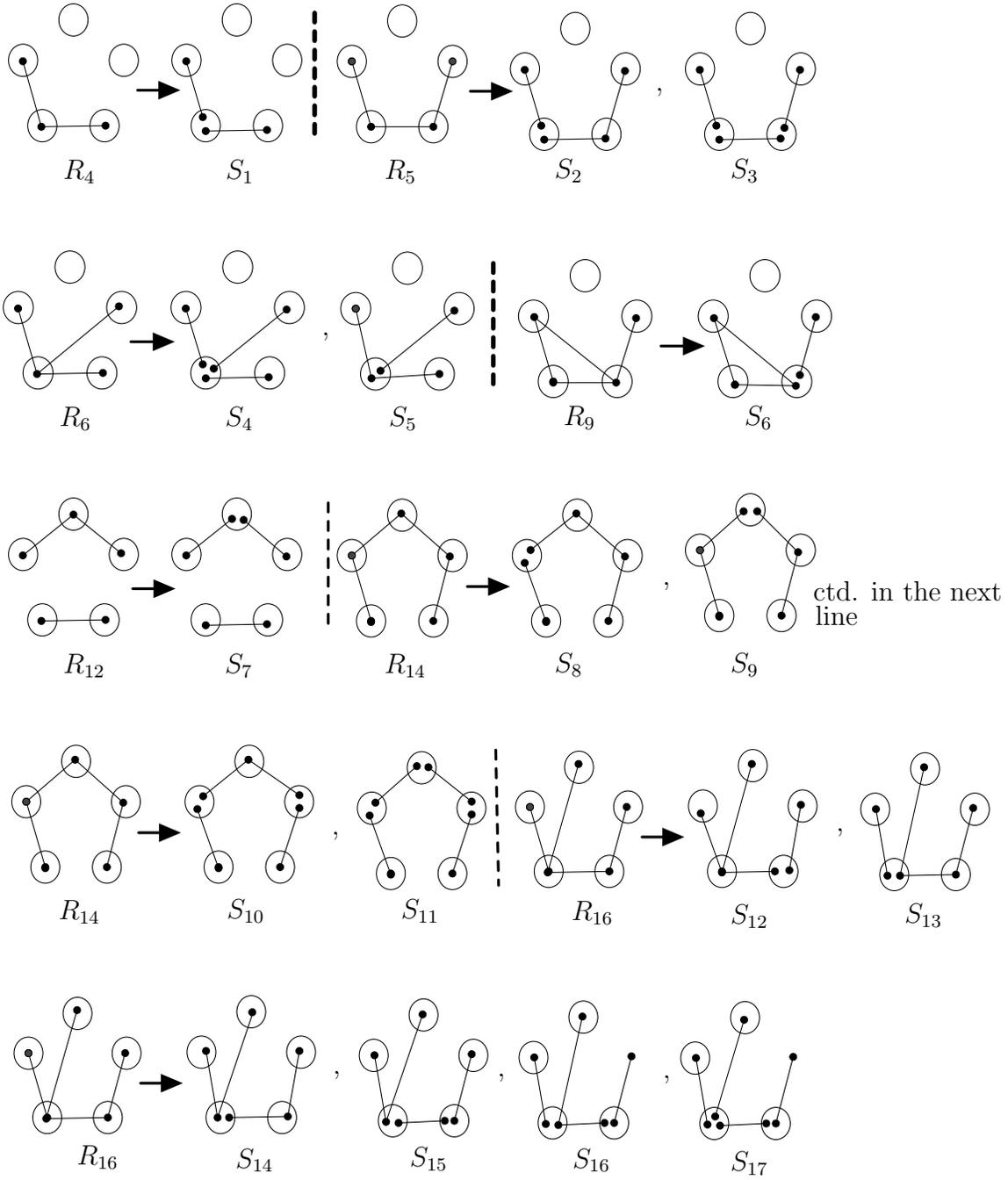
\begin{center}
\definecolor{uuuuuu}{rgb}{0.26666666666666666,0.26666666666666666,0.26666666666666666}
\begin{tikzpicture}[line cap=round,line join=round,>=triangle 45,x=0.7363281250000001cm,y=0.8138297872340424cm]
\clip(-5.600000000000002,-2.8599999999999888) rectangle (14.880000000000006,0.9000000000000069);
\draw(-4.2,0.42971746379216136) ellipse (0.22089843750000002cm and 0.24414893617021272cm);
\draw(-5.235541752799932,-0.3226476591422037) ellipse (0.22089843750000002cm and 0.24414893617021272cm);
\draw(-3.1644582472000677,-0.3226476591422038) ellipse (0.22089843750000002cm and 0.24414893617021272cm);
\draw(-3.56,-1.54) ellipse (0.22089843750000002cm and 0.24414893617021272cm);
\draw(-4.84,-1.54) ellipse (0.22089843750000002cm and 0.24414893617021272cm);
\draw(-0.8399999999999999,0.4297174637921619) ellipse (0.22089843750000002cm and 0.24414893617021272cm);
\draw(-1.8755417527999323,-0.32264765914220334) ellipse (0.22089843750000002cm and 0.24414893617021272cm);
\draw(0.1955417527999329,-0.32264765914220384) ellipse (0.22089843750000002cm and 0.24414893617021272cm);
\draw(-0.2,-1.54) ellipse (0.22089843750000002cm and 0.24414893617021272cm);
\draw(-1.48,-1.54) ellipse (0.22089843750000002cm and 0.24414893617021272cm);
\draw (1.92,-1.56)-- (3.2,-1.56);
\draw(2.5600000000000005,0.40971746379216223) ellipse (0.22089843750000002cm and 0.24414893617021272cm);
\draw(1.5244582472000674,-0.34264765914220313) ellipse (0.22089843750000002cm and 0.24414893617021272cm);
\draw(3.595541752799933,-0.3426476591422037) ellipse (0.22089843750000002cm and 0.24414893617021272cm);
\draw(3.2,-1.56) ellipse (0.22089843750000002cm and 0.24414893617021272cm);
\draw(1.92,-1.56) ellipse (0.22089843750000002cm and 0.24414893617021272cm);
\draw (-4.600000000000001,-1.9999999999999896) node[anchor=north west] {$R_4$};
\draw (-1.2600000000000005,-1.9999999999999896) node[anchor=north west] {$S_1$};
\draw (2.000000000000001,-1.9999999999999896) node[anchor=north west] {$R_5$};
\draw (5.500000000000002,-1.9799999999999898) node[anchor=north west] {$S_2$};
\draw (9.120000000000003,-1.9799999999999898) node[anchor=north west] {$S_3$};
\draw [->,line width=1.2000000000000002pt] (3.96,-0.9) -- (4.84,-0.9);
\draw (7.580000000000003,-0.6599999999999913) node[anchor=north west] {,};
\draw [line width=2.0pt,dash pattern=on 4pt off 4pt] (0.76,0.56)-- (0.76,-1.84);
\draw (-5.24,-0.34)-- (-4.86,-1.56);
\draw (-4.86,-1.56)-- (-3.54,-1.54);
\draw (-1.88,-0.34)-- (-1.54,-1.38);
\draw (-0.22,-1.62)-- (-1.48,-1.64);
\draw (1.5244582472000674,-0.34264765914220313)-- (1.92,-1.56);
\draw (3.2,-1.56)-- (3.595541752799933,-0.3426476591422037);
\draw(6.12,0.2697174637921611) ellipse (0.22089843750000002cm and 0.24414893617021272cm);
\draw(5.084458247200068,-0.48264765914220387) ellipse (0.22089843750000002cm and 0.24414893617021272cm);
\draw(7.155541752799932,-0.48264765914220437) ellipse (0.22089843750000002cm and 0.24414893617021272cm);
\draw(6.76,-1.7) ellipse (0.22089843750000002cm and 0.24414893617021272cm);
\draw(5.48,-1.7) ellipse (0.22089843750000002cm and 0.24414893617021272cm);
\draw (5.08,-0.5)-- (5.42,-1.54);
\draw (6.74,-1.78)-- (5.48,-1.8);
\draw (6.74,-1.78)-- (7.14,-0.52);
\draw(9.719999999999999,0.26971746379216044) ellipse (0.22089843750000002cm and 0.24414893617021272cm);
\draw(8.684458247200068,-0.48264765914220464) ellipse (0.22089843750000002cm and 0.24414893617021272cm);
\draw(10.75554175279993,-0.4826476591422038) ellipse (0.22089843750000002cm and 0.24414893617021272cm);
\draw(10.36,-1.7) ellipse (0.22089843750000002cm and 0.24414893617021272cm);
\draw(9.08,-1.7) ellipse (0.22089843750000002cm and 0.24414893617021272cm);
\draw (8.68,-0.5)-- (9.02,-1.54);
\draw (10.34,-1.78)-- (9.08,-1.8);
\draw (10.74,-0.5)-- (10.42,-1.58);
\draw [->,line width=1.2000000000000002pt] (-2.88,-0.88) -- (-2.0,-0.88);
\begin{scriptsize}
\draw [fill=black] (1.92,-1.56) circle (1.5pt);
\draw [fill=black] (3.2,-1.56) circle (1.5pt);
\draw [fill=uuuuuu] (1.5244582472000674,-0.34264765914220313) circle (1.5pt);
\draw [fill=uuuuuu] (3.595541752799933,-0.3426476591422037) circle (1.5pt);
\draw [fill=uuuuuu] (1.5244582472000674,-0.34264765914220313) circle (1.5pt);
\draw [fill=black] (-5.24,-0.34) circle (1.5pt);
\draw [fill=black] (-4.86,-1.56) circle (1.5pt);
\draw [fill=black] (-3.54,-1.54) circle (1.5pt);
\draw [fill=black] (-1.88,-0.34) circle (1.5pt);
\draw [fill=black] (-1.54,-1.38) circle (1.5pt);
\draw [fill=black] (-0.22,-1.62) circle (1.5pt);
\draw [fill=black] (-1.48,-1.64) circle (1.5pt);
\draw [fill=black] (5.08,-0.5) circle (1.5pt);
\draw [fill=black] (5.42,-1.54) circle (1.5pt);
\draw [fill=black] (6.74,-1.78) circle (1.5pt);
\draw [fill=black] (5.48,-1.8) circle (1.5pt);
\draw [fill=black] (7.14,-0.52) circle (1.5pt);
\draw [fill=black] (8.68,-0.5) circle (1.5pt);
\draw [fill=black] (9.02,-1.54) circle (1.5pt);
\draw [fill=black] (10.34,-1.78) circle (1.5pt);
\draw [fill=black] (9.08,-1.8) circle (1.5pt);
\draw [fill=black] (10.74,-0.5) circle (1.5pt);
\draw [fill=black] (10.42,-1.58) circle (1.5pt);
\end{scriptsize}
\end{tikzpicture} 
\vspace{4pt}
\definecolor{uuuuuu}{rgb}{0.26666666666666666,0.26666666666666666,0.26666666666666666}
\begin{tikzpicture}[line cap=round,line join=round,>=triangle 45,x=0.7558593750000001cm,y=0.8138297872340424cm]
\clip(-5.600000000000002,-2.8599999999999888) rectangle (14.880000000000006,0.9000000000000069);
\draw(-4.2,0.42971746379216136) ellipse (0.22675781250000002cm and 0.24414893617021272cm);
\draw(-5.235541752799932,-0.3226476591422037) ellipse (0.22675781250000002cm and 0.24414893617021272cm);
\draw(-3.1644582472000677,-0.3226476591422038) ellipse (0.22675781250000002cm and 0.24414893617021272cm);
\draw(-3.56,-1.54) ellipse (0.22675781250000002cm and 0.24414893617021272cm);
\draw(-4.84,-1.54) ellipse (0.22675781250000002cm and 0.24414893617021272cm);
\draw(-0.8399999999999999,0.4297174637921619) ellipse (0.22675781250000002cm and 0.24414893617021272cm);
\draw(-1.8755417527999323,-0.32264765914220334) ellipse (0.22675781250000002cm and 0.24414893617021272cm);
\draw(0.1955417527999329,-0.32264765914220384) ellipse (0.22675781250000002cm and 0.24414893617021272cm);
\draw(-0.2,-1.54) ellipse (0.22675781250000002cm and 0.24414893617021272cm);
\draw(-1.48,-1.54) ellipse (0.22675781250000002cm and 0.24414893617021272cm);
\draw(2.5600000000000005,0.40971746379216223) ellipse (0.22675781250000002cm and 0.24414893617021272cm);
\draw(1.5244582472000674,-0.34264765914220313) ellipse (0.22675781250000002cm and 0.24414893617021272cm);
\draw(3.595541752799933,-0.3426476591422037) ellipse (0.22675781250000002cm and 0.24414893617021272cm);
\draw(3.2,-1.56) ellipse (0.22675781250000002cm and 0.24414893617021272cm);
\draw(1.92,-1.56) ellipse (0.22675781250000002cm and 0.24414893617021272cm);
\draw (-4.600000000000001,-1.9999999999999896) node[anchor=north west] {$R_6$};
\draw (-1.2800000000000005,-1.9999999999999896) node[anchor=north west] {$S_4$};
\draw (2.000000000000001,-1.9999999999999896) node[anchor=north west] {$S_5$};
\draw (5.500000000000002,-1.9799999999999898) node[anchor=north west] {$R_9$};
\draw (9.120000000000003,-1.9799999999999898) node[anchor=north west] {$S_6$};
\draw [->,line width=1.2000000000000002pt] (7.64,-1.04) -- (8.52,-1.04);
\draw (0.6600000000000003,-0.5999999999999913) node[anchor=north west] {,};
\draw [line width=2.0pt,dash pattern=on 4pt off 4pt] (4.28,0.5)-- (4.28,-1.9);
\draw (-5.24,-0.34)-- (-4.86,-1.56);
\draw (-4.86,-1.56)-- (-3.54,-1.54);
\draw (-1.88,-0.34)-- (-1.54,-1.38);
\draw (-0.22,-1.62)-- (-1.48,-1.64);
\draw(6.12,0.2697174637921611) ellipse (0.22675781250000002cm and 0.24414893617021272cm);
\draw(5.084458247200068,-0.48264765914220387) ellipse (0.22675781250000002cm and 0.24414893617021272cm);
\draw(7.155541752799932,-0.48264765914220437) ellipse (0.22675781250000002cm and 0.24414893617021272cm);
\draw(6.76,-1.7) ellipse (0.22675781250000002cm and 0.24414893617021272cm);
\draw(5.48,-1.7) ellipse (0.22675781250000002cm and 0.24414893617021272cm);
\draw (5.1,-0.5)-- (5.48,-1.7);
\draw (6.74,-1.72)-- (5.48,-1.72);
\draw (6.74,-1.72)-- (7.14,-0.52);
\draw(9.719999999999999,0.26971746379216044) ellipse (0.22675781250000002cm and 0.24414893617021272cm);
\draw(8.684458247200068,-0.48264765914220464) ellipse (0.22675781250000002cm and 0.24414893617021272cm);
\draw(10.75554175279993,-0.4826476591422038) ellipse (0.22675781250000002cm and 0.24414893617021272cm);
\draw(10.36,-1.7) ellipse (0.22675781250000002cm and 0.24414893617021272cm);
\draw(9.08,-1.7) ellipse (0.22675781250000002cm and 0.24414893617021272cm);
\draw (8.7,-0.52)-- (9.12,-1.76);
\draw (10.34,-1.78)-- (9.12,-1.76);
\draw (10.74,-0.5)-- (10.42,-1.58);
\draw (-3.22,-0.3)-- (-4.86,-1.56);
\draw (0.14,-0.36)-- (-1.32,-1.46);
\draw (10.34,-1.78)-- (8.7,-0.52);
\draw (5.1,-0.5)-- (6.74,-1.72);
\draw [->,line width=1.2000000000000002pt] (-2.98,-0.96) -- (-2.1,-0.96);
\draw (3.2,-1.56)-- (1.84,-1.64);
\draw (1.84,-1.64)-- (1.5244582472000674,-0.34264765914220313);
\draw (2.02,-1.5)-- (3.5,-0.38);
\begin{scriptsize}
\draw [fill=black] (3.2,-1.56) circle (1.5pt);
\draw [fill=uuuuuu] (1.5244582472000674,-0.34264765914220313) circle (1.5pt);
\draw [fill=uuuuuu] (1.5244582472000674,-0.34264765914220313) circle (1.5pt);
\draw [fill=black] (-5.24,-0.34) circle (1.5pt);
\draw [fill=black] (-4.86,-1.56) circle (1.5pt);
\draw [fill=black] (-3.54,-1.54) circle (1.5pt);
\draw [fill=black] (-1.88,-0.34) circle (1.5pt);
\draw [fill=black] (-1.54,-1.38) circle (1.5pt);
\draw [fill=black] (-0.22,-1.62) circle (1.5pt);
\draw [fill=black] (-1.48,-1.64) circle (1.5pt);
\draw [fill=black] (5.1,-0.5) circle (1.5pt);
\draw [fill=black] (5.48,-1.7) circle (1.5pt);
\draw [fill=black] (6.74,-1.72) circle (1.5pt);
\draw [fill=black] (5.48,-1.72) circle (1.5pt);
\draw [fill=black] (7.14,-0.52) circle (1.5pt);
\draw [fill=black] (8.7,-0.52) circle (1.5pt);
\draw [fill=black] (9.12,-1.76) circle (1.5pt);
\draw [fill=black] (10.34,-1.78) circle (1.5pt);
\draw [fill=black] (9.12,-1.76) circle (1.5pt);
\draw [fill=black] (10.74,-0.5) circle (1.5pt);
\draw [fill=black] (10.42,-1.58) circle (1.5pt);
\draw [fill=black] (-3.22,-0.3) circle (1.5pt);
\draw [fill=black] (0.14,-0.36) circle (1.5pt);
\draw [fill=black] (-1.32,-1.46) circle (1.5pt);
\draw [fill=black] (3.5,-0.38) circle (1.5pt);
\draw [fill=black] (2.02,-1.5) circle (1.5pt);
\draw [fill=black] (1.84,-1.64) circle (1.5pt);
\end{scriptsize}
\end{tikzpicture} 
\vspace{4pt}
\definecolor{uuuuuu}{rgb}{0.26666666666666666,0.26666666666666666,0.26666666666666666}
\begin{tikzpicture}[line cap=round,line join=round,>=triangle 45,x=0.7363281250000001cm,y=0.8138297872340424cm]
\clip(-5.600000000000002,-2.8599999999999888) rectangle (14.880000000000006,0.9000000000000069);
\draw(-4.2,0.42971746379216136) ellipse (0.22089843750000002cm and 0.24414893617021272cm);
\draw(-5.235541752799932,-0.3226476591422037) ellipse (0.22089843750000002cm and 0.24414893617021272cm);
\draw(-3.1644582472000677,-0.3226476591422038) ellipse (0.22089843750000002cm and 0.24414893617021272cm);
\draw(-3.56,-1.54) ellipse (0.22089843750000002cm and 0.24414893617021272cm);
\draw(-4.84,-1.54) ellipse (0.22089843750000002cm and 0.24414893617021272cm);
\draw(-0.8399999999999999,0.4297174637921619) ellipse (0.22089843750000002cm and 0.24414893617021272cm);
\draw(-1.8755417527999323,-0.32264765914220334) ellipse (0.22089843750000002cm and 0.24414893617021272cm);
\draw(0.1955417527999329,-0.32264765914220384) ellipse (0.22089843750000002cm and 0.24414893617021272cm);
\draw(-0.2,-1.54) ellipse (0.22089843750000002cm and 0.24414893617021272cm);
\draw(-1.48,-1.54) ellipse (0.22089843750000002cm and 0.24414893617021272cm);
\draw(2.5600000000000005,0.40971746379216223) ellipse (0.22089843750000002cm and 0.24414893617021272cm);
\draw(1.5244582472000674,-0.34264765914220313) ellipse (0.22089843750000002cm and 0.24414893617021272cm);
\draw(3.595541752799933,-0.3426476591422037) ellipse (0.22089843750000002cm and 0.24414893617021272cm);
\draw(3.2,-1.56) ellipse (0.22089843750000002cm and 0.24414893617021272cm);
\draw(1.92,-1.56) ellipse (0.22089843750000002cm and 0.24414893617021272cm);
\draw (-4.600000000000001,-1.9999999999999896) node[anchor=north west] {$R_{12}$};
\draw (-1.2800000000000005,-1.9999999999999896) node[anchor=north west] {$S_7$};
\draw (2.000000000000001,-1.9999999999999896) node[anchor=north west] {$R_{14}$};
\draw (5.500000000000002,-1.9799999999999898) node[anchor=north west] {$S_8$};
\draw (9.120000000000003,-1.9799999999999898) node[anchor=north west] {$S_9$};
\draw [->,line width=1.2000000000000002pt] (3.9,-0.92) -- (4.78,-0.92);
\draw (7.720000000000002,-0.6599999999999913) node[anchor=north west] {,};
\draw [line width=1.2000000000000002pt,dash pattern=on 4pt off 4pt] (1.04,0.64)-- (1.04,-1.76);
\draw (-4.86,-1.56)-- (-3.54,-1.54);
\draw (-0.22,-1.62)-- (-1.48,-1.64);
\draw (1.5244582472000674,-0.34264765914220313)-- (1.92,-1.56);
\draw [->,line width=1.2000000000000002pt] (-2.98,-0.96) -- (-2.1,-0.96);
\draw (-5.24,-0.34)-- (-4.2,0.42);
\draw (-4.2,0.42)-- (-3.22,-0.3);
\draw (0.18,-0.36)-- (-0.7,0.32);
\draw (-0.94,0.34)-- (-1.88,-0.34);
\draw (1.5244582472000674,-0.34264765914220313)-- (2.54,0.44);
\draw (2.54,0.44)-- (3.54,-0.36);
\draw (3.54,-0.36)-- (3.2,-1.56);
\draw(6.160000000000001,0.40971746379216245) ellipse (0.22089843750000002cm and 0.24414893617021272cm);
\draw(5.124458247200067,-0.3426476591422027) ellipse (0.22089843750000002cm and 0.24414893617021272cm);
\draw(7.195541752799933,-0.34264765914220396) ellipse (0.22089843750000002cm and 0.24414893617021272cm);
\draw(6.8,-1.56) ellipse (0.22089843750000002cm and 0.24414893617021272cm);
\draw(5.52,-1.56) ellipse (0.22089843750000002cm and 0.24414893617021272cm);
\draw (6.14,0.44)-- (7.14,-0.36);
\draw (7.14,-0.36)-- (6.8,-1.56);
\draw (6.14,0.44)-- (5.2,-0.24);
\draw (5.52,-1.58)-- (5.08,-0.48);
\draw(9.719999999999999,0.5097174637921604) ellipse (0.22089843750000002cm and 0.24414893617021272cm);
\draw(8.684458247200068,-0.24264765914220465) ellipse (0.22089843750000002cm and 0.24414893617021272cm);
\draw(10.75554175279993,-0.24264765914220382) ellipse (0.22089843750000002cm and 0.24414893617021272cm);
\draw(10.36,-1.46) ellipse (0.22089843750000002cm and 0.24414893617021272cm);
\draw(9.08,-1.46) ellipse (0.22089843750000002cm and 0.24414893617021272cm);
\draw (8.684458247200068,-0.24264765914220465)-- (9.08,-1.46);
\draw (10.7,-0.28)-- (10.36,-1.46);
\draw (9.58,0.48)-- (8.684458247200068,-0.24264765914220465);
\draw (10.7,-0.28)-- (9.86,0.48);
\draw (10.820000000000004,-0.6399999999999912) node[anchor=north west] {ctd. in the next };
\draw (10.860000000000005,-1.0999999999999908) node[anchor=north west] {line};
\begin{scriptsize}
\draw [fill=black] (1.92,-1.56) circle (1.5pt);
\draw [fill=black] (3.2,-1.56) circle (1.5pt);
\draw [fill=uuuuuu] (1.5244582472000674,-0.34264765914220313) circle (1.5pt);
\draw [fill=uuuuuu] (1.5244582472000674,-0.34264765914220313) circle (1.5pt);
\draw [fill=black] (-5.24,-0.34) circle (1.5pt);
\draw [fill=black] (-4.86,-1.56) circle (1.5pt);
\draw [fill=black] (-3.54,-1.54) circle (1.5pt);
\draw [fill=black] (-1.88,-0.34) circle (1.5pt);
\draw [fill=black] (-0.22,-1.62) circle (1.5pt);
\draw [fill=black] (-1.48,-1.64) circle (1.5pt);
\draw [fill=black] (-3.22,-0.3) circle (1.5pt);
\draw [fill=black] (0.18,-0.36) circle (1.5pt);
\draw [fill=black] (1.92,-1.58) circle (1.5pt);
\draw [fill=black] (3.54,-0.36) circle (1.5pt);
\draw [fill=black] (-4.2,0.42) circle (1.5pt);
\draw [fill=black] (-0.7,0.32) circle (1.5pt);
\draw [fill=black] (-0.94,0.34) circle (1.5pt);
\draw [fill=black] (2.54,0.44) circle (1.5pt);
\draw [fill=black] (5.52,-1.56) circle (1.5pt);
\draw [fill=black] (6.8,-1.56) circle (1.5pt);
\draw [fill=black] (5.52,-1.58) circle (1.5pt);
\draw [fill=black] (7.14,-0.36) circle (1.5pt);
\draw [fill=black] (6.14,0.44) circle (1.5pt);
\draw [fill=black] (5.2,-0.24) circle (1.5pt);
\draw [fill=black] (5.08,-0.48) circle (1.5pt);
\draw [fill=black] (9.08,-1.46) circle (1.5pt);
\draw [fill=black] (10.36,-1.46) circle (1.5pt);
\draw [fill=uuuuuu] (8.684458247200068,-0.24264765914220465) circle (1.5pt);
\draw [fill=uuuuuu] (8.684458247200068,-0.24264765914220465) circle (1.5pt);
\draw [fill=black] (9.08,-1.48) circle (1.5pt);
\draw [fill=black] (10.7,-0.28) circle (1.5pt);
\draw [fill=black] (9.58,0.48) circle (1.5pt);
\draw [fill=black] (9.86,0.48) circle (1.5pt);
\end{scriptsize}
\end{tikzpicture} 
\vspace{4pt}
\definecolor{uuuuuu}{rgb}{0.26666666666666666,0.26666666666666666,0.26666666666666666}
\begin{tikzpicture}[line cap=round,line join=round,>=triangle 45,x=0.7286821705426357cm,y=0.8138297872340424cm]
\clip(-5.600000000000003,-2.8599999999999888) rectangle (15.040000000000006,0.9000000000000069);
\draw(-4.16,0.42971746379216136) ellipse (0.2186046511627907cm and 0.24414893617021272cm);
\draw(-5.195541752799932,-0.3226476591422037) ellipse (0.2186046511627907cm and 0.24414893617021272cm);
\draw(-3.1244582472000677,-0.3226476591422038) ellipse (0.2186046511627907cm and 0.24414893617021272cm);
\draw(-3.52,-1.54) ellipse (0.2186046511627907cm and 0.24414893617021272cm);
\draw(-4.8,-1.54) ellipse (0.2186046511627907cm and 0.24414893617021272cm);
\draw (-4.720000000000003,-1.9799999999999898) node[anchor=north west] {$R_{14}$};
\draw (-1.2200000000000015,-1.9599999999999898) node[anchor=north west] {$S_{10}$};
\draw (2.4,-1.9599999999999898) node[anchor=north west] {$S_{11}$};
\draw [->,line width=1.2000000000000002pt] (-2.86,-0.9) -- (-1.98,-0.9);
\draw (0.9599999999999992,-0.6999999999999912) node[anchor=north west] {,};
\draw (-5.195541752799932,-0.3226476591422037)-- (-4.8,-1.54);
\draw (-5.195541752799932,-0.3226476591422037)-- (-4.18,0.46);
\draw (-4.18,0.46)-- (-3.18,-0.34);
\draw (-3.18,-0.34)-- (-3.52,-1.54);
\draw(-0.5599999999999998,0.4297174637921619) ellipse (0.2186046511627907cm and 0.24414893617021272cm);
\draw(-1.5955417527999325,-0.32264765914220334) ellipse (0.2186046511627907cm and 0.24414893617021272cm);
\draw(0.47554175279993294,-0.32264765914220384) ellipse (0.2186046511627907cm and 0.24414893617021272cm);
\draw(0.08,-1.54) ellipse (0.2186046511627907cm and 0.24414893617021272cm);
\draw(-1.2,-1.54) ellipse (0.2186046511627907cm and 0.24414893617021272cm);
\draw (-0.58,0.46)-- (0.48,-0.2);
\draw (-0.58,0.46)-- (-1.52,-0.22);
\draw (-1.2,-1.56)-- (-1.64,-0.46);
\draw (0.48,-0.44)-- (0.08,-1.54);
\draw(3.0200000000000005,0.30971746379216236) ellipse (0.2186046511627907cm and 0.24414893617021272cm);
\draw(1.9844582472000674,-0.442647659142203) ellipse (0.2186046511627907cm and 0.24414893617021272cm);
\draw(4.0555417527999325,-0.44264765914220355) ellipse (0.2186046511627907cm and 0.24414893617021272cm);
\draw(3.66,-1.66) ellipse (0.2186046511627907cm and 0.24414893617021272cm);
\draw(2.38,-1.66) ellipse (0.2186046511627907cm and 0.24414893617021272cm);
\draw (3.16,0.32)-- (4.06,-0.32);
\draw (2.38,-1.68)-- (1.94,-0.58);
\draw (4.06,-0.56)-- (3.66,-1.66);
\draw (2.06,-0.34)-- (2.92,0.34);
\draw(6.29,0.31894062842040893) ellipse (0.2186046511627907cm and 0.24414893617021272cm);
\draw(5.270638587087567,-0.42166878946810676) ellipse (0.2186046511627907cm and 0.24414893617021272cm);
\draw(7.309361412912434,-0.4216687894681072) ellipse (0.2186046511627907cm and 0.24414893617021272cm);
\draw(6.92,-1.62) ellipse (0.2186046511627907cm and 0.24414893617021272cm);
\draw(5.66,-1.62) ellipse (0.2186046511627907cm and 0.24414893617021272cm);
\draw (9.240000000000002,-2.0599999999999894) node[anchor=north west] {$S_{12}$};
\draw (12.860000000000005,-2.03999999999999) node[anchor=north west] {$S_{13}$};
\draw [->,line width=1.2000000000000002pt] (7.6,-0.98) -- (8.48,-0.98);
\draw (5.270638587087567,-0.42166878946810676)-- (5.66,-1.62);
\draw (7.28,-0.42)-- (6.92,-1.62);
\draw(9.899999999999999,0.3497174637921603) ellipse (0.2186046511627907cm and 0.24414893617021272cm);
\draw(8.864458247200067,-0.4026476591422048) ellipse (0.2186046511627907cm and 0.24414893617021272cm);
\draw(10.93554175279993,-0.40264765914220396) ellipse (0.2186046511627907cm and 0.24414893617021272cm);
\draw(10.54,-1.62) ellipse (0.2186046511627907cm and 0.24414893617021272cm);
\draw(9.26,-1.62) ellipse (0.2186046511627907cm and 0.24414893617021272cm);
\draw (9.88,0.38)-- (9.26,-1.64);
\draw (9.26,-1.64)-- (8.82,-0.54);
\draw (10.9,-0.38)-- (10.66,-1.6);
\draw [line width=1.2000000000000002pt,dash pattern=on 4pt off 4pt] (4.56,0.64)-- (4.62,-1.88);
\draw (11.440000000000005,-0.6399999999999912) node[anchor=north west] {,};
\draw (5.940000000000001,-1.9999999999999896) node[anchor=north west] {$R_{16}$};
\draw (5.64,-1.64)-- (6.92,-1.62);
\draw (6.28,0.38)-- (5.64,-1.64);
\draw (9.26,-1.64)-- (10.36,-1.62);
\draw(13.480000000000002,0.28971746379216157) ellipse (0.2186046511627907cm and 0.24414893617021272cm);
\draw(12.444458247200068,-0.46264765914220296) ellipse (0.2186046511627907cm and 0.24414893617021272cm);
\draw(14.515541752799933,-0.46264765914220485) ellipse (0.2186046511627907cm and 0.24414893617021272cm);
\draw(14.12,-1.68) ellipse (0.2186046511627907cm and 0.24414893617021272cm);
\draw(12.84,-1.68) ellipse (0.2186046511627907cm and 0.24414893617021272cm);
\draw (14.48,-0.44)-- (14.12,-1.68);
\draw (12.46,-0.46)-- (12.7,-1.7);
\draw (13.46,0.32)-- (12.96,-1.7);
\draw (12.96,-1.7)-- (14.12,-1.68);
\begin{scriptsize}
\draw [fill=black] (-4.8,-1.54) circle (1.5pt);
\draw [fill=black] (-3.52,-1.54) circle (1.5pt);
\draw [fill=uuuuuu] (-5.195541752799932,-0.3226476591422037) circle (1.5pt);
\draw [fill=uuuuuu] (-5.195541752799932,-0.3226476591422037) circle (1.5pt);
\draw [fill=black] (-4.8,-1.56) circle (1.5pt);
\draw [fill=black] (-3.18,-0.34) circle (1.5pt);
\draw [fill=black] (-4.18,0.46) circle (1.5pt);
\draw [fill=black] (-1.2,-1.54) circle (1.5pt);
\draw [fill=black] (-1.2,-1.56) circle (1.5pt);
\draw [fill=black] (0.48,-0.2) circle (1.5pt);
\draw [fill=black] (-0.58,0.46) circle (1.5pt);
\draw [fill=black] (-1.52,-0.22) circle (1.5pt);
\draw [fill=black] (-1.64,-0.46) circle (1.5pt);
\draw [fill=black] (0.08,-1.54) circle (1.5pt);
\draw [fill=black] (0.48,-0.44) circle (1.5pt);
\draw [fill=black] (2.38,-1.66) circle (1.5pt);
\draw [fill=black] (2.38,-1.68) circle (1.5pt);
\draw [fill=black] (4.06,-0.32) circle (1.5pt);
\draw [fill=black] (3.16,0.32) circle (1.5pt);
\draw [fill=black] (2.06,-0.34) circle (1.5pt);
\draw [fill=black] (1.94,-0.58) circle (1.5pt);
\draw [fill=black] (3.66,-1.66) circle (1.5pt);
\draw [fill=black] (4.06,-0.56) circle (1.5pt);
\draw [fill=black] (2.92,0.34) circle (1.5pt);
\draw [fill=black] (5.66,-1.62) circle (1.5pt);
\draw [fill=black] (6.92,-1.62) circle (1.5pt);
\draw [fill=uuuuuu] (5.270638587087567,-0.42166878946810676) circle (1.5pt);
\draw [fill=uuuuuu] (5.270638587087567,-0.42166878946810676) circle (1.5pt);
\draw [fill=black] (5.64,-1.64) circle (1.5pt);
\draw [fill=black] (7.28,-0.42) circle (1.5pt);
\draw [fill=black] (6.28,0.38) circle (1.5pt);
\draw [fill=black] (9.26,-1.62) circle (1.5pt);
\draw [fill=black] (9.26,-1.64) circle (1.5pt);
\draw [fill=black] (9.88,0.38) circle (1.5pt);
\draw [fill=black] (8.82,-0.54) circle (1.5pt);
\draw [fill=black] (10.66,-1.6) circle (1.5pt);
\draw [fill=black] (10.9,-0.38) circle (1.5pt);
\draw [fill=black] (10.36,-1.62) circle (1.5pt);
\draw [fill=black] (13.46,0.32) circle (1.5pt);
\draw [fill=black] (12.46,-0.46) circle (1.5pt);
\draw [fill=black] (14.12,-1.68) circle (1.5pt);
\draw [fill=black] (14.48,-0.44) circle (1.5pt);
\draw [fill=black] (12.96,-1.7) circle (1.5pt);
\draw [fill=black] (12.7,-1.7) circle (1.5pt);
\end{scriptsize}
\end{tikzpicture} 
\vspace{4pt}
\definecolor{uuuuuu}{rgb}{0.26666666666666666,0.26666666666666666,0.26666666666666666}
\begin{tikzpicture}[line cap=round,line join=round,>=triangle 45,x=0.7286821705426357cm,y=0.8138297872340424cm]
\clip(-5.600000000000003,-2.8599999999999888) rectangle (15.040000000000006,0.9000000000000069);
\draw (9.240000000000002,-2.0599999999999894) node[anchor=north west] {$S_{17}$};
\draw (5.960000000000001,-1.9999999999999896) node[anchor=north west] {$S_{16}$};
\draw(-4.13,0.33894062842040895) ellipse (0.2186046511627907cm and 0.24414893617021272cm);
\draw(-5.149361412912433,-0.40166878946810675) ellipse (0.2186046511627907cm and 0.24414893617021272cm);
\draw(-3.1106385870875664,-0.4016687894681072) ellipse (0.2186046511627907cm and 0.24414893617021272cm);
\draw(-3.5,-1.6) ellipse (0.2186046511627907cm and 0.24414893617021272cm);
\draw(-4.76,-1.6) ellipse (0.2186046511627907cm and 0.24414893617021272cm);
\draw [->,line width=1.2000000000000002pt] (-2.82,-0.96) -- (-1.94,-0.96);
\draw (-5.149361412912433,-0.40166878946810675)-- (-4.76,-1.6);
\draw (-3.14,-0.4)-- (-3.5,-1.6);
\draw(-0.5199999999999998,0.3697174637921615) ellipse (0.2186046511627907cm and 0.24414893617021272cm);
\draw(-1.5555417527999325,-0.38264765914220356) ellipse (0.2186046511627907cm and 0.24414893617021272cm);
\draw(0.5155417527999324,-0.38264765914220394) ellipse (0.2186046511627907cm and 0.24414893617021272cm);
\draw(0.12,-1.6) ellipse (0.2186046511627907cm and 0.24414893617021272cm);
\draw(-1.16,-1.6) ellipse (0.2186046511627907cm and 0.24414893617021272cm);
\draw (0.48,-0.36)-- (0.24,-1.58);
\draw (-4.78,-1.62)-- (-3.5,-1.6);
\draw (-4.14,0.4)-- (-4.78,-1.62);
\draw(3.0600000000000005,0.30971746379216236) ellipse (0.2186046511627907cm and 0.24414893617021272cm);
\draw(2.0244582472000676,-0.442647659142203) ellipse (0.2186046511627907cm and 0.24414893617021272cm);
\draw(4.095541752799933,-0.44264765914220355) ellipse (0.2186046511627907cm and 0.24414893617021272cm);
\draw(3.7,-1.66) ellipse (0.2186046511627907cm and 0.24414893617021272cm);
\draw(2.42,-1.66) ellipse (0.2186046511627907cm and 0.24414893617021272cm);
\draw (4.06,-0.42)-- (3.7,-1.66);
\draw (2.04,-0.44)-- (2.28,-1.68);
\draw (-1.0400000000000016,-1.9999999999999896) node[anchor=north west] {$S_{14}$};
\draw (2.5799999999999996,-1.9799999999999898) node[anchor=north west] {$S_{15}$};
\draw (-4.3400000000000025,-1.9399999999999897) node[anchor=north west] {$R_{16}$};
\draw (0.9799999999999992,-0.5599999999999914) node[anchor=north west] {,};
\draw (4.4,-0.5999999999999913) node[anchor=north west] {,};
\draw (7.840000000000002,-0.6399999999999912) node[anchor=north west] {,};
\draw (0.24,-1.58)-- (-0.98,-1.6);
\draw (3.04,0.32)-- (2.28,-1.68);
\draw (2.54,-1.7)-- (3.5,-1.66);
\draw(6.38,0.2697174637921611) ellipse (0.2186046511627907cm and 0.24414893617021272cm);
\draw(5.344458247200068,-0.48264765914220387) ellipse (0.2186046511627907cm and 0.24414893617021272cm);
\draw(7.02,-1.7) ellipse (0.2186046511627907cm and 0.24414893617021272cm);
\draw(5.74,-1.7) ellipse (0.2186046511627907cm and 0.24414893617021272cm);
\draw (7.38,-0.46)-- (7.02,-1.7);
\draw (5.36,-0.48)-- (5.6,-1.72);
\draw (5.86,-1.74)-- (6.82,-1.7);
\draw(9.74,0.2594349275843163) ellipse (0.2186046511627907cm and 0.24414893617021272cm);
\draw(8.704458247200067,-0.4929301953500495) ellipse (0.2186046511627907cm and 0.24414893617021272cm);
\draw(10.38,-1.7102825362078473) ellipse (0.2186046511627907cm and 0.24414893617021272cm);
\draw(9.1,-1.7102825362078473) ellipse (0.2186046511627907cm and 0.24414893617021272cm);
\draw (10.74,-0.4702825362078473)-- (10.38,-1.7102825362078473);
\draw (8.72,-0.49028253620784734)-- (8.96,-1.7302825362078473);
\draw (9.22,-1.7502825362078474)-- (10.18,-1.7102825362078473);
\draw (6.36,0.28)-- (5.86,-1.74);
\draw (9.72,0.22)-- (9.12,-1.58);
\draw (-1.46,-0.36)-- (-1.22,-1.6);
\draw (-0.45999999999999996,0.42000000000001747)-- (-1.22,-1.6);
\begin{scriptsize}
\draw [fill=black] (-4.76,-1.6) circle (1.5pt);
\draw [fill=black] (-3.5,-1.6) circle (1.5pt);
\draw [fill=uuuuuu] (-5.149361412912433,-0.40166878946810675) circle (1.5pt);
\draw [fill=uuuuuu] (-5.149361412912433,-0.40166878946810675) circle (1.5pt);
\draw [fill=black] (-4.78,-1.62) circle (1.5pt);
\draw [fill=black] (-3.14,-0.4) circle (1.5pt);
\draw [fill=black] (-4.14,0.4) circle (1.5pt);
\draw [fill=black] (-0.5,0.36) circle (1.5pt);
\draw [fill=black] (0.24,-1.58) circle (1.5pt);
\draw [fill=black] (0.48,-0.36) circle (1.5pt);
\draw [fill=black] (3.04,0.32) circle (1.5pt);
\draw [fill=black] (2.04,-0.44) circle (1.5pt);
\draw [fill=black] (3.7,-1.66) circle (1.5pt);
\draw [fill=black] (4.06,-0.42) circle (1.5pt);
\draw [fill=black] (2.54,-1.7) circle (1.5pt);
\draw [fill=black] (2.28,-1.68) circle (1.5pt);
\draw [fill=black] (-0.98,-1.6) circle (1.5pt);
\draw [fill=black] (3.5,-1.66) circle (1.5pt);
\draw [fill=black] (6.36,0.28) circle (1.5pt);
\draw [fill=black] (5.36,-0.48) circle (1.5pt);
\draw [fill=black] (7.02,-1.7) circle (1.5pt);
\draw [fill=black] (7.38,-0.46) circle (1.5pt);
\draw [fill=black] (5.86,-1.74) circle (1.5pt);
\draw [fill=black] (5.6,-1.72) circle (1.5pt);
\draw [fill=black] (6.82,-1.7) circle (1.5pt);
\draw [fill=black] (8.72,-0.49028253620784734) circle (1.5pt);
\draw [fill=black] (10.38,-1.7102825362078473) circle (1.5pt);
\draw [fill=black] (10.74,-0.4702825362078473) circle (1.5pt);
\draw [fill=black] (9.22,-1.7502825362078474) circle (1.5pt);
\draw [fill=black] (8.96,-1.7302825362078473) circle (1.5pt);
\draw [fill=black] (10.18,-1.7102825362078473) circle (1.5pt);
\draw [fill=black] (9.72,0.22) circle (1.5pt);
\draw [fill=black] (9.12,-1.58) circle (1.5pt);
\draw [fill=black] (-1.46,-0.36) circle (1.5pt);
\draw [fill=black] (-1.22,-1.6) circle (1.5pt);
\end{scriptsize}
\end{tikzpicture} 
\vspace{4pt}
\definecolor{uuuuuu}{rgb}{0.26666666666666666,0.26666666666666666,0.26666666666666666}
\begin{tikzpicture}[line cap=round,line join=round,>=triangle 45,x=0.7558593750000001cm,y=0.8138297872340424cm]
\clip(-5.600000000000002,-2.8599999999999888) rectangle (14.880000000000006,0.9000000000000069);
\draw(-4.16,0.42971746379216136) ellipse (0.22675781250000002cm and 0.24414893617021272cm);
\draw(-5.195541752799932,-0.3226476591422037) ellipse (0.22675781250000002cm and 0.24414893617021272cm);
\draw(-3.1244582472000677,-0.3226476591422038) ellipse (0.22675781250000002cm and 0.24414893617021272cm);
\draw(-3.52,-1.54) ellipse (0.22675781250000002cm and 0.24414893617021272cm);
\draw(-4.8,-1.54) ellipse (0.22675781250000002cm and 0.24414893617021272cm);
\draw (-4.7200000000000015,-1.9799999999999898) node[anchor=north west] {$R_{15}$};
\draw (-1.2200000000000004,-1.9599999999999898) node[anchor=north west] {$S_{18}$};
\draw (2.400000000000001,-1.9599999999999898) node[anchor=north west] {$S_{19}$};
\draw [->,line width=1.2000000000000002pt] (-2.86,-0.9) -- (-1.98,-0.9);
\draw (0.9000000000000004,-0.8199999999999911) node[anchor=north west] {,};
\draw (-5.195541752799932,-0.3226476591422037)-- (-4.18,0.46);
\draw (-4.18,0.46)-- (-3.18,-0.34);
\draw(3.0200000000000005,0.30971746379216236) ellipse (0.22675781250000002cm and 0.24414893617021272cm);
\draw(1.9844582472000674,-0.442647659142203) ellipse (0.22675781250000002cm and 0.24414893617021272cm);
\draw(4.0555417527999325,-0.44264765914220355) ellipse (0.22675781250000002cm and 0.24414893617021272cm);
\draw(3.66,-1.66) ellipse (0.22675781250000002cm and 0.24414893617021272cm);
\draw(2.38,-1.66) ellipse (0.22675781250000002cm and 0.24414893617021272cm);
\draw (3.16,0.4)-- (4.06,-0.32);
\draw (2.38,-1.68)-- (2.92,0.4);
\draw (3.16,0.36)-- (3.66,-1.66);
\draw (2.06,-0.34)-- (2.9,0.42);
\draw (-4.18,0.46)-- (-4.8,-1.56);
\draw (-4.18,0.46)-- (-3.52,-1.54);
\draw(-0.5399999999999998,0.34971746379216184) ellipse (0.22675781250000002cm and 0.24414893617021272cm);
\draw(-1.5755417527999325,-0.4026476591422034) ellipse (0.22675781250000002cm and 0.24414893617021272cm);
\draw(0.49554175279993296,-0.4026476591422039) ellipse (0.22675781250000002cm and 0.24414893617021272cm);
\draw(0.1,-1.62) ellipse (0.22675781250000002cm and 0.24414893617021272cm);
\draw(-1.18,-1.62) ellipse (0.22675781250000002cm and 0.24414893617021272cm);
\draw (-0.42,0.42)-- (0.5,-0.28);
\draw (-1.18,-1.64)-- (-0.42,0.38);
\draw (-0.4,0.44)-- (0.1,-1.62);
\draw (-1.5,-0.3)-- (-0.68,0.38);
\draw (5.700000000000002,-1.9799999999999898) node[anchor=north west] {$S_{20}$};
\draw (9.320000000000004,-1.9799999999999898) node[anchor=north west] {$S_{21}$};
\draw(9.940000000000003,0.28971746379216157) ellipse (0.22675781250000002cm and 0.24414893617021272cm);
\draw(8.904458247200068,-0.46264765914220296) ellipse (0.22675781250000002cm and 0.24414893617021272cm);
\draw(10.975541752799934,-0.46264765914220485) ellipse (0.22675781250000002cm and 0.24414893617021272cm);
\draw(10.58,-1.68) ellipse (0.22675781250000002cm and 0.24414893617021272cm);
\draw(9.3,-1.68) ellipse (0.22675781250000002cm and 0.24414893617021272cm);
\draw (10.08,0.38)-- (10.98,-0.34);
\draw (9.3,-1.7)-- (9.86,0.16);
\draw (10.08,0.16)-- (10.58,-1.68);
\draw (8.98,-0.36)-- (9.8,0.38);
\draw(6.38,0.32971746379216116) ellipse (0.22675781250000002cm and 0.24414893617021272cm);
\draw(5.344458247200068,-0.4226476591422038) ellipse (0.22675781250000002cm and 0.24414893617021272cm);
\draw(7.415541752799932,-0.4226476591422043) ellipse (0.22675781250000002cm and 0.24414893617021272cm);
\draw(7.02,-1.64) ellipse (0.22675781250000002cm and 0.24414893617021272cm);
\draw(5.74,-1.64) ellipse (0.22675781250000002cm and 0.24414893617021272cm);
\draw (6.52,0.42)-- (7.42,-0.3);
\draw (5.74,-1.66)-- (6.3,0.2);
\draw (6.52,0.4)-- (7.02,-1.64);
\draw (5.42,-0.32)-- (6.24,0.42);
\draw (4.540000000000002,-0.8199999999999911) node[anchor=north west] {,};
\draw (7.980000000000003,-0.859999999999991) node[anchor=north west] {,};
\begin{scriptsize}
\draw [fill=black] (-4.8,-1.54) circle (1.5pt);
\draw [fill=black] (-3.52,-1.54) circle (1.5pt);
\draw [fill=uuuuuu] (-5.195541752799932,-0.3226476591422037) circle (1.5pt);
\draw [fill=uuuuuu] (-5.195541752799932,-0.3226476591422037) circle (1.5pt);
\draw [fill=black] (-4.8,-1.56) circle (1.5pt);
\draw [fill=black] (-3.18,-0.34) circle (1.5pt);
\draw [fill=black] (-4.18,0.46) circle (1.5pt);
\draw [fill=black] (2.38,-1.66) circle (1.5pt);
\draw [fill=black] (2.38,-1.68) circle (1.5pt);
\draw [fill=black] (4.06,-0.32) circle (1.5pt);
\draw [fill=black] (3.16,0.4) circle (1.5pt);
\draw [fill=black] (2.06,-0.34) circle (1.5pt);
\draw [fill=black] (2.92,0.4) circle (1.5pt);
\draw [fill=black] (3.66,-1.66) circle (1.5pt);
\draw [fill=black] (-1.18,-1.62) circle (1.5pt);
\draw [fill=black] (-1.18,-1.64) circle (1.5pt);
\draw [fill=black] (0.5,-0.28) circle (1.5pt);
\draw [fill=black] (-1.5,-0.3) circle (1.5pt);
\draw [fill=black] (0.1,-1.62) circle (1.5pt);
\draw [fill=black] (-0.68,0.38) circle (1.5pt);
\draw [fill=black] (9.3,-1.68) circle (1.5pt);
\draw [fill=black] (9.3,-1.7) circle (1.5pt);
\draw [fill=black] (10.98,-0.34) circle (1.5pt);
\draw [fill=black] (10.08,0.38) circle (1.5pt);
\draw [fill=black] (8.98,-0.36) circle (1.5pt);
\draw [fill=black] (9.86,0.16) circle (1.5pt);
\draw [fill=black] (10.58,-1.68) circle (1.5pt);
\draw [fill=black] (10.08,0.16) circle (1.5pt);
\draw [fill=black] (9.8,0.38) circle (1.5pt);
\draw [fill=black] (5.74,-1.64) circle (1.5pt);
\draw [fill=black] (5.74,-1.66) circle (1.5pt);
\draw [fill=black] (7.42,-0.3) circle (1.5pt);
\draw [fill=black] (6.52,0.42) circle (1.5pt);
\draw [fill=black] (5.42,-0.32) circle (1.5pt);
\draw [fill=black] (6.3,0.2) circle (1.5pt);
\draw [fill=black] (7.02,-1.64) circle (1.5pt);
\draw [fill=black] (6.24,0.42) circle (1.5pt);
\draw [fill=uuuuuu] (-0.3895159202692209,0.3968055915091898) circle (1.5pt);
\end{scriptsize}
\end{tikzpicture} 
\vspace{4pt}
\definecolor{uuuuuu}{rgb}{0.26666666666666666,0.26666666666666666,0.26666666666666666}
\begin{tikzpicture}[line cap=round,line join=round,>=triangle 45,x=0.7577519379844962cm,y=0.8138297872340424cm]
\clip(-5.600000000000003,-2.8599999999999888) rectangle (15.040000000000006,0.9000000000000069);
\draw(-4.16,0.42971746379216136) ellipse (0.22732558139534886cm and 0.24414893617021272cm);
\draw(-5.195541752799932,-0.3226476591422037) ellipse (0.22732558139534886cm and 0.24414893617021272cm);
\draw(-3.1244582472000677,-0.3226476591422038) ellipse (0.22732558139534886cm and 0.24414893617021272cm);
\draw(-3.52,-1.54) ellipse (0.22732558139534886cm and 0.24414893617021272cm);
\draw(-4.8,-1.54) ellipse (0.22732558139534886cm and 0.24414893617021272cm);
\draw (-4.720000000000003,-1.9799999999999898) node[anchor=north west] {$R_{17}$};
\draw (-1.2200000000000015,-1.9599999999999898) node[anchor=north west] {$S_{22}$};
\draw (2.4,-1.9599999999999898) node[anchor=north west] {$R_{19}$};
\draw [->,line width=1.2000000000000002pt] (-2.86,-0.9) -- (-1.98,-0.9);
\draw (-5.195541752799932,-0.3226476591422037)-- (-4.8,-1.54);
\draw (-5.195541752799932,-0.3226476591422037)-- (-4.18,0.46);
\draw (-3.18,-0.34)-- (-3.52,-1.54);
\draw(-0.5599999999999998,0.4297174637921619) ellipse (0.22732558139534886cm and 0.24414893617021272cm);
\draw(-1.5955417527999325,-0.32264765914220334) ellipse (0.22732558139534886cm and 0.24414893617021272cm);
\draw(0.47554175279993294,-0.32264765914220384) ellipse (0.22732558139534886cm and 0.24414893617021272cm);
\draw(0.08,-1.54) ellipse (0.22732558139534886cm and 0.24414893617021272cm);
\draw(-1.2,-1.54) ellipse (0.22732558139534886cm and 0.24414893617021272cm);
\draw (-0.58,0.46)-- (-1.52,-0.22);
\draw (-1.2,-1.56)-- (-1.64,-0.46);
\draw (0.48,-0.44)-- (0.08,-1.54);
\draw(3.0200000000000005,0.30971746379216236) ellipse (0.22732558139534886cm and 0.24414893617021272cm);
\draw(1.9844582472000674,-0.442647659142203) ellipse (0.22732558139534886cm and 0.24414893617021272cm);
\draw(4.0555417527999325,-0.44264765914220355) ellipse (0.22732558139534886cm and 0.24414893617021272cm);
\draw(3.66,-1.66) ellipse (0.22732558139534886cm and 0.24414893617021272cm);
\draw(2.38,-1.66) ellipse (0.22732558139534886cm and 0.24414893617021272cm);
\draw (3.02,0.34)-- (4.06,-0.32);
\draw (2.38,-1.68)-- (2.04,-0.36);
\draw (2.38,-1.68)-- (3.66,-1.66);
\draw (2.06,-0.34)-- (3.04,0.32);
\draw (9.240000000000002,-2.0599999999999894) node[anchor=north west] {$S_{24}$};
\draw (12.860000000000005,-2.03999999999999) node[anchor=north west] {$S_{25}$};
\draw [->,line width=1.2000000000000002pt] (4.26,-1.0) -- (5.14,-1.0);
\draw [line width=1.2000000000000002pt,dash pattern=on 4pt off 4pt] (1.2,0.68)-- (1.26,-1.84);
\draw (11.440000000000005,-0.6399999999999912) node[anchor=north west] {,};
\draw (5.940000000000001,-1.9999999999999896) node[anchor=north west] {$S_{23}$};
\draw (-5.195541752799932,-0.3226476591422037)-- (-3.18,-0.34);
\draw (-3.52,-1.54)-- (-4.8,-1.56);
\draw (-1.2,-1.56)-- (0.08,-1.54);
\draw (0.48,-0.44)-- (-1.64,-0.46);
\draw (4.06,-0.32)-- (2.06,-0.34);
\draw (7.700000000000002,-0.5599999999999914) node[anchor=north west] {,};
\draw(6.34,0.409717463792162) ellipse (0.22732558139534886cm and 0.24414893617021272cm);
\draw(5.304458247200067,-0.3426476591422036) ellipse (0.22732558139534886cm and 0.24414893617021272cm);
\draw(7.375541752799933,-0.34264765914220346) ellipse (0.22732558139534886cm and 0.24414893617021272cm);
\draw(6.98,-1.56) ellipse (0.22732558139534886cm and 0.24414893617021272cm);
\draw(5.7,-1.56) ellipse (0.22732558139534886cm and 0.24414893617021272cm);
\draw (6.34,0.44)-- (7.38,-0.22);
\draw (5.7,-1.58)-- (6.98,-1.56);
\draw (5.38,-0.24)-- (6.36,0.42);
\draw (7.38,-0.22)-- (5.38,-0.24);
\draw(13.399999999999999,0.3897174637921603) ellipse (0.22732558139534886cm and 0.24414893617021272cm);
\draw(12.364458247200067,-0.36264765914220476) ellipse (0.22732558139534886cm and 0.24414893617021272cm);
\draw(14.43554175279993,-0.3626476591422039) ellipse (0.22732558139534886cm and 0.24414893617021272cm);
\draw(14.04,-1.58) ellipse (0.22732558139534886cm and 0.24414893617021272cm);
\draw(12.76,-1.58) ellipse (0.22732558139534886cm and 0.24414893617021272cm);
\draw (13.4,0.42)-- (14.44,-0.24);
\draw (12.66,-1.6)-- (12.42,-0.28);
\draw (12.44,-0.26)-- (13.42,0.4);
\draw (14.44,-0.24)-- (12.44,-0.26);
\draw (5.7,-1.58)-- (5.28,-0.48);
\draw(9.640000000000002,0.42971746379216147) ellipse (0.22732558139534886cm and 0.24414893617021272cm);
\draw(8.604458247200068,-0.32264765914220306) ellipse (0.22732558139534886cm and 0.24414893617021272cm);
\draw(10.675541752799933,-0.32264765914220495) ellipse (0.22732558139534886cm and 0.24414893617021272cm);
\draw(10.28,-1.54) ellipse (0.22732558139534886cm and 0.24414893617021272cm);
\draw(9.0,-1.54) ellipse (0.22732558139534886cm and 0.24414893617021272cm);
\draw (9.64,0.46)-- (10.68,-0.2);
\draw (8.68,-0.22)-- (9.66,0.44);
\draw (10.68,-0.2)-- (8.68,-0.22);
\draw (10.26,-1.52)-- (9.14,-1.54);
\draw (8.58,-0.46)-- (8.94,-1.52);
\draw (14.04,-1.58)-- (12.86,-1.6);
\begin{scriptsize}
\draw [fill=black] (-4.8,-1.54) circle (1.5pt);
\draw [fill=black] (-3.52,-1.54) circle (1.5pt);
\draw [fill=uuuuuu] (-5.195541752799932,-0.3226476591422037) circle (1.5pt);
\draw [fill=uuuuuu] (-5.195541752799932,-0.3226476591422037) circle (1.5pt);
\draw [fill=black] (-4.8,-1.56) circle (1.5pt);
\draw [fill=black] (-3.18,-0.34) circle (1.5pt);
\draw [fill=black] (-4.18,0.46) circle (1.5pt);
\draw [fill=black] (-1.2,-1.54) circle (1.5pt);
\draw [fill=black] (-1.2,-1.56) circle (1.5pt);
\draw [fill=black] (-0.58,0.46) circle (1.5pt);
\draw [fill=black] (-1.52,-0.22) circle (1.5pt);
\draw [fill=black] (-1.64,-0.46) circle (1.5pt);
\draw [fill=black] (0.08,-1.54) circle (1.5pt);
\draw [fill=black] (0.48,-0.44) circle (1.5pt);
\draw [fill=black] (2.38,-1.66) circle (1.5pt);
\draw [fill=black] (2.38,-1.68) circle (1.5pt);
\draw [fill=black] (4.06,-0.32) circle (1.5pt);
\draw [fill=black] (3.02,0.34) circle (1.5pt);
\draw [fill=black] (2.06,-0.34) circle (1.5pt);
\draw [fill=black] (3.66,-1.66) circle (1.5pt);
\draw [fill=black] (5.7,-1.56) circle (1.5pt);
\draw [fill=black] (5.7,-1.58) circle (1.5pt);
\draw [fill=black] (7.38,-0.22) circle (1.5pt);
\draw [fill=black] (6.34,0.44) circle (1.5pt);
\draw [fill=black] (5.38,-0.24) circle (1.5pt);
\draw [fill=black] (6.98,-1.56) circle (1.5pt);
\draw [fill=black] (12.66,-1.6) circle (1.5pt);
\draw [fill=black] (14.44,-0.24) circle (1.5pt);
\draw [fill=black] (13.4,0.42) circle (1.5pt);
\draw [fill=black] (12.44,-0.26) circle (1.5pt);
\draw [fill=black] (5.28,-0.48) circle (1.5pt);
\draw [fill=black] (10.68,-0.2) circle (1.5pt);
\draw [fill=black] (9.64,0.46) circle (1.5pt);
\draw [fill=black] (8.68,-0.22) circle (1.5pt);
\draw [fill=black] (8.58,-0.46) circle (1.5pt);
\draw [fill=black] (10.26,-1.52) circle (1.5pt);
\draw [fill=black] (9.14,-1.54) circle (1.5pt);
\draw [fill=black] (8.94,-1.52) circle (1.5pt);
\draw [fill=black] (12.86,-1.6) circle (1.5pt);
\draw [fill=black] (14.04,-1.58) circle (1.5pt);
\end{scriptsize}
\end{tikzpicture} 
\vspace{4pt}
\definecolor{uuuuuu}{rgb}{0.26666666666666666,0.26666666666666666,0.26666666666666666}
\begin{tikzpicture}[line cap=round,line join=round,>=triangle 45,x=0.7286821705426357cm,y=0.8138297872340424cm]
\clip(-5.600000000000003,-2.8599999999999888) rectangle (15.040000000000006,0.9000000000000069);
\draw(-4.16,0.42971746379216136) ellipse (0.2186046511627907cm and 0.24414893617021272cm);
\draw(-5.195541752799932,-0.3226476591422037) ellipse (0.2186046511627907cm and 0.24414893617021272cm);
\draw(-3.1244582472000677,-0.3226476591422038) ellipse (0.2186046511627907cm and 0.24414893617021272cm);
\draw(-3.52,-1.54) ellipse (0.2186046511627907cm and 0.24414893617021272cm);
\draw(-4.8,-1.54) ellipse (0.2186046511627907cm and 0.24414893617021272cm);
\draw (-4.720000000000003,-1.9799999999999898) node[anchor=north west] {$R_{20}$};
\draw (-1.2200000000000015,-1.9599999999999898) node[anchor=north west] {$S_{26}$};
\draw (2.4,-1.9599999999999898) node[anchor=north west] {$S_{27}$};
\draw [->,line width=1.2000000000000002pt] (-2.86,-0.9) -- (-1.98,-0.9);
\draw (-5.195541752799932,-0.3226476591422037)-- (-4.8,-1.54);
\draw (-5.195541752799932,-0.3226476591422037)-- (-4.18,0.46);
\draw (-3.18,-0.34)-- (-3.52,-1.54);
\draw(-0.5599999999999998,0.4297174637921619) ellipse (0.2186046511627907cm and 0.24414893617021272cm);
\draw(-1.5955417527999325,-0.32264765914220334) ellipse (0.2186046511627907cm and 0.24414893617021272cm);
\draw(0.47554175279993294,-0.32264765914220384) ellipse (0.2186046511627907cm and 0.24414893617021272cm);
\draw(0.08,-1.54) ellipse (0.2186046511627907cm and 0.24414893617021272cm);
\draw(-1.2,-1.54) ellipse (0.2186046511627907cm and 0.24414893617021272cm);
\draw (-0.58,0.46)-- (-1.62,-0.22);
\draw (-1.2,-1.56)-- (-1.64,-0.46);
\draw(3.0200000000000005,0.30971746379216236) ellipse (0.2186046511627907cm and 0.24414893617021272cm);
\draw(1.9844582472000674,-0.442647659142203) ellipse (0.2186046511627907cm and 0.24414893617021272cm);
\draw(4.0555417527999325,-0.44264765914220355) ellipse (0.2186046511627907cm and 0.24414893617021272cm);
\draw(3.66,-1.66) ellipse (0.2186046511627907cm and 0.24414893617021272cm);
\draw(2.38,-1.66) ellipse (0.2186046511627907cm and 0.24414893617021272cm);
\draw (3.02,0.34)-- (4.06,-0.32);
\draw (2.06,-0.34)-- (3.04,0.32);
\draw (9.240000000000002,-2.0599999999999894) node[anchor=north west] {$S_{28}$};
\draw [->,line width=1.2000000000000002pt] (7.6,-0.9) -- (8.48,-0.9);
\draw [line width=1.2000000000000002pt,dash pattern=on 4pt off 4pt] (4.6,0.58)-- (4.66,-1.94);
\draw (5.940000000000001,-1.9999999999999896) node[anchor=north west] {$R_{23}$};
\draw (-5.195541752799932,-0.3226476591422037)-- (-3.18,-0.34);
\draw (4.06,-0.32)-- (2.06,-0.34);
\draw (0.9999999999999992,-0.6199999999999913) node[anchor=north west] {,};
\draw(6.34,0.409717463792162) ellipse (0.2186046511627907cm and 0.24414893617021272cm);
\draw(5.304458247200067,-0.3426476591422036) ellipse (0.2186046511627907cm and 0.24414893617021272cm);
\draw(7.375541752799933,-0.34264765914220346) ellipse (0.2186046511627907cm and 0.24414893617021272cm);
\draw(6.98,-1.56) ellipse (0.2186046511627907cm and 0.24414893617021272cm);
\draw(5.7,-1.56) ellipse (0.2186046511627907cm and 0.24414893617021272cm);
\draw (6.34,0.44)-- (7.36,-0.34);
\draw (5.7,-1.58)-- (5.32,-0.36);
\draw(9.640000000000002,0.42971746379216147) ellipse (0.2186046511627907cm and 0.24414893617021272cm);
\draw(8.604458247200068,-0.32264765914220306) ellipse (0.2186046511627907cm and 0.24414893617021272cm);
\draw(10.675541752799933,-0.32264765914220495) ellipse (0.2186046511627907cm and 0.24414893617021272cm);
\draw(10.28,-1.54) ellipse (0.2186046511627907cm and 0.24414893617021272cm);
\draw(9.0,-1.54) ellipse (0.2186046511627907cm and 0.24414893617021272cm);
\draw (9.76,0.4)-- (10.64,-0.34);
\draw (-4.18,0.46)-- (-3.18,-0.34);
\draw (0.5,-0.28)-- (-1.62,-0.22);
\draw (-0.58,0.46)-- (0.5,-0.28);
\draw (2.38,-1.66)-- (1.94,-0.56);
\draw (3.66,-1.66)-- (4.04,-0.54);
\draw (0.5,-0.28)-- (0.08,-1.54);
\draw (6.34,0.44)-- (5.32,-0.36);
\draw (6.34,0.44)-- (5.7,-1.58);
\draw (6.34,0.44)-- (7.02,-1.58);
\draw (7.02,-1.58)-- (7.36,-0.34);
\draw (8.62,-0.3)-- (9.5,0.42);
\draw (9.5,0.42)-- (9.02,-1.54);
\draw (9.02,-1.54)-- (8.62,-0.3);
\draw (10.64,-0.34)-- (10.26,-1.52);
\draw (10.26,-1.52)-- (9.76,0.4);
\begin{scriptsize}
\draw [fill=black] (-4.8,-1.54) circle (1.5pt);
\draw [fill=black] (-3.52,-1.54) circle (1.5pt);
\draw [fill=uuuuuu] (-5.195541752799932,-0.3226476591422037) circle (1.5pt);
\draw [fill=uuuuuu] (-5.195541752799932,-0.3226476591422037) circle (1.5pt);
\draw [fill=black] (-4.8,-1.56) circle (1.5pt);
\draw [fill=black] (-3.18,-0.34) circle (1.5pt);
\draw [fill=black] (-4.18,0.46) circle (1.5pt);
\draw [fill=black] (-1.2,-1.54) circle (1.5pt);
\draw [fill=black] (-1.2,-1.56) circle (1.5pt);
\draw [fill=black] (-0.58,0.46) circle (1.5pt);
\draw [fill=black] (-1.62,-0.22) circle (1.5pt);
\draw [fill=black] (-1.64,-0.46) circle (1.5pt);
\draw [fill=black] (0.08,-1.54) circle (1.5pt);
\draw [fill=black] (0.5,-0.28) circle (1.5pt);
\draw [fill=black] (2.38,-1.66) circle (1.5pt);
\draw [fill=black] (2.38,-1.68) circle (1.5pt);
\draw [fill=black] (4.06,-0.32) circle (1.5pt);
\draw [fill=black] (3.02,0.34) circle (1.5pt);
\draw [fill=black] (2.06,-0.34) circle (1.5pt);
\draw [fill=black] (3.66,-1.66) circle (1.5pt);
\draw [fill=black] (5.7,-1.56) circle (1.5pt);
\draw [fill=black] (5.7,-1.58) circle (1.5pt);
\draw [fill=black] (7.36,-0.34) circle (1.5pt);
\draw [fill=black] (6.34,0.44) circle (1.5pt);
\draw [fill=black] (7.02,-1.58) circle (1.5pt);
\draw [fill=black] (5.32,-0.36) circle (1.5pt);
\draw [fill=black] (10.64,-0.34) circle (1.5pt);
\draw [fill=black] (9.76,0.4) circle (1.5pt);
\draw [fill=black] (8.62,-0.3) circle (1.5pt);
\draw [fill=black] (10.26,-1.52) circle (1.5pt);
\draw [fill=black] (9.02,-1.54) circle (1.5pt);
\draw [fill=black] (1.94,-0.56) circle (1.5pt);
\draw [fill=black] (4.04,-0.54) circle (1.5pt);
\draw [fill=black] (9.5,0.42) circle (1.5pt);
\end{scriptsize}
\end{tikzpicture} 
\vspace{4pt}
\begin{tikzpicture}[line cap=round,line join=round,>=triangle 45,x=0.7286821705426368cm,y=0.8138297872340423cm]
\clip(1.5400000000000056,-2.9199999999999884) rectangle (22.180000000000092,0.8400000000000071);
\draw (2.4000000000000092,-1.9599999999999898) node[anchor=north west] {$R_{21}$};
\draw(3.0200000000000005,0.30971746379216236) ellipse (0.21860465116279104cm and 0.2441489361702127cm);
\draw(1.9844582472000674,-0.442647659142203) ellipse (0.21860465116279104cm and 0.2441489361702127cm);
\draw(4.0555417527999325,-0.44264765914220355) ellipse (0.21860465116279104cm and 0.2441489361702127cm);
\draw(3.66,-1.66) ellipse (0.21860465116279104cm and 0.2441489361702127cm);
\draw(2.38,-1.66) ellipse (0.21860465116279104cm and 0.2441489361702127cm);
\draw (3.02,0.34)-- (4.06,-0.32);
\draw (2.38,-1.68)-- (2.04,-0.36);
\draw (2.06,-0.34)-- (3.04,0.32);
\draw (9.240000000000038,-2.0599999999999894) node[anchor=north west] {$S_{30}$};
\draw [->,line width=1.2000000000000002pt] (4.26,-1.0) -- (5.14,-1.0);
\draw (5.9400000000000235,-1.9999999999999896) node[anchor=north west] {$S_{29}$};
\draw (4.06,-0.32)-- (2.06,-0.34);
\draw (7.720000000000031,-0.5599999999999913) node[anchor=north west] {,};
\draw(6.34,0.409717463792162) ellipse (0.21860465116279104cm and 0.2441489361702127cm);
\draw(5.304458247200067,-0.3426476591422036) ellipse (0.21860465116279104cm and 0.2441489361702127cm);
\draw(7.375541752799933,-0.34264765914220346) ellipse (0.21860465116279104cm and 0.2441489361702127cm);
\draw(6.98,-1.56) ellipse (0.21860465116279104cm and 0.2441489361702127cm);
\draw(5.7,-1.56) ellipse (0.21860465116279104cm and 0.2441489361702127cm);
\draw (6.34,0.44)-- (7.38,-0.22);
\draw (5.38,-0.24)-- (6.36,0.42);
\draw (7.38,-0.22)-- (5.38,-0.24);
\draw (5.7,-1.58)-- (5.28,-0.48);
\draw(9.640000000000002,0.42971746379216147) ellipse (0.21860465116279104cm and 0.2441489361702127cm);
\draw(8.604458247200068,-0.32264765914220306) ellipse (0.21860465116279104cm and 0.2441489361702127cm);
\draw(10.675541752799933,-0.32264765914220495) ellipse (0.21860465116279104cm and 0.2441489361702127cm);
\draw(10.28,-1.54) ellipse (0.21860465116279104cm and 0.2441489361702127cm);
\draw(9.0,-1.54) ellipse (0.21860465116279104cm and 0.2441489361702127cm);
\draw (9.64,0.46)-- (10.68,-0.2);
\draw (8.68,-0.22)-- (9.66,0.44);
\draw (10.68,-0.2)-- (8.68,-0.22);
\draw (8.48,-0.42)-- (8.94,-1.52);
\draw (3.66,-1.66)-- (2.06,-0.34);
\draw (6.98,-1.56)-- (5.38,-0.24);
\draw (10.26,-1.52)-- (8.72,-0.42);
\draw (12.400000000000052,-2.01999999999999) node[anchor=north west] {$S_{31}$};
\draw(12.800000000000002,0.4697174637921615) ellipse (0.21860465116279104cm and 0.2441489361702127cm);
\draw(11.764458247200068,-0.282647659142203) ellipse (0.21860465116279104cm and 0.2441489361702127cm);
\draw(13.835541752799934,-0.2826476591422049) ellipse (0.21860465116279104cm and 0.2441489361702127cm);
\draw(13.44,-1.5) ellipse (0.21860465116279104cm and 0.2441489361702127cm);
\draw(12.16,-1.5) ellipse (0.21860465116279104cm and 0.2441489361702127cm);
\draw (12.8,0.5)-- (13.84,-0.16);
\draw (11.84,-0.18)-- (12.82,0.48);
\draw (13.84,-0.16)-- (11.84,-0.18);
\draw (11.82,-0.4)-- (12.1,-1.48);
\draw (13.42,-1.48)-- (11.82,-0.42);
\begin{scriptsize}
\draw [fill=black] (2.38,-1.66) circle (1.5pt);
\draw [fill=black] (2.38,-1.68) circle (1.5pt);
\draw [fill=black] (4.06,-0.32) circle (1.5pt);
\draw [fill=black] (3.02,0.34) circle (1.5pt);
\draw [fill=black] (2.06,-0.34) circle (1.5pt);
\draw [fill=black] (3.66,-1.66) circle (1.5pt);
\draw [fill=black] (5.7,-1.56) circle (1.5pt);
\draw [fill=black] (5.7,-1.58) circle (1.5pt);
\draw [fill=black] (7.38,-0.22) circle (1.5pt);
\draw [fill=black] (6.34,0.44) circle (1.5pt);
\draw [fill=black] (5.38,-0.24) circle (1.5pt);
\draw [fill=black] (6.98,-1.56) circle (1.5pt);
\draw [fill=black] (5.28,-0.48) circle (1.5pt);
\draw [fill=black] (10.68,-0.2) circle (1.5pt);
\draw [fill=black] (9.64,0.46) circle (1.5pt);
\draw [fill=black] (8.68,-0.22) circle (1.5pt);
\draw [fill=black] (8.48,-0.42) circle (1.5pt);
\draw [fill=black] (10.26,-1.52) circle (1.5pt);
\draw [fill=black] (8.94,-1.52) circle (1.5pt);
\draw [fill=black] (8.72,-0.42) circle (1.5pt);
\draw [fill=black] (13.84,-0.16) circle (1.5pt);
\draw [fill=black] (12.8,0.5) circle (1.5pt);
\draw [fill=black] (11.84,-0.18) circle (1.5pt);
\draw [fill=black] (11.82,-0.4) circle (1.5pt);
\draw [fill=black] (13.42,-1.48) circle (1.5pt);
\draw [fill=black] (12.1,-1.48) circle (1.5pt);
\end{scriptsize}
\end{tikzpicture} 
\vspace{4pt}
\definecolor{uuuuuu}{rgb}{0.26666666666666666,0.26666666666666666,0.26666666666666666}
\begin{tikzpicture}[line cap=round,line join=round,>=triangle 45,x=0.7577519379844962cm,y=0.8138297872340424cm]
\clip(-5.600000000000003,-2.8599999999999888) rectangle (15.040000000000006,0.9000000000000069);
\draw(-4.16,0.42971746379216136) ellipse (0.22732558139534886cm and 0.24414893617021272cm);
\draw(-5.195541752799932,-0.3226476591422037) ellipse (0.22732558139534886cm and 0.24414893617021272cm);
\draw(-3.1244582472000677,-0.3226476591422038) ellipse (0.22732558139534886cm and 0.24414893617021272cm);
\draw(-3.52,-1.54) ellipse (0.22732558139534886cm and 0.24414893617021272cm);
\draw(-4.8,-1.54) ellipse (0.22732558139534886cm and 0.24414893617021272cm);
\draw (-4.720000000000003,-1.9799999999999898) node[anchor=north west] {$R_{24}$};
\draw (-1.2200000000000015,-1.9599999999999898) node[anchor=north west] {$S_{32}$};
\draw [->,line width=1.2000000000000002pt] (-2.86,-0.9) -- (-1.98,-0.9);
\draw (-5.195541752799932,-0.3226476591422037)-- (-4.8,-1.54);
\draw (-3.18,-0.34)-- (-3.52,-1.54);
\draw(-0.5599999999999998,0.4297174637921619) ellipse (0.22732558139534886cm and 0.24414893617021272cm);
\draw(-1.5955417527999325,-0.32264765914220334) ellipse (0.22732558139534886cm and 0.24414893617021272cm);
\draw(0.47554175279993294,-0.32264765914220384) ellipse (0.22732558139534886cm and 0.24414893617021272cm);
\draw(0.08,-1.54) ellipse (0.22732558139534886cm and 0.24414893617021272cm);
\draw(-1.2,-1.54) ellipse (0.22732558139534886cm and 0.24414893617021272cm);
\draw (-0.58,0.46)-- (0.48,-0.2);
\draw (-1.2,-1.56)-- (-1.64,-0.46);
\draw (0.48,-0.44)-- (0.08,-1.54);
\draw [line width=1.2000000000000002pt,dash pattern=on 4pt off 4pt] (1.2,0.68)-- (1.26,-1.84);
\draw (-5.195541752799932,-0.3226476591422037)-- (-3.18,-0.34);
\draw (-3.52,-1.54)-- (-4.8,-1.56);
\draw (-1.2,-1.56)-- (0.08,-1.54);
\draw (0.48,-0.44)-- (-1.64,-0.46);
\draw (-4.18,0.46)-- (-3.18,-0.34);
\draw (-3.18,-0.34)-- (-4.8,-1.56);
\draw (0.48,-0.44)-- (-1.2,-1.56);
\draw(2.8600000000000008,0.36971746379216175) ellipse (0.22732558139534886cm and 0.24414893617021272cm);
\draw(1.8244582472000679,-0.38264765914220333) ellipse (0.22732558139534886cm and 0.24414893617021272cm);
\draw(3.895541752799933,-0.38264765914220406) ellipse (0.22732558139534886cm and 0.24414893617021272cm);
\draw(3.5,-1.6) ellipse (0.22732558139534886cm and 0.24414893617021272cm);
\draw(2.22,-1.6) ellipse (0.22732558139534886cm and 0.24414893617021272cm);
\draw (2.3,-2.03999999999999) node[anchor=north west] {$R_{25}$};
\draw (5.800000000000001,-2.01999999999999) node[anchor=north west] {$S_{33}$};
\draw [->,line width=1.2000000000000002pt] (4.16,-0.96) -- (5.04,-0.96);
\draw (1.8244582472000679,-0.38264765914220333)-- (2.22,-1.6);
\draw (3.84,-0.4)-- (3.5,-1.6);
\draw(6.46,0.36971746379216097) ellipse (0.22732558139534886cm and 0.24414893617021272cm);
\draw(5.424458247200068,-0.382647659142204) ellipse (0.22732558139534886cm and 0.24414893617021272cm);
\draw(7.495541752799932,-0.3826476591422045) ellipse (0.22732558139534886cm and 0.24414893617021272cm);
\draw(7.1,-1.6) ellipse (0.22732558139534886cm and 0.24414893617021272cm);
\draw(5.82,-1.6) ellipse (0.22732558139534886cm and 0.24414893617021272cm);
\draw (6.44,0.4)-- (5.44,-0.32);
\draw (5.82,-1.62)-- (5.38,-0.52);
\draw (7.5,-0.5)-- (7.1,-1.6);
\draw [line width=1.2000000000000002pt,dash pattern=on 4pt off 4pt] (8.22,0.62)-- (8.28,-1.9);
\draw (1.8244582472000679,-0.38264765914220333)-- (3.84,-0.4);
\draw (3.5,-1.6)-- (2.22,-1.62);
\draw (5.82,-1.62)-- (7.1,-1.6);
\draw (7.5,-0.5)-- (5.38,-0.52);
\draw (3.84,-0.4)-- (2.22,-1.62);
\draw (7.5,-0.5)-- (5.82,-1.62);
\draw(9.899999999999999,0.3297174637921605) ellipse (0.22732558139534886cm and 0.24414893617021272cm);
\draw(8.864458247200067,-0.4226476591422046) ellipse (0.22732558139534886cm and 0.24414893617021272cm);
\draw(10.93554175279993,-0.42264765914220376) ellipse (0.22732558139534886cm and 0.24414893617021272cm);
\draw(10.54,-1.64) ellipse (0.22732558139534886cm and 0.24414893617021272cm);
\draw(9.26,-1.64) ellipse (0.22732558139534886cm and 0.24414893617021272cm);
\draw (9.340000000000002,-2.0799999999999894) node[anchor=north west] {$R_{30}$};
\draw (12.840000000000005,-2.0599999999999894) node[anchor=north west] {$S_{34}$};
\draw [->,line width=1.2000000000000002pt] (11.2,-1.0) -- (12.08,-1.0);
\draw (8.864458247200067,-0.4226476591422046)-- (9.26,-1.64);
\draw (10.88,-0.44)-- (10.54,-1.64);
\draw(13.5,0.3297174637921637) ellipse (0.22732558139534886cm and 0.24414893617021272cm);
\draw(12.464458247200067,-0.4226476591422021) ellipse (0.22732558139534886cm and 0.24414893617021272cm);
\draw(14.535541752799933,-0.4226476591422026) ellipse (0.22732558139534886cm and 0.24414893617021272cm);
\draw(14.14,-1.64) ellipse (0.22732558139534886cm and 0.24414893617021272cm);
\draw(12.86,-1.64) ellipse (0.22732558139534886cm and 0.24414893617021272cm);
\draw (13.48,0.36)-- (12.48,-0.34);
\draw (12.86,-1.66)-- (12.42,-0.56);
\draw (14.5,-0.5)-- (14.14,-1.64);
\draw [line width=1.2000000000000002pt,dash pattern=on 4pt off 4pt] (15.26,0.58)-- (15.32,-1.94);
\draw (8.864458247200067,-0.4226476591422046)-- (10.88,-0.44);
\draw (10.54,-1.64)-- (9.26,-1.66);
\draw (12.86,-1.66)-- (14.14,-1.64);
\draw (14.5,-0.5)-- (12.42,-0.56);
\draw (10.88,-0.44)-- (9.26,-1.66);
\draw (14.5,-0.5)-- (12.86,-1.66);
\draw (2.84,0.4)-- (1.8244582472000679,-0.38264765914220333);
\draw (14.14,-1.64)-- (12.42,-0.56);
\draw (10.54,-1.64)-- (8.864458247200067,-0.4226476591422046);
\draw (9.88,0.36)-- (8.864458247200067,-0.4226476591422046);
\begin{scriptsize}
\draw [fill=black] (-4.8,-1.54) circle (1.5pt);
\draw [fill=black] (-3.52,-1.54) circle (1.5pt);
\draw [fill=uuuuuu] (-5.195541752799932,-0.3226476591422037) circle (1.5pt);
\draw [fill=uuuuuu] (-5.195541752799932,-0.3226476591422037) circle (1.5pt);
\draw [fill=black] (-4.8,-1.56) circle (1.5pt);
\draw [fill=black] (-3.18,-0.34) circle (1.5pt);
\draw [fill=black] (-4.18,0.46) circle (1.5pt);
\draw [fill=black] (-1.2,-1.54) circle (1.5pt);
\draw [fill=black] (-1.2,-1.56) circle (1.5pt);
\draw [fill=black] (-0.58,0.46) circle (1.5pt);
\draw [fill=black] (0.48,-0.2) circle (1.5pt);
\draw [fill=black] (-1.64,-0.46) circle (1.5pt);
\draw [fill=black] (0.08,-1.54) circle (1.5pt);
\draw [fill=black] (0.48,-0.44) circle (1.5pt);
\draw [fill=black] (2.22,-1.6) circle (1.5pt);
\draw [fill=black] (3.5,-1.6) circle (1.5pt);
\draw [fill=uuuuuu] (1.8244582472000679,-0.38264765914220333) circle (1.5pt);
\draw [fill=uuuuuu] (1.8244582472000679,-0.38264765914220333) circle (1.5pt);
\draw [fill=black] (2.22,-1.62) circle (1.5pt);
\draw [fill=black] (3.84,-0.4) circle (1.5pt);
\draw [fill=black] (2.84,0.4) circle (1.5pt);
\draw [fill=black] (5.82,-1.6) circle (1.5pt);
\draw [fill=black] (5.82,-1.62) circle (1.5pt);
\draw [fill=black] (6.44,0.4) circle (1.5pt);
\draw [fill=black] (5.44,-0.32) circle (1.5pt);
\draw [fill=black] (5.38,-0.52) circle (1.5pt);
\draw [fill=black] (7.1,-1.6) circle (1.5pt);
\draw [fill=black] (7.5,-0.5) circle (1.5pt);
\draw [fill=black] (10.54,-1.64) circle (1.5pt);
\draw [fill=uuuuuu] (8.864458247200067,-0.4226476591422046) circle (1.5pt);
\draw [fill=uuuuuu] (8.864458247200067,-0.4226476591422046) circle (1.5pt);
\draw [fill=black] (9.26,-1.66) circle (1.5pt);
\draw [fill=black] (10.88,-0.44) circle (1.5pt);
\draw [fill=black] (9.88,0.36) circle (1.5pt);
\draw [fill=black] (12.86,-1.64) circle (1.5pt);
\draw [fill=black] (12.86,-1.66) circle (1.5pt);
\draw [fill=black] (13.48,0.36) circle (1.5pt);
\draw [fill=black] (12.48,-0.34) circle (1.5pt);
\draw [fill=black] (12.42,-0.56) circle (1.5pt);
\draw [fill=black] (14.14,-1.64) circle (1.5pt);
\draw [fill=black] (14.5,-0.5) circle (1.5pt);
\end{scriptsize}
\end{tikzpicture} 
\vspace{12pt}

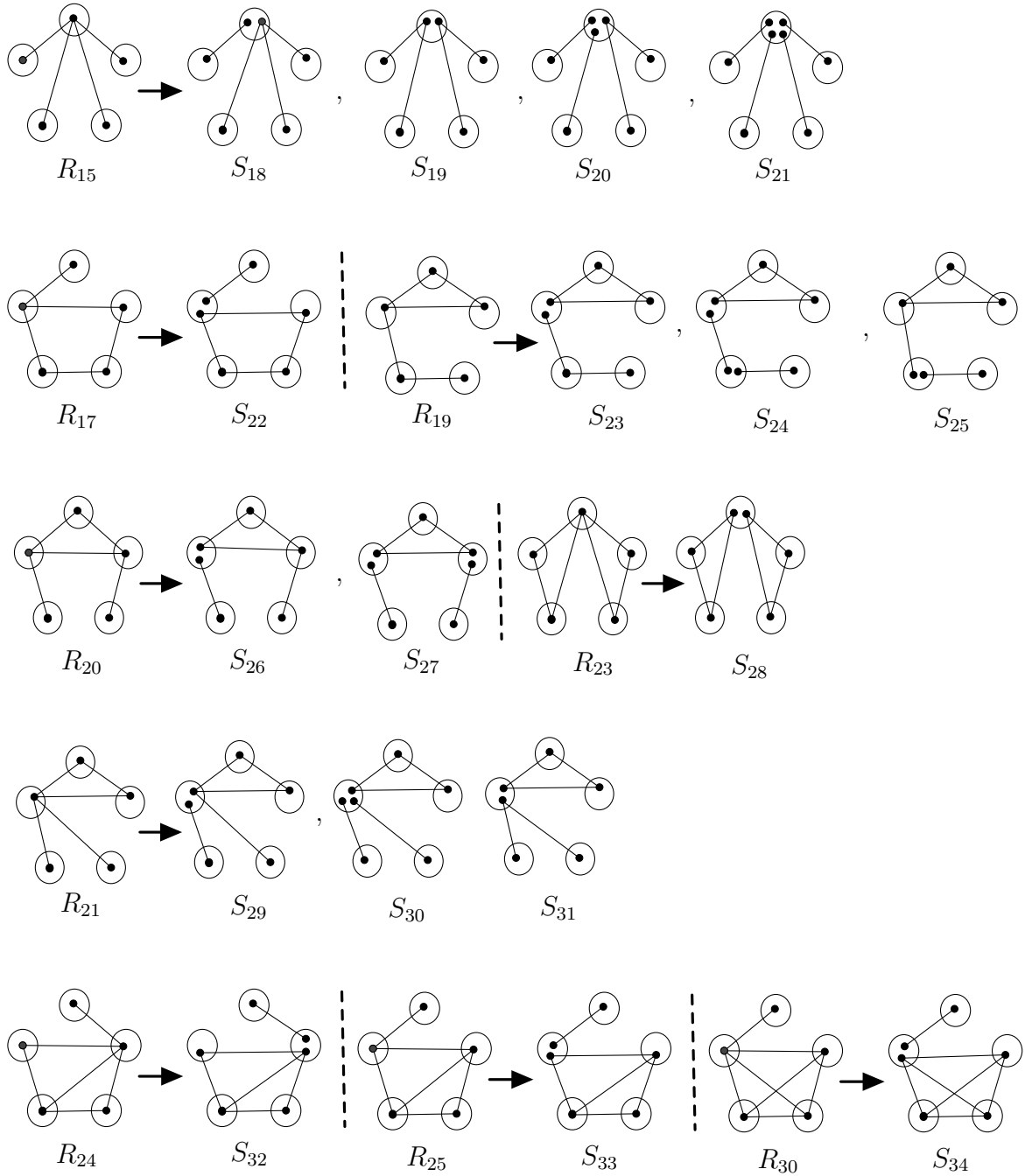
Figure 4. Ramsey $(C_n,K_6)$ critical graphs of Type2 ($S_i$,  $1 \leq i \leq 34$)
\end{center}
\vspace{8pt}

\noindent Henceforth, we conclude that there are exactly 68  Ramsey $(C_n,K_6)$ critical graphs out of which 34 are  categorized as Type1 critical graphs (labeled $R_i$, $1 \leq i \leq 34$ ) and the balance 34 are categorized as Type2 critical graphs (labeled $S_i$, $1 \leq i \leq 34$).

\end{document}